\documentclass{informs1-arxiv}              

\usepackage{natbib}
 \NatBibNumeric
 \def\bibfont{\small}%
 \bibpunct[, ]{[}{]}{,}{n}{}{,}%

\usepackage[hidelinks=true,bookmarksopen=false,draft=false]{hyperref}

\usepackage[utf8]{inputenc}
\usepackage[T1]{fontenc}
\usepackage[english]{babel}
\usepackage{enumitem}
\usepackage{xcolor}

\usepackage{xfrac}
\usepackage{cancel}
\usepackage{bbm}

\usepackage{subcaption}
\usepackage{multirow}
\usepackage{cleveref}

\def\EMAIL#1{\href{mailto:#1}{#1}}
\def\URL#1{\href{#1}{#1}}         

\TheoremsNumberedThrough     

\EquationsNumberedThrough    

\newcommand{\ds}{\displaystyle}
\newcommand{\defn}{\stackrel{\triangle}{=}}
\newcommand{\mymod}{\hspace{-0.2cm}\mod}

\newcommand{\tth}{{}^{\text{th}}}
\newcommand{\amax}{A_{\max}}
\newcommand{\smax}{S_{\max}}
\newcommand{\dto}{\downarrow}

\newcommand{\bI}{\mathbb{I}}

\newcommand{\bR}{\mathbb{R}}
\newcommand{\bZ}{\mathbb{Z}}

\newcommand{\cA}{\mathcal{A}}
\newcommand{\cC}{\mathcal{C}}
\newcommand{\cF}{\mathcal{F}}
\newcommand{\cH}{\mathcal{H}}
\newcommand{\cK}{\mathcal{K}}
\newcommand{\cL}{\mathcal{L}}
\newcommand{\cM}{\mathcal{M}}
\newcommand{\cN}{\mathcal{N}}
\newcommand{\cS}{\mathcal{S}}
\newcommand{\cT}{\mathcal{T}}
\newcommand{\cX}{\mathcal{X}}
\newcommand{\cY}{\mathcal{Y}}

\newcommand{\va}{\boldsymbol{a}}
\newcommand{\vb}{\boldsymbol{b}}
\newcommand{\vc}{\boldsymbol{c}}
\newcommand{\ve}{\boldsymbol{e}}
\newcommand{\vq}{\boldsymbol{q}}
\newcommand{\vr}{\boldsymbol{r}}
\renewcommand{\vs}{\boldsymbol{s}}
\newcommand{\vu}{\boldsymbol{u}}
\newcommand{\vv}{\boldsymbol{v}}
\newcommand{\vw}{\boldsymbol{w}}
\newcommand{\vx}{\boldsymbol{x}}
\newcommand{\vy}{\boldsymbol{y}}
\newcommand{\vz}{\boldsymbol{z}}

\newcommand{\vone}{\boldsymbol{1}}
\newcommand{\vzero}{\boldsymbol{0}}

\newcommand{\valpha}{\boldsymbol{\alpha}}
\newcommand{\vsigma}{\boldsymbol{\sigma}}

\newcommand{\vnu}{\boldsymbol{\nu}}
\newcommand{\vmu}{\boldsymbol{\mu}}
\newcommand{\vpsi}{\boldsymbol{\psi}}
\newcommand{\vlambda}{\boldsymbol{\lambda}}

\newcommand{\qbar}{\overline{q}}
\newcommand{\vabar}{\overline{\va}}
\newcommand{\abar}{\overline{a}}
\newcommand{\vsqm}{\vs(\vqbar^\peps,\Mbar)}
\newcommand{\vuqma}{\vu(\vqbar^\peps,\Mbar,\vabar^\peps)}
\newcommand{\vsbar}{\overline{\vs}}
\newcommand{\sbar}{\overline{s}}
\newcommand{\vubar}{\overline{\vu}}
\newcommand{\ubar}{\overline{u}}
\newcommand{\Mbar}{\overline{M}}
\newcommand{\Bbar}{\overline{B}}

\newcommand{\vqbarpar}{\vqbar_{\parallel}}
\newcommand{\vqbarperp}{\vqbar_{\perp}}
\newcommand{\qbarpari}[1]{\qbar_{\parallel #1}}
\newcommand{\abarpari}[1]{\abar_{\parallel #1}}
\newcommand{\sbarpari}[1]{\sbar_{\parallel #1}}
\newcommand{\ubarpari}[1]{\ubar_{\parallel #1}}

\newcommand{\vqpark}{\vq_{\parallel \cK}}
\newcommand{\vqperpk}{\vq_{\perp \cK}}
\newcommand{\vqbar}{\overline{\vq}}
\newcommand{\vqbarpark}{\vqbar_{\parallel \cK}}
\newcommand{\vqbarperpk}{\vqbar_{\perp \cK}}

\newcommand{\vqparh}{\vq_{\parallel \cH}}
\newcommand{\vqperph}{\vq_{\perp \cH}}
\newcommand{\vqbarparh}{\vqbar_{\parallel \cH}}
\newcommand{\vqbarperph}{\vqbar_{\perp \cH}}

\newcommand{\vabarparh}{\vabar_{\parallel \cH}}
\newcommand{\vsbarparh}{\vsbar_{\parallel \cH}}
\newcommand{\vubarparh}{\vubar_{\parallel \cH}}

\renewcommand{\pm}{{(m)}}
\newcommand{\pl}{{(\ell)}}
\newcommand{\pml}{{(m,\ell)}}
\newcommand{\peps}{{(\epsilon)}}

\newcommand{\vcl}{\vc^\pl}
\newcommand{\bl}{b^\pl}
\newcommand{\bml}{b^\pml}
\newcommand{\piml}{\pi^\pml}
\newcommand{\gammam}{\gamma^\pm}

\newcommand{\vclone}{\vc^{(\ell_1)}}
\newcommand{\vcltwo}{\vc^{(\ell_2)}}

\newcommand{\vcltilde}{\widetilde{\vc}^\pl}
\newcommand{\vqbartildel}{\widetilde{\vqbar}^\pl}
\newcommand{\vubartildel}{\widetilde{\vubar}^\pl}
\newcommand{\vqbarparhtilde}{\vqbartildel_{\parallel\cH}}
\newcommand{\vqbarperphtilde}{\vqbartildel_{\perp\cH}}
\newcommand{\vubarparhtilde}{\vubartildel_{\parallel\cH}}
\newcommand{\vubarperphtilde}{\vubartildel_{\perp\cH}}
\newcommand{\ctildel}{\tilde{c}^\pl}

\newcommand{\ind}[1]{\mathbbm{1}_{\left\{#1 \right\}}}
\newcommand{\Prob}[1]{P\left[#1\right]}
\newcommand{\E}[1]{\mathbb{E}\left[#1 \right]}
\newcommand{\Evq}[1]{\mathbb{E}_{\vq}\left[#1 \right]}
\newcommand{\Em}[1]{\mathbb{E}_m\left[#1 \right]}
\newcommand{\Var}[1]{\text{Var}\left[#1 \right]}
\newcommand{\Cov}[1]{\text{Cov}\left[#1 \right]}

\newcommand{\Phibar}{\overline{\Phi}}
\newcommand{\alphabar}{\overline{\alpha}}
\newcommand{\betabar}{\overline{\beta}}
\newcommand{\chibar}{\overline{\chi}}

\makeatletter
\newcommand{\pushright}[1]{\ifmeasuring@#1\else\omit\hfill$\displaystyle#1$\fi\ignorespaces}

\newcommand{\numberthis}{\addtocounter{equation}{1}\tag{\theequation}}

\def\ba#1\ea{\begin{align*}#1\end{align*}}
\def\ban#1\ean{\begin{align}#1\end{align}}

\allowdisplaybreaks

\begin{document}


\RUNAUTHOR{Hurtado-Lange and Maguluri}

\RUNTITLE{Generalized Switch under Multidimensional State Space Collapse}

\TITLE{Heavy-traffic Analysis of Queueing Systems with no Complete
	Resource Pooling}

\ARTICLEAUTHORS{%
\AUTHOR{Daniela Hurtado-Lange}
\AFF{Georgia Institute of Technology, \EMAIL{d.hurtado@gatech.edu}, \URL{https://sites.google.com/view/daniela-hurtado-lange}}
\AUTHOR{Siva Theja Maguluri}
\AFF{Georgia Institute of Technology, \EMAIL{siva.theja@gatech.edu}, \URL{https://sites.google.com/site/sivatheja/}}
} 

\ABSTRACT{%
We study the heavy-traffic limit of the generalized switch operating under MaxWeight, without assuming that the CRP condition is satisfied and allowing for correlated arrivals. The main contribution of this paper is the steady-state mean of linear combinations of queue lengths in heavy traffic. We showcase the generality of our result by presenting  various stochastic networks as corollaries, each of which is a contribution by itself. In particular, we study the input-queued switch with correlated arrivals and we show that if the state space collapses to a full-dimensional subspace, the correlation among the arrival processes does not matter in heavy traffic. We exemplify this last case with a parallel-server system, an $\cN$-system, and an ad hoc wireless network.

While the above results are obtained using the drift method, we additionally present a negative result showing a limitation of the drift method. We show that it is not possible to obtain the individual queue lengths using the drift method with polynomial test functions.  We do this by presenting an alternate view of the drift method in terms of a system of linear equations, and we use this system of equations to obtain bounds on arbitrary linear combinations of the queue lengths.
}%


\MSCCLASS{60K25
	; 68M20
	; 90B22
	; 60H99
}
\ORMSCLASS{Primary: Probability/Statistics ; secondary: Stochastic processes}
\HISTORY{}

\maketitle

%

\section{Introduction.}\label{sec:introduction}
Resource allocation problems arise in a variety of settings such as wireless networks, wired networks, data centers, cloud computing, ride-hailing systems, call centers, etc. One way to analyze them is to study the delay, using tools from queueing theory. In such case they are modeled as Stochastic Processing Networks (SPNs).
A major challenge is that the analysis of such queueing models is usually not tractable in general settings, and so asymptotic analysis is a popular methodology. Heavy-traffic analysis is an asymptotic approach where the system is loaded very close to its capacity, and the corresponding queueing (delay) behavior of various resource allocation algorithms is studied.

Heavy-traffic limits have been obtained in a wide variety of systems using a program
based on fluid limits, diffusion limits, and Reflected Brownian Motion (RBM) processes, as shown by \citeauthor{harrison_2013_book} \cite{harrison_2013_book}. In this approach, the queueing process is scaled appropriately, and the limiting fluid or diffusion process is studied. The limit of a diffusion scaled process is shown to converge to an RBM process. Typically this RBM lives in a lower-dimensional subspace. This phenomenon is known as State Space Collapse (SSC), and it makes the heavy-traffic analysis tractable since one can study a lower-dimensional RBM. Several systems where the state space collapses to a line (i.e., to a one-dimensional subspace) have been extensively studied in the literature using this approach (see the survey by \citeauthor{williams_survey_SPN} \cite{williams_survey_SPN} for a rigorous list). Typically this happens when there is a unique outer normal vector to the point of the boundary of the capacity region that is being approached in heavy traffic. Such systems are said to satisfy the Complete Resource Pooling (CRP) condition. For a formal definition of the CRP condition, the reader is referred to \citeauthor{stolyar2004maxweight}'s work \cite{stolyar2004maxweight}. As shown by \citeauthor{harlop_state_space} \cite{harlop_state_space} and \citeauthor{dai2008max_pressure} \cite{dai2008max_pressure}, the intuitive meaning of the CRP is that , in the heavy-traffic limit, there is a single bottleneck resource and, hence, the queueing system behaves as a single server queue. Under the CRP condition, in the diffusion limit one obtains an RBM on a line, which is well understood. However, a major challenge is in using this program for SPNs where the CRP condition is not satisfied (i.e., when there are multiple resources that are simultaneously in heavy traffic). In such case, one needs to solve for the steady-state distribution of a RBM in a multidimensional subset of $\bR^n$, and this is not known in general, as shown by \citeauthor{kang2012diffusion} \cite{kang2012diffusion}. The focus of this paper is to study systems that do not satisfy the CRP condition.

According to \citeauthor{kang2012diffusion} \cite{kang2012diffusion} and \citeauthor{shah_switch_open} \cite{shah_switch_open}, one of the simplest queueing systems where the CRP condition is not satisfied is an input-queued switch, and \citeauthor{williams_survey_SPN} \cite{williams_survey_SPN} identifies it as a focus of study in the SPN literature since it serves as a guiding example to study more general systems that do not satisfy CRP.
Recently, the drift method was developed by \citeauthor{atilla} \cite{atilla}, as an alternate way to study heavy-traffic limits of queueing systems, based on a generalization of Kingman's bound \citep{kingman} in a $G/G/1$ queue. The drift method was used to characterize the heavy-traffic scaled sum of queue lengths in input-queued switches by \citeauthor{MagSri_SSY16_Switch} \cite{MagSri_SSY16_Switch}, and \citeauthor{QUESTA_switch} \cite{QUESTA_switch}.

In order to study SPNs when the CRP condition is not met, in this paper we consider a very general queueing model, called \emph{generalized switch}, that subsumes several (single-hop) SPNs with control in the service process, and that was first proposed by \citeauthor{stolyar2004maxweight} \cite{stolyar2004maxweight}. A detailed description of the model is provided in \Cref{sec:gen.switch.model}. Particular cases of the generalized switch are ad hoc wireless networks, wireless networks in presence of fading, input-queued switches, and parallel-server systems.

In this paper we study the generalized switch operating under MaxWeight scheduling algorithm, which we describe in detail in \Cref{sec:gen.switch.model}. MaxWeight algorithm was first proposed by \citeauthor{TasEph_92} \cite{TasEph_92}, in the context of down-link in wireless base stations, and has been used in a variety of queueing systems. For example, in the work by \citeauthor{stolyar2004maxweight} \cite{stolyar2004maxweight}, \citeauthor{gupta2010delay} \cite{gupta2010delay} and \citeauthor{Mey_08}  \cite{Mey_08}. Some of its advantages are that it is a throughput optimal algorithm (i.e., it keeps the system stable for all arrival rates in the capacity region), and it only requires information about the state of the system (and not parameters such as the arrival rates).

The generalized switch has been studied under the CRP condition and independent arrivals assumption 
using both, the diffusion limits approach by \citeauthor{stolyar2004maxweight} \cite{stolyar2004maxweight} and the drift method by \citeauthor{atilla} \cite{atilla}. In this paper, we focus on the case when the CRP condition is not necessarily met, and so SSC may occur to a multidimensional subspace. Also, we assume that the arrival process to each queue is a sequence of independent and identically distributed (i.i.d.) random variables, but we do not require that these sequences are independent of each other. The main contributions of this paper are:
\begin{enumerate}[label=(\roman*)]
	\item In \Cref{gs.thm:bounds} we characterize the heavy-traffic scaled mean of certain linear combinations of the queue lengths in steady state under the MaxWeight algorithm. Moreover, we obtain lower and upper bounds that are valid in all regimes (not necessarily heavy traffic), but are tight in the heavy-traffic regime. This result is immediately applicable in several systems, as we showcase in \Cref{sec:interpretation}, and it includes both, the CRP and the non-CRP cases. Little is known about SPNs that do not satisfy CRP, since the most common approach in the literature is the use of diffusion limits, and solving a multidimensional RBM is an open question. In this paper we contribute to understanding the heavy-traffic behavior of non-CRP systems by providing the mean of some linear combinations of the queue lengths.
	
	\item In \Cref{cor.nxn.switch} we compute the heavy-traffic limit of the total queue length in an input-queued switch with correlated arrivals. As mentioned above, the input-queued switch has had considerable attention in the literature. However, it has only been studied under independent arrival processes by \citeauthor{MagSri_SSY16_Switch} \cite{MagSri_SSY16_Switch}, and \citeauthor{QUESTA_switch} \cite{QUESTA_switch}. The input-queued switch is a model for an ideal data center network, and independent arrivals is an unrealistic assumption in this setting. In fact, as shown by \citeauthor{datacenter_traffic} \cite{datacenter_traffic} and \citeauthor{kandula_datacenter_traffic} \cite{kandula_datacenter_traffic}, for example, 
	data centers experience hot-spots, and, hence, the arrivals to different queues are highly correlated.
	
	\item We illustrate how \Cref{gs.thm:bounds} can be immediately applied to a variety of systems. Specifically, we show how to apply it to parallel-server systems (\Cref{gen.switch.cor.parallel-gral}, and \Cref{gen.switch.cor.n-queues}), the so-called $\cN$-system (\Cref{cor:Nsystem}) and ad hoc wireless networks (\Cref{cor:adhoc}). 
	
	\item In \Cref{subsec:full-dim-SSC} we show that, if SSC is full-dimensional, then the heavy-traffic limit of the mean queue lengths does not depend on the correlation among arrival processes. In other words, if the systems experience full-dimensional SSC, the expected linear combination of queue lengths behaves as if the queues were independent. This result is rather surprising, and it was not known.
	
	\item In \Cref{s.theorem:2switch} we show that using the drift method with polynomial test functions, it is impossible to obtain the moments of all linear combinations of the queue lengths. We prove this result by presenting an alternate way of thinking of the drift method. Traditionally, the key step in using this approach, is to design the correct test function to obtain all the moments. However, it is not clear a priori if there are test functions that give all these moments. Instead of trying to guess the right test function, this point of view shows that one can think about solving a set of linear equations. This system of linear equations turns out to be under-determined, and the major challenge is to obtain more equations using the constraints in the system, in order to solve for all the unknowns and obtain the complete joint distribution of the queue lengths when the CRP condition is not satisfied.
	
	\item In \Cref{s.2.thm.lp} we obtain lower and upper bounds on the steady-state mean of an arbitrary linear combination of queue lengths. We do this by formulating a Linear Program (LP) using the under-determined system of equations from \Cref{s.theorem:2switch}. We present numerical results in the case of Bernoulli arrivals, for different values of the traffic intensity. For simplicity of exposition, we do this only in the special case of an input-queued switch, and the same approach can be used for the generalized switch.
\end{enumerate}

The second, third and fourth contributions described above are proved as corollaries of the main theorem (\Cref{gs.thm:bounds}). However, they answer questions that, to the best of our knowledge, were open and, thus, they are contributions by themselves. This shows the versatility and power of \Cref{gs.thm:bounds}.

The system of equations we propose in the fifth contribution presents an alternate view of the drift method, and it explains its success. While it is known that it is notoriously hard to solve the stationary distribution of a multidimensional RBM, it has been a little surprising that simple drift-based arguments give the mean of sum of the queue lengths in several systems, such as the ones studied by \citeauthor{MagSri_SSY16_Switch} \cite{MagSri_SSY16_Switch}, \citeauthor{QUESTA_switch} \cite{QUESTA_switch} and \citeauthor{WeinaBandwidthJournal} \cite{WeinaBandwidthJournal}. The system of equations shows that due to the difficulty of the underlying problem, it is not possible to get all the mean queue lengths individually. However, because of the structure of the system of equations, it is possible to obtain certain linear combinations. In the case of input-queued switch and the bandwidth sharing system, the sum of the queue lengths was one of the linear combinations that is easy to obtain.

In the proof of \Cref{gs.thm:bounds} we use the drift method, and we work directly with the original queueing system, without any fluid or diffusion scaling. The drift method consists of two main steps: (1) Prove SSC and (2) Compute asymptotically tight bounds. We additionally compute a Universal Lower Bound (ULB) for a linear combination of the queue lengths, i.e., a lower bound on the queue lengths that does not depend on the scheduling policy. We establish SSC in terms of certain moment bounds using a Lyapunov drift argument. By definition of steady state, the drift of any function with finite expectation is zero. We pick a quadratic test function and set its drift to zero in steady state to obtain the result. The choice of this test function is important in the drift method. 
We use the norm of the projection of the queue length vector into the space of the SSC as our test function, which was also used by \citeauthor{atilla} \cite{atilla}, \citeauthor{MagSri_SSY16_Switch} \cite{MagSri_SSY16_Switch}, and \citeauthor{QUESTA_switch} \cite{QUESTA_switch}. In this paper we use this method in a discrete time system, but it can be also used in continuous time systems. For example,\citeauthor{WeinaBandwidthJournal}  \cite{WeinaBandwidthJournal} use the drift method in the context of a bandwidth sharing network, which operates in continuous time.

The organization of the rest of this paper is as follows. We start in \Cref{subsec:notation} establishing the main notation that will be used in this paper. In \Cref{sec:gen.switch.model} we describe the generalized switch model.   In \Cref{sec:gen.switch.heavy.traffic} we present the main result of this paper, along with SSC and the ULB. Specifically, in \Cref{sec:gen.switch.ulb} we present the ULB, in \Cref{sec:gen.switch.ssc} we show SSC, and in \Cref{sec:gen.switch.bounds} we present the main result of this paper (\Cref{gs.thm:bounds}), together with remarks that will help its interpretation. Then, in \Cref{sec:interpretation} we present relevant applications our result, including the study of the input-queued switch under correlated arrivals (\Cref{cor.nxn.switch}), full-dimensional SSC (\Cref{gen.switch.cor.full.ssc}) and the $\cN$-system (\Cref{cor:Nsystem}). In \Cref{sec:system.of.equations} we present the alternate view of the drift method in the context of an input-queued switch (\Cref{s.theorem:2switch}), and the linear programs to obtain bounds (\Cref{s.2.thm.lp}).  In \Cref{sec:proofs} we present the proof of \Cref{gs.thm:bounds} and \Cref{s.theorem:2switch}, and in \Cref{sec:conclusion} we present our conclusions and future work.

\subsection{Notation.}\label{subsec:notation}

In this subsection we introduce the notation that we use along the paper. We use $[n]$ to denote the set of integer numbers between 1 and $n$, both included. 
We use $\bR$ to denote the set of real numbers and $\bZ$ to denote the set of integer numbers. We add a subscript $+$ to denote nonnegativity, and a superscript with a number to denote the dimension. 
We use bold letters to denote vectors and non-bold letters with a subscript to denote their elements. 
We write $\vx=(x_1,x_2,\ldots,x_n)$ for convenience, but we treat vectors as column vectors unless otherwise stated. Given two vectors $\vx$ and $\vy$, we write $\langle\vx,\vy\rangle$ to denote the dot product between $\vx$ and $\vy$, and $\|\vx\|$ to denote the Euclidean norm. For a matrix $A$, we write $A^T$ to denote its transpose.  Given two matrices $A$ and $B$, we write $A\circ B$ to denote the Hadamard's product between $A$ and $B$, i.e., the matrix that results from multiplying term by term the elements of $A$ and $B$. 
Let $\mathbb{I}_n$ be the identity matrix of $n\times n$. We use $\ve^{(i,n)}$ to denote the $i\tth$ canonical vector in $\bR^n$, i.e., a vector with a 1 in the $i\tth$ element and zeros in all other entries. I the dimension is clear from the context we omit the $n$. 

%

For an irreducible and aperiodic Markov Chain  $\left\{X(k):k\in\bZ_+\right\}$ over a countable state space $\cX$, suppose $Z:\cX\to \bR_+$ is a Lyapunov function. Define the drift of $Z$ at $x$ as
\begin{align}\label{gs.eq.def.drift}
	\Delta Z(x)\defn \big[Z\left(X(k+1)\right)-Z\left(X(k)\right) \big]\ind{X(k)=x}.
\end{align}
Thus, $\Delta Z(x)$ is a random variable that measures the amount of change in the value of $Z$ in one time slot, starting from the state $x$.


\section{Generalized switch model.}\label{sec:gen.switch.model}
In this section we describe the model in detail and we state known stability results. Consider $n$ queues operating in discrete time, and let $\vq(k)$ be the vector of queue lengths at the beginning of time slot $k$, for each $k\in\bZ_+$. For each $i\in[n]$, let $\{a_i(k):k\in \bZ_+\}$ be a sequence of i.i.d. random variables such that $a_i(k)$ is the number of arrivals to the $i\tth$ queue in time slot $k$. For each $i\in[n]$ let $\lambda_i\defn\E{a_i(1)}$, and $\amax$ be a finite constant such that $a_i(1)\leq \amax$ with probability 1 for all $i\in[n]$. Let $\Sigma_a$ be the covariance matrix of the vector $\va(1)$. 

The servers interfere with each other, so in each time slot a set of interference constraints must be satisfied. Additionally, there are conditions of the environment that affect these constraints, which we group in a single random variable called channel state. In other words, given the channel state, the set of interference constraints is known and it may change if the channel state changes. Let $\{M(k): k\in\bZ_+\}$ be a sequence of i.i.d. random variables, such that $M(k)$ is the channel state in time slot $k$. Let $\cM$ be the state space of the channel state and $\vpsi$ be the probability mass function of $M(1)$, i.e., for each $m\in \cM$ we define $\psi_m\defn \Prob{M(1)=m}$. For each $m\in\cM$, define $\cS^\pm$ as the set of feasible service rate vectors in channel state $m$, that is, the set of vectors that satisfy the interference constraints in channel state $m$. For each $m\in\cM$ the set $\cS^\pm$ contains the projection on the coordinate axes of all its vectors. Formally, we assume that if $\vx\in\cS^\pm$ for some $m\in\cM$, then $\vx-x_i\ve^{(i)}\in \cS^\pm$ for all $i\in[n]$. We assume that $\cM$ is a finite set and that, for each $m\in\cM$, the set $\cS^\pm$ is finite. Therefore, the potential service offered to each queue in each time slot is bounded. Let $\smax$ be a finite upper bound.

Let $\vs\left(\vq(k),M(k)\right)$ be the vector of potential service in time slot $k$. In other words, for each $i\in[n]$, $s_i\left(\vq(k),M(k)\right)$ is the number of jobs from queue $i$ that would be processed if there were enough jobs in line. Let $\vu\left(\vq(k),M(k),\va(k)\right)$ be the vector of unused service in time slot $k$, i.e., for each $i\in[n]$, $u_i\left(\vq(k),M(k),\va(k) \right)$ is the difference between the potential service and the number of packets that are actually processed from queue $i$ in time slot $k$. For ease of exposition, and with a slight abuse of notation, from now on we use $\vs(k)$ and $\vu(k)$ to denote $\vs\left(\vq(k),M(k)\right)$ and $\vu\left(\vq(k),M(k),\va(k)\right)$, respectively.

In each time slot a scheduling problem must be solved to decide which queues will be served and at which rate. If queue $i$ is not scheduled to receive service in time slot $k$, then $s_i(k)=0$. In our model, the order of events in one time slot is as follows. First, the channel state is observed; second, a schedule is selected; third, arrivals occur and, at the end of each time slot, the jobs are processed according to the selected schedule. Hence, the dynamics of the queues are as follows. For each $k\in\bZ_+$ and $i\in[n]$ we have
\begin{align}\label{gs.eq.dynamics.queues}
	q_i(k+1)=q_i(k)+a_i(k)-s_i(k)+u_i(k).
\end{align}

In each queue, the unused service is nonzero only when the respective potential service is greater than the number of packets available (packets in line and arrivals). In such case, the queue is empty in the next time slot. Therefore, the following equation is satisfied with probability 1
\begin{align}\label{gs.eq.qu}
	q_i(k+1)u_i(k)=0\qquad\forall k\in\bZ_+,\;\forall i\in[n].
\end{align}
However, if $i\neq j$, then $q_i(k+1)u_j(k)$ is not necessarily zero.

In this paper we consider the generalized switch operating under MaxWeight algorithm, which means that in each time slot the schedule with the longest total weighted queue length is selected, where the possible weight vectors are the feasible service rate vectors. Formally, if $M(k)=m$ then
\begin{align}\label{gs.eq.MW}
	\vs(k)\in\argmax_{\vx\in \cS^\pm}\langle \vq(k),\vx\rangle,
\end{align}
and ties are broken at random.

It was proved by \citeauthor{atilla} \cite{atilla} that the capacity region of the generalized switch is
\begin{align}
	\cC=& \sum_{m\in \cM} \psi_m\,ConvexHull\left(\cS^\pm \right) \label{gs.eq.cap.reg.def.weighted.CH} \\
	=& ConvexHull\left(\left\{\sum_{m\in \cM}\psi_m \vx^\pm: \vx^\pm\in\cS^\pm\;\;\forall m\in\cM \right\}\right), \label{gs.eq.cap.reg.def.CH.weighted.sum}
\end{align}
and that MaxWeight algorithm is throughput optimal. Then, the generalized switch operating under MaxWeight is positive recurrent for all $\vlambda$ in the interior of $\cC$.

Since the sets $\cM$ and $\cS^\pm$ are finite for all $m\in\cM$, the capacity region $\cC$ is a polytope (bounded polyhedron) in $\bR^n_+$. Then, we can describe it as the intersection of finitely many half-spaces. Let $L$ be the minimum number of half-spaces that is needed to describe $\cC$ and, for each $\ell\in[L]$, let $\vcl$ and $\bl$ be the parameters that define the $\ell\tth$ half-space. Then,
\begin{align}\label{gs.eq.cap.reg.pol}
	\cC= \left\{\vx\in\bR^n_+:\langle \vcl,\vx \rangle\leq\bl\;,\,\ell=1,\ldots,L  \right\}.
\end{align}
Since the sets $\cS^\pm$ contain the projection of their elements on the coordinate axes, the capacity region $\cC$ is coordinate convex. Then, without loss of generality, we assume $\vcl\geq \vzero$ and $\bl>0$ for all $\ell\in [L]$. For ease of exposition, we also assume $\|\vcl\|=1$ for all $\ell\in[L]$. For each $\ell\in[L]$, let $\cF^\pl$ be the $\ell\tth$ facet of $\cC$, i.e., we define $\cF^\pl \defn\left\{\vx\in\cC:\langle\vcl,\vx\rangle=\bl \right\}$.

Observe that the schedules selected by MaxWeight algorithm do not necessarily belong to the capacity region $\cC$. This can be seen from \eqref{gs.eq.MW} and \eqref{gs.eq.cap.reg.def.weighted.CH} because $\psi_m\leq 1$ for all $m\in \cM$. However, the expected service rate vector does belong to the capacity region. We prove this result formally in \Cref{gs.lemma:schedule.in.C}.
\begin{lemma}\label{gs.lemma:schedule.in.C}
	Consider a generalized switch operating under MaxWeight as described above, and let $\Evq{\,\cdot\,}\defn\E{\,\cdot\,|\vq(k)=\vq}$. Then, $\Evq{\langle \vq(k),\vs(k)\rangle} = \max_{\vx\in \cC} \langle\vq,\vx\rangle$.
\end{lemma}
\proof{Proof of Lemma 1.}
Since $\vs(k)$ is selected using MaxWeight algorithm (see \eqref{gs.eq.MW}), we have
\begin{align*}
	\Evq{\langle \vq(k),\vs(k)\rangle}=& \Evq{\max_{\vx\in \cS^{(M(k))}}\langle \vq(k),\vx\rangle} \\
	\stackrel{(a)}{=}& \sum_{m\in \cM}\psi_m \max_{\vx\in \cS^\pm}\langle \vq,\vx\rangle
	\stackrel{(b)}{=} \max_{\vx\in\cC}\langle \vq,\vx\rangle,
\end{align*}
where $(a)$ holds because the channel state process is independent from the queue lengths process; and $(b)$ holds by definition of the capacity region $\cC$ presented in \eqref{gs.eq.cap.reg.def.weighted.CH}.
\Halmos \endproof

For technical reasons that will be apparent in \Cref{sec:gen.switch.bounds}, we introduce the following definition. For each $\ell\in[L]$ and $m\in\cM$ define the \textit{maximum $\vcl$-weighted service rate available when channel state is $m$} as
\begin{align}\label{gs.eq.cl.weighted.def}
	\bml=\max_{\vx\in\cS^\pm}\langle \vcl,\vx\rangle.
\end{align}
Observe that $\vcl$ and $\bml$ define a half-space that passes through the boundary of  $ConvexHull\left(\cS^\pm\right)$, but this half-space does not necessarily define a facet of $ConvexHull\left(\cS^\pm   \right)$. For each $\ell\in[L]$ and $k\in\bZ_+$, let $B_\ell(k)\defn b^{\left(M(k),\ell\right)}$. Notice that $B_\ell(k)$ is an i.i.d. sequence of random variables that satisfies $\Prob{B_\ell(1)=\bml}=\psi_m$ for each $m\in\cM$.  Let $\Sigma_B$ be the covariance matrix of the vector $\boldsymbol{B}(1)\defn \left(B_1(1),\ldots,B_L(1)\right)$, i.e., for each $\ell_1,\ell_2\in[L]$ we have
\begin{align*}
	\left(\Sigma_B\right)_{\ell_1,\ell_2} \defn \E{B_{\ell_1}(k)B_{\ell_2}(k)} - \E{B_{\ell_1}(k)}\E{B_{\ell_2}(k)}.
\end{align*}

We model heavy traffic as follows. We fix a vector $\vnu$ in the boundary of $\cC$ and we consider a set of generalized switches operating under MaxWeight as described above, parametrized by $\epsilon\in(0,1)$. The heavy-traffic limit is the limit as $\epsilon\downarrow0$ and, as $\epsilon$ gets small, the vector of mean arrival rates approaches $\vnu$. Formally, we parametrize the queueing system in the following way. We let $\vq^\peps(k)$, $\va^\peps(k)$, $\vs^\peps(k)$ and $\vu^\peps(k)$ be the vectors of queue lengths, arrivals, potential service and unused service, respectively, in time slot $k$, in the system parametrized by $\epsilon$. The parametrization is such that the vector of mean arrival rate is $\vlambda^\peps\defn \E{\va^\peps(1)}=(1-\epsilon)\vnu$. Therefore, $\vlambda^\peps$ belongs to the interior of $\cC$ for each $\epsilon\in(0,1)$ and, as $\epsilon\downarrow 0$, the arrival rate vector $\vlambda^\peps$ approaches the boundary of the capacity region at the point $\vnu$.

Heavy-traffic analysis of the generalized switch has been performed in the past, using the diffusion limits approach by \citeauthor{stolyar2004maxweight} \cite{stolyar2004maxweight}, and the Drift method by \citeauthor{atilla} \cite{atilla}. However, in both cases, the analysis is under the assumption that SSC occurs into a one-dimensional subspace (CRP condition), i.e., when the vector $\vnu$ is in the interior of a facet of the capacity region $\cC$. In this paper, we focus on cases where the vector $\vnu$ may live at the intersection of facets. Define $P\defn \left\{\ell\in [L]: \vnu\in\cF^\pl \right\}$, that is, $P$ is the set of indices of all the facets that intersect at $\vnu$. Observe that, if $P$ has only one element, we are under the CRP condition, and our results in this case agree with the results proved by \citeauthor{atilla} \cite{atilla}. In this paper we focus on the case where $P$ is allowed to have more than one element.

For each $\epsilon\in(0,1)$, let $\vqbar^\peps$ be a steady-state random vector such that the Markov chain $\left\{\vq^\peps(k):k\in\bZ_+ \right\}$ converges in distribution to $\vqbar^\peps$ as $k\uparrow\infty$. Since MaxWeight is throughput optimal, the Markov chain $\left\{\vq^\peps(k):k\in\bZ_+ \right\}$ is positive recurrent for each $\epsilon\in(0,1)$, so $\vqbar^\peps$ is well defined. Let $\vabar^\peps$ be a steady-state vector which is equal in distribution to $\va^\peps(1)$. Then, $\E{\vabar^\peps}=\vlambda^\peps$ and for each $i\in[n]$ we have $\abar^\peps_i\leq \amax$ with probability 1. Let $\Sigma_a^\peps$ be the covariance matrix of the vector $\vabar^\peps$. Let $\Mbar$ and $\Bbar_\ell$ be steady-state random variables that are equal in distribution to $M(1)$ and $B_\ell(1)$ for each $\ell\in[L]$, respectively. Let $\vsbar^\peps\defn \vsqm$ be the vector of potential service in steady-state, and $\vubar^\peps\defn\vuqma$ be the vector of unused service. Define $\left(\vqbar^\peps\right)^+\defn \vqbar^\peps+\vabar^\peps-\vsbar^\peps+\vubar^\peps$ as the vector of queue lengths one time slot after $\vqbar^\peps$ is observed, given that the vectors of arrivals and potential service are $\vabar^\peps$ and $\vsbar^\peps$, respectively.

In \Cref{sec:gen.switch.ssc} we prove that the state space collapses into the cone $\cK$ described below. In other words, we show that the vector of queue lengths can be approximated by a vector in $\cK$ in heavy traffic. Let $\cK$ be the cone generated by $\left\{\vcl:\ell\in P\right\}$ and $\cH$ be the subspace generated by the same set of vectors. Formally,
\begin{align}
	& \cK=\left\{\vx\in\bR^n_+:\;\vx=\sum_{\ell\in P}\xi_\ell \vcl\;, \;\xi_\ell\geq 0 \;\forall \ell\in P \right\}. \label{gs.eq.coneK}
\end{align}

\begin{figure}
	\centering
	\caption{Example of capacity region $\cC$ and cone $\cK$.}\label{gs.fig:cone}
	\includegraphics[width=0.3\linewidth]{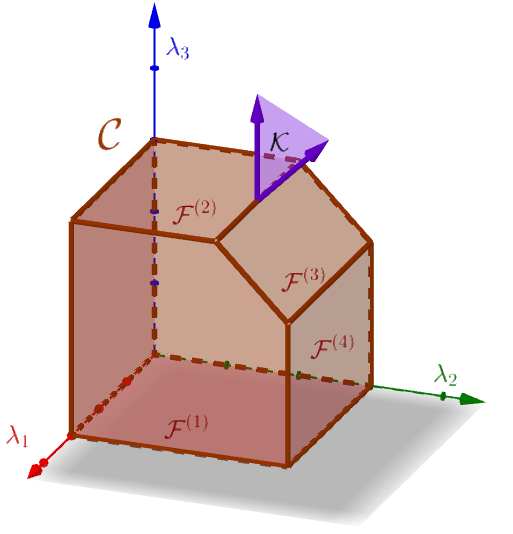}
\end{figure}

A pictorial example of the capacity region $\cC$ and the cone $\cK$ when $n=3$ is presented in \Cref{gs.fig:cone}.
Let $\tilde{P}\subset P$ be a set of indices such that the set $\left\{\vcl:\ell\in \tilde{P} \right\}$ is linearly independent, and let $C=\left[\vcl\right]_{\ell\in \tilde{P}}$ be a matrix where the columns are a linearly independent subset of the vectors that generate the cone $\cK$. Observe that the column space of the matrix $C$ is exactly the subspace $\cH$.

\section{Heavy-traffic analysis of the generalized switch.}\label{sec:gen.switch.heavy.traffic}

In this section we perform heavy-traffic analysis of the generalized switch. In \Cref{sec:gen.switch.ulb} we present a ULB, that is independent of the scheduling policy; in \Cref{sec:gen.switch.ssc} we present the SSC result formally; and in \Cref{sec:gen.switch.bounds} we present the main result of this paper (\Cref{gs.thm:bounds}), where we compute asymptotically tight bounds on linear combinations of the queue lengths.

\subsection{Universal lower bound.}\label{sec:gen.switch.ulb}
In this section we compute a ULB for certain linear combinations of the vector of queue lengths. The bound is universal in the sense that it remains valid for all scheduling policies.

\begin{proposition}\label{gs.prop.ulb}
	Consider a generalized switch parametrized by $\epsilon\in(0,1)$, as described in  \Cref{sec:gen.switch.model}. Let $\vz\in\cK$ with $\vz\neq \vzero$ and $\vr\in\bR_+^{|P|}$ be such that $\vz=\langle \vr, \vcl\rangle$. Then, for each $\epsilon\in(0,1)$ we have
	\begin{align*}
		& \E{\langle\vz,\vqbar^\peps\rangle}
		\geq  \dfrac{1}{2\epsilon\langle\vz,\vnu\rangle}\left( \vz^T\Sigma_a^\peps\vz + \vr^T \Sigma_B\vr \right)-f(\epsilon),
	\end{align*}
	where $f(\epsilon)=\dfrac{b_{\max}\langle \vone, \vr\rangle}{2}-\dfrac{\epsilon\langle\vz,\vnu\rangle}{2}$ is $o\left(\frac{1}{\epsilon} \right)$ (i.e., $\lim_{\epsilon\dto 0}\epsilon f(\epsilon)=0$) and 
	$b_{\max}=\max_{m\in\cM,\ell\in P}\bml$.
\end{proposition}

The proposition is proved by coupling the queue length vector of the generalized switch with a single server queue 
$\left\{\Phi^\peps(k): k\in\bZ_+ \right\}$ constructed as follows. We let $\alpha^\peps(k)\defn \langle\vz,\va^\peps(k)\rangle$ be the number of arrivals in time slot $k$ and $\beta(k)$ be the potential service, where $\Prob{\beta(k)=\sum_{\ell\in P}r_\ell\bml}=\psi_m$ for each $m\in \cM$. Then, it is easy to see that $\Phi^\peps(k)$ is stochastically smaller than $\langle\vz,\vq^\peps(k)\rangle$ (by definition of $\bml$ in \eqref{gs.eq.cl.weighted.def}). Therefore, a lower bound to the expected value of $\Phi^\peps(k)$ in steady state is also a lower bound to $\E{\langle\vz,\vqbar^\peps\rangle}$. The last step in the proof is to compute such lower bound, which we do by setting to zero the drift of $V_{ULB}(\Phi)=\Phi^2$. 
Note that it is essential that the weights $r_\ell$ are nonnegative to obtain a lower bound in the proof. This is the reason why $\vz\in\cH$ is not enough and we require $\vz\in\cK$. The rest of the proof is presented in Appendix \ref{app:prop.gs.ulb}.

\subsection{State space collapse.}\label{sec:gen.switch.ssc}

We prove SSC into the cone $\cK$ defined in \eqref{gs.eq.coneK} in heavy traffic. We start introducing the notation. For each $\epsilon\in(0,1)$, let $\vqpark^\peps(k)$ be the projection of $\vq^\peps(k)$ on $\cK$ and $\vqperpk^\peps(k)\defn \vq^\peps(k)-\vqpark^\peps(k)$. Similarly, define $\vqparh^\peps(k)$ as the projection of $\vq^\peps(k)$ on $\cH$ and $\vqperph^\peps(k)\defn \vq^\peps(k)-\vqparh^\peps(k)$. We know the Markov chain $\{\vq^\peps(k): k\in\bZ_+\}$ is positive recurrent for each $\epsilon\in(0,1)$, so by definition of projection we also have that $\{\vqpark^\peps(k): k\in\bZ_+\}$, $\{\vqperpk^\peps(k): k\in\bZ_+\}$, $\{\vqparh^\peps(k): k\in\bZ_+\}$ and $\{\vqperph^\peps(k): k\in\bZ_+\}$ are positive recurrent for each $\epsilon\in(0,1)$. Then, we define $\vqbarpark^\peps$, $\vqbarperpk^\peps$, $\vqbarparh^\peps$ and $\vqbarperph^\peps$ as steady-state vectors which are limit in distribution of each them, respectively. In the next proposition we state SSC formally.

\begin{proposition}\label{gs.prop:gen.switch.SSC}
	Given a vector $\vnu$ in the boundary of $\cC$ and $\epsilon\in(0,1)$, consider a generalized switch operating under MaxWeight, parametrized by $\epsilon$ as described in \Cref{sec:gen.switch.model}, and let $P$ be defined as in \Cref{sec:gen.switch.model} as well. Let $\delta>0$ be such that $\delta\leq \bl-\langle\vcl,\vnu\rangle$ for all $\ell\in [L]\setminus P$ if $[L]\setminus P\neq\emptyset$, and $\delta=1$ if $[L]\setminus P=\emptyset$. If $\epsilon<\frac{\delta}{2\|\vnu\|}$, then for each $t=1,2,\ldots$ there exists a constant $T_t$ such that $\E{\|\vqbarperph^\peps \|^t}\leq\E{\|\vqbarperpk^\peps \|^t}\leq T_t$.
\end{proposition}

We provide the proof of \Cref{gs.prop:gen.switch.SSC} in Appendix \ref{app:prop.gs.ssc}. We adopt the technique introduced by \citeauthor{atilla} \cite{atilla}, so our proof is similar to theirs. The challenges in obtaining our result arise in the second step of the drift method, which corresponds to \Cref{gs.thm:bounds}.

SSC is a consequence of \Cref{gs.prop:gen.switch.SSC} for the following reason. As $\epsilon\downarrow0$, $\left\|\vqbar^\peps\right\|$ goes to infinity (this can be easily concluded from \Cref{gs.thm:bounds}). Therefore, \Cref{gs.prop:gen.switch.SSC} implies that as $\epsilon$ gets small, we can approximate $\vqbar^\peps\approx\vqbarpark^\peps$ because all the moments of $\|\vqbarperpk^\peps\|$ are bounded.

Observe that the cone $\cK$ is determined by the facets that intersect at $\vnu$. Moreover, the dimension of the cone is $n-d_{\vnu}$, where $d_{\vnu}$ is the dimension of the face of $\cC$ where $\vnu$ is. For example, if $\vnu$ is in the relative interior of a facet then $d_{\vnu}=n-1$, and this implies that $\cK$ is one-dimensional. This is the CRP case, which was studied by \citeauthor{atilla} \cite{atilla} and \citeauthor{stolyar2004maxweight} \cite{stolyar2004maxweight}. Similarly, if $\vnu$ is a vertex of $\cC$ then $d_{\vnu}=0$ and, hence, $\cK$ is $n$-dimensional. In the last case, we say that SSC is full dimensional. We study the full-dimensional case in \Cref{subsec:full-dim-SSC}.

\subsection{Asymptotically tight bounds.}\label{sec:gen.switch.bounds}

In \Cref{sec:gen.switch.ssc} we showed SSC into the cone $\cK$, which implies SSC into the subspace $\cH$. In this section we present the main result of this paper (\Cref{gs.thm:bounds}), where we provide asymptotically tight bounds to the expected value of certain linear combinations of the queue lengths in steady state. After the statement of the theorem we present some remarks and applications, and we delay the proof to \Cref{subsec:proof.gen.switch}.

\begin{theorem}\label{gs.thm:bounds}
	Given a vector $\vnu$ in the boundary of $\cC$, let $P$ be defined as in \Cref{sec:gen.switch.model}. Consider a set of generalized switches operating under MaxWeight, indexed by the heavy-traffic parameter $\epsilon\in(0,1)$ as described in \Cref{sec:gen.switch.model}. Then, for any vector $\vw\in\cap_{\ell\in P}\cF^\pl$ we have
	\begin{align}\label{eq.gen.switch.thm.prelimit}
		\left|\E{\langle\vqbar^\peps, \vw\rangle}- \dfrac{1}{2\epsilon} \vone^T\left(H\circ \Sigma_a^\peps\right)\vone - \dfrac{1}{2\epsilon} \vone^T\left((C^TC)^{-1}\circ \Sigma_B \right)\vone \right|&\leq K(\epsilon),
	\end{align}
	where $\epsilon K(\epsilon)$ converges to 0 as $\epsilon\downarrow 0$ and $H\defn C(C^T C)^{-1}C^T$ is the projection matrix into $\cH$. Further, suppose $\lim_{\epsilon\dto 0}\Sigma_a^\peps= \Sigma_a$ component-wise. Then,
	\begin{align}\label{eq.gen.switch.thm.limit}
		\begin{aligned}
			& \lim_{\epsilon\dto 0} \epsilon\, \E{\langle \vqbar^\peps,\vw\rangle}
			= \dfrac{1}{2}\left(\vone^T\left(H\circ \Sigma_a\right)\vone +\vone^T\left((C^TC)^{-1}\circ \Sigma_B \right)\vone \right).
		\end{aligned}
	\end{align}
\end{theorem}

First observe that \eqref{eq.gen.switch.thm.prelimit} gives bounds that are valid for all regimes, not necessarily heavy traffic. Additionally, it shows that the queue lengths grow to infinity as the traffic intensity grows (i.e., as $\epsilon\dto 0$). 

In \eqref{eq.gen.switch.thm.limit}, observe that the right-hand side has two terms: one corresponding to randomness in the arrival process, and the other one to randomness in the service process. The first term is a linear combination of the covariance matrix of the arrival process, and the weights of the linear combination are determined by the projection matrix on the subspace $\cH$, which is where SSC occurs. The second term is a linear combination of the elements of a covariance matrix which is related to the channel state. Since the potential service rate vector is selected using MaxWeight algorithm (see \eqref{gs.eq.MW}), it is not actually random once queue lengths and channel state are observed. However, the channel state is a random variable that defines the feasible set where MaxWeight is solved. Hence, the second term in \eqref{eq.gen.switch.thm.limit}, which includes a covariance matrix related to channel state, represents the randomness on the service process.

A third observation is that, in order to project on the subspace $\cH$ generated by the cone $\cK$, we had to drop the vectors $\vcl$ with $\ell\in P$ that are linearly dependent (recall that the columns of the matrix $C$ are a linearly independent subset of the vectors that generate $\cK$). Clearly, the cone generated by the columns of $C$ is not equal to $\cK$. However, projecting on the subspace $\cH$ is sufficient, and we do not need to worry about these linearly dependent vectors that we dropped. 

In the next remark we write \eqref{eq.gen.switch.thm.limit} in different ways to facilitate interpretation of the result.

\begin{remark}\label{gs.remark.limit.epar}
	\Cref{eq.gen.switch.thm.limit} can be also written as
	\begin{align}
		\lim_{\epsilon\dto 0} \epsilon\E{\langle\vqbar^\peps,\vw\rangle} =& \dfrac{1}{2}\left(\sum_{i=1}^n\sum_{j=1}^n \langle\ve^{(i)}, \ve^{(j)}_{\parallel \cH}\rangle (\Sigma_{a})_{i,j} + \sum_{\ell_1\in \tilde{P} }\sum_{\ell_2\in \tilde{P}} (C^TC)^{-1}_{\ell_1,\ell_2}(\Sigma_B)_{\ell_1,\ell_2} \right) \label{eq.gen.switch.thm.limit.epar} \\
		=& \dfrac{1}{2}\Bigg(Trace\left(H\Sigma_a^T \right) + Trace\left((C^TC)^{-1}\Sigma_B^T \right) \Bigg) \label{eq.gen.switch.thm.limit.trace},
	\end{align}
	where the subscript $\parallel\hspace{-0.05in}\cH$ denotes projection on the subspace $\cH$, $(\Sigma_a)_{i,j}$ is the element $(i,j)$ of the covariance matrix $\Sigma_a$ for each $i,j\in[n]$, and $(\Sigma_B)_{\ell_1,\ell_2}$ is the element $(\ell_1,\ell_2)$ of $\Sigma_B$ for each $\ell_1,\ell_2\in \tilde{P}$.
	
	In some cases, the projection of a vector on $\cH$ is known in closed form, and it is simpler to work with than the projection matrix. For example, in the case of a completely saturated input-queued switch, \citeauthor{MagSri_SSY16_Switch} \cite{MagSri_SSY16_Switch} directly compute the projections, but writing down the projection matrix is more involved. 
\end{remark}

We present the proof of \Cref{gs.remark.limit.epar} below.
\proof{Proof of Remark \ref{gs.remark.limit.epar}.}
If we expand the products on the right-hand side of \eqref{eq.gen.switch.thm.limit} we obtain
\begin{align*}
& \dfrac{1}{2}\left(\vone^T\left(H\circ \Sigma_a\right)\vone + \vone^T\left((C^TC)^{-1}\circ \Sigma_B \right)\vone \right) \\
\stackrel{(a)}{=}& \dfrac{1}{2} \left(\sum_{i=1}^n \sum_{j=1}^n h_{i,j}(\Sigma_a)_{i,j} + \sum_{\ell_1\in \tilde{P} }\sum_{\ell_2\in \tilde{P}} (C^TC)^{-1}_{\ell_1,\ell_2}(\Sigma_B)_{\ell_1,\ell_2} \right) \\
\stackrel{(b)}{=}& \dfrac{1}{2}\left(\sum_{i=1}^n \sum_{j=1}^n \left(\ve^{(i)}\right)^T H\ve^{(j)} (\Sigma_a)_{i,j} + \sum_{\ell_1\in \tilde{P} }\sum_{\ell_2\in \tilde{P}} (C^TC)^{-1}_{\ell_1,\ell_2}(\Sigma_B)_{\ell_1,\ell_2} \right) \\
\stackrel{(c)}{=}& \dfrac{1}{2}\left(\sum_{i=1}^n \sum_{j=1}^n \langle\ve^{(i)},\ve^{(j)}_{\parallel \cH}\rangle (\Sigma_a)_{i,j} + \sum_{\ell_1\in \tilde{P} }\sum_{\ell_2\in \tilde{P}} (C^TC)^{-1}_{\ell_1,\ell_2}(\Sigma_B)_{\ell_1,\ell_2} \right),
\end{align*}
where $(a)$ holds by definition of Hadamard's product; $(b)$ holds by definition of the canonical vectors $\ve^{(i)}$ and by definition of matrix product; and $(c)$ holds by definition of inner product and because $H\ve^{(j)}$ is the projection of $\ve^{(j)}$ on the subspace $\cH$.

The proof of \eqref{eq.gen.switch.thm.limit.trace} holds by properties of Hadamard's product and trace, and we omit it.
\Halmos
\endproof

Observe that the bounds presented in \Cref{gs.prop.ulb} and \Cref{gs.thm:bounds} may be for different linear combinations of the vector of queue lengths. In \Cref{gs.prop.ulb} the vector of weights is $\vz\in\cK$ and in \Cref{gs.thm:bounds} it is $\vw\in\cap_{\ell\in P}\cF^\pl$. In the next remark we give sufficient conditions under which these bounds correspond to the same linear combination of the queue lengths.

\begin{remark}\label{gs.remark.farkas}
	Let $A$ be a matrix with columns $\vcl$ for $\ell\in P$ and $\vb_P$ be a vector with elements $\bl$ for $\ell\in P$. Observe that the column space of $A$ is equal to the column space of $C$, but the columns of $A$ may not be linearly independent. In fact, if the columns of $A$ are linearly independent, then $A=C$. Then, \Cref{gs.prop.ulb} and \Cref{gs.thm:bounds} give bounds to the same linear combination of the queue lengths if the set $\cA\defn \left\{\vx\in\bR^{|P|}:\vx^T A^TA\geq 0\;,\; \vx^T\vb_P<0 \right\}$ is empty.
\end{remark}

\proof{Proof of \Cref{gs.remark.farkas}.}
	We can obtain bounds to the same linear combination of the queue lengths if there exists a vector $\vy\in \cK\cap\left(\cap_{\ell\in P}\cF^\pl\right)$. In other words, if the set $\mathcal{Y}\defn \left\{\vy\in\bR_+^{|P|}: AA^T\vy=\vb_P \right\}$ is nonempty. By Farkas' lemma (Theorem 4.6 by \citeauthor{bertsimas_LPbook} \cite{bertsimas_LPbook}), proving that $\cY\neq\emptyset$ is equivalent to proving that $\cA=\emptyset$.
\Halmos\endproof

In the proof of \Cref{gs.thm:bounds} we use the drift method, which is a two-step procedure to compute bounds on linear combinations of the queue lengths that are tight in heavy traffic. The first step is to prove SSC, which we did in \Cref{gs.prop:gen.switch.SSC}, and the second step is to set to zero the drift of $V(\vq)=\|\vq_{\parallel \cH}\|^2$. While these steps are standard for the drift method as developed by \citeauthor{atilla} \cite{atilla}, \citeauthor{MagSri_SSY16_Switch} \cite{MagSri_SSY16_Switch}, \citeauthor{QUESTA_switch} \cite{QUESTA_switch}, and \citeauthor{WeinaBandwidthJournal} \cite{WeinaBandwidthJournal}, different challenges arise in each case depending on the system one is studying. In this case we are working with the generalized switch, which is a very general model. Hence, we overcome difficulties that were not part of the work listed above. We summarize these below:
\begin{enumerate}[label=(\alph*)]
	\item Since the effective capacity region is the average of several individual capacity regions (see the definition of the capacity region in \eqref{gs.eq.cap.reg.def.weighted.CH}), the vector of potential service does not necessarily belong to the effective capacity region. Then, it is not obvious how to deal with the terms that involve the service vector.
	\item 
	In the case of an input-queued switch as studied by \citeauthor{MagSri_SSY16_Switch} \cite{MagSri_SSY16_Switch}, the projected service vector $\vsbarparh$ is constant due to the structure of the system. In the case of the generalized switch, this is not the case, and this leads to significant challenges. In particular, the computation of the term $\E{\langle\vqbarparh,\vsbarparh\rangle}$ is not trivial. We used the properties of the system and MaxWeight algorithm to bound this term.
	\item The final closed-form expression that we obtain for the steady-state expectation of the queue lengths is novel, and is a contribution in itself. To compute this expression (the term $\cT_2$ in the proof), we use the least squares problem  to obtain an expression that is valid for any generalized switch. 
	\citeauthor{atilla} \cite{atilla}, \citeauthor{MagSri_SSY16_Switch} \cite{MagSri_SSY16_Switch} and \citeauthor{QUESTA_switch} \cite{QUESTA_switch} explicitly use the underlying symmetry of the specific systems that were studied, and therefore, it is not clear how to generalize.
\end{enumerate}

Challenge (a) was addressed in Lemmas \ref{gs.lemma:bjl.bl} and \ref{gs.lemma:cl.s.bl}. These lemmas form an important part of the entire proof and are used repeatedly. Challenge (b) was addressed in \Cref{gs.claim.T1.partial}. Finally, overcoming challenge (c) using the least square problem gave us the closed form expression for the right-hand side in \Cref{gs.thm:bounds}.

\section{Applications of \texorpdfstring{\Cref{gs.thm:bounds}}{Theorem 1}.}\label{sec:interpretation}

The generalized switch is a model that subsumes several SPNs, such as ad hoc wireless networks, the input-queued switch, down-link base stations and the parallel-server system.  In this section we elaborate on a few applications to give examples of the use of \Cref{gs.thm:bounds}, and it is by no means an exhaustive list.
We start with an input-queued switch in \Cref{subsec:switch}, and then, in \Cref{subsec:full-dim-SSC}, we present examples where full-dimensional SSC is observed.


\subsection{Input-queued switch.}\label{subsec:switch}

The drift method has been used to perform heavy-traffic analysis of the input-queued switch operating under MaxWeight in both, completely and incompletely saturated cases by \citeauthor{MagSri_SSY16_Switch} \cite{MagSri_SSY16_Switch} and \citeauthor{QUESTA_switch} \cite{QUESTA_switch}, respectively. In both scenarios, the analysis is performed under the assumption that the arrivals to different queues are independent. However, this is an unrealistic assumption in data center networks. Indeed, it has been shown that the traffic exhibits hot-spots, i.e., there are subsets of queues that simultaneously perceive a surge on traffic, as shown by  \citeauthor{datacenter_traffic} \cite{datacenter_traffic} and \citeauthor{kandula_datacenter_traffic}  \cite{kandula_datacenter_traffic}. This implies that the arrival processes are highly correlated.
In this section we focus on the completely saturated input-queued switch, and we obtain the heavy-traffic limit of the scaled total queue length when the arrivals are correlated, as a corollary of \Cref{gs.thm:bounds}. \Cref{cor.nxn.switch} generalizes the main result proved by \citeauthor{MagSri_SSY16_Switch} \cite{MagSri_SSY16_Switch} and it is of special interest by itself, given the nature of the arrival processes to data center networks observed in reality. We start specifying the model.

Consider a system with $N^2$ queues operating in discrete time. There are $N$ input ports, $N$ output ports, and there is a different queue for each input/output pair. Each of these pairs has its own arrival process and all the arriving packets have the same size, which is equal to one time slot. The service process must satisfy the following feasibility constraints. In each time slot, at most one packet can be transmitted from each input port, and each output port can process at most one packet. We can think of this system as a matrix of input/output pairs, where rows represent inputs and columns represent outputs. Then, the constraint described above can be also stated as follows. In each time slot, at most one queue can be active (i.e., processing jobs) in each row and each column.

\begin{figure}
	\centering
	\caption{Diagram of the queue length vector for the input-queued switch.\label{pic:switch}}
	\begin{subfigure}[t]{0.2\linewidth}
		\centering
		\includegraphics[width=\linewidth]{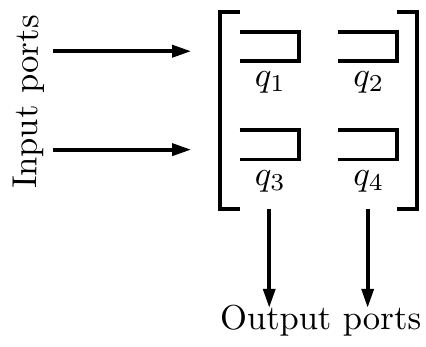}
		\caption{$2\times 2$ switch.\label{pic:switch.2x2}}
	\end{subfigure}
	\hspace{.5in}
	\begin{subfigure}[t]{0.25\linewidth}
		\centering
		\includegraphics[width=\linewidth]{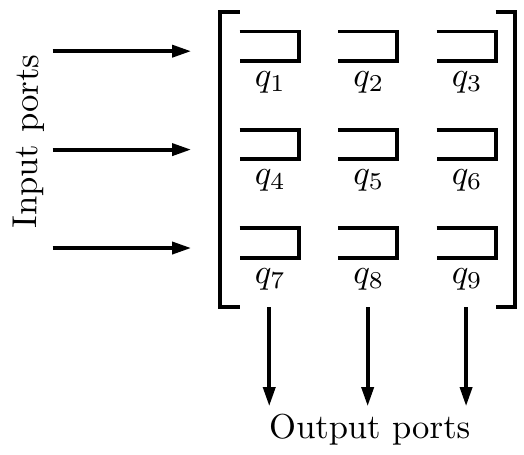}
		\caption{$3\times 3$ switch.\label{pic:switch.3x3}}
	\end{subfigure}
\end{figure}

This model corresponds to a generalized switch with $n=N^2$ queues, where the channel state is constant over time.
As mentioned above, the input-queued switch has a natural matrix-shape interpretation. \citeauthor{MagSri_SSY16_Switch}\cite{MagSri_SSY16_Switch} and \citeauthor{QUESTA_switch} \cite{QUESTA_switch} represent the vectors of queue lengths, arrivals and services by $N\times N$ matrices, but they are treated as vectors. Specifically, dot products and norms are computed as if these matrices were column vectors. In this paper, however, we will write them as column vectors to be consistent with the notation we introduced in \Cref{sec:gen.switch.model}. We enumerate the elements of the vectors row by row. For each $i\in[n]$ we have that $q_i(k)$ is the number of packets in line in input port $\left\lceil\frac{i}{N} \right\rceil$, waiting for service from output $(i\mymod N)$ if $i$ is not a multiple of $N$, and output $N$ otherwise. Similarly for the vectors of arrivals, potential service and unused service. In \Cref{pic:switch} we show how to build the vectors in the case of $N=2$ (\Cref{pic:switch.2x2}) and $N=3$ (\Cref{pic:switch.3x3}).

For ease of exposition, we introduce the following notation. For each $i\in[N^2]$ let
\begin{align*}
	& row(i)\defn \left\{\left(\left\lceil\frac{i}{N} \right\rceil-1 \right)N+j: j\in[N] \right\}\setminus \{i\} \\[3pt]
	& col(i)\defn \left\{j\in\left[N^2\right]: i\mymod N=j\mymod N \right\}\setminus \{i\} \\[3pt]
	& other(i)\defn \left[N^2\right]\setminus \big(row(i)\cup col(i)\cup \{i\} \big).
\end{align*}
In words, the set $row(i)$ contains the index of all elements in the same row as $i$, except by $i$; $col(i)$ contains the index of the elements in the same column as $i$, except by $i$; and $other(i)$ contains all indexes that do not correspond to the same row or column as $i$, or $i$ itself.

We explicitly know the feasibility constraints in the input-queued switch. Then, we can compute the set of feasible service rate vectors $\cS$ and the capacity region $\cC$. We obtain
\begin{align*}
	\cS= \left\{\vx\in\{0,1\}^{N^2}:\; \sum_{i=1}^N x_{N(j-1)+i}\leq 1\;\forall j\in [N]\; \text{and}\; \sum_{j=1}^N x_{N(j-1)+i}\leq 1\;\forall i\in[N] \right\},
\end{align*}
and
\begin{align}
	\cC& =ConvexHull(\cS) \nonumber\\
	&= \left\{\vx\in\bR_+^{N^2}:\; \sum_{i=1}^N x_{N(j-1)+i}\leq 1\;\forall j\in [N]\;\text{and}\; \sum_{j=1}^N x_{N(j-1)+i}\leq 1\;\forall i\in[N] \right\}. \label{eq.soe.N.C}
\end{align}
Then, the number of hyperplanes that define the capacity region is $L=2N$, the right-hand side parameters are $b^\pl=1$ for all $\ell\in[2N]$, and the left-hand side vectors $\vc^\pl$ are defined as follows.
\begin{align}\label{gs.eq.def.cl.switch}
	\vcl&= \begin{cases}
		\ds\sum_{i=N(\ell-1)+1}^{N\ell} \ve^{(i)} & \text{, if }\ell\in [N] \\
		\ds\sum_{i\in\{i':i'\mymod N = \ell\mymod N\}} \ve^{(i)} & \text{, if }\ell\in[2N]\setminus [N].
	\end{cases}
\end{align}

Completely saturated switch means that the vector $\vnu$ that we approach in the heavy-traffic limit satisfies all the inequalities in \eqref{eq.soe.N.C} at equality. Formally, $\vnu$ satisfies $\langle\vcl,\vnu\rangle=b^\pl$ for all $\ell\in[2N]$. Then, $P=[2N]$. If $\vnu$ does not satisfy all the inequalities at equality, it is said that the switch is incompletely saturated. We do not study the incompletely saturated case in this paper.

Recall that the cone $\cK$ where SSC occurs is the cone generated by the vectors $\vcl$ with $\ell\in P$. In this case, since $P=[2N]$ and since we explicitly know the vectors $\vcl$, it can be easily proved that the cone $\cK$ can be described as
\begin{align}\label{eq.switch.cone}
	\cK=& \left\{\vx\in\bR^{N^2}_+:\; x_i= \dfrac{1}{N}\sum_{j\in row(i)\cup \{i\}} x_j + \dfrac{1}{N}\sum_{j\in col(i)\cup \{i\}} x_j - \dfrac{1}{N^2}\sum_{j=1}^{N^2}x_j \right\}.
\end{align}
The proof of this claim is just algebra, and we omit it for brevity. In this case, it can be also proved that the subspace $\cH$ generated by the cone $\cK$ satisfies $\cK=\cH\cap \bR^{N^2}_+$.

Now we present the heavy-traffic limit of the scaled sum of the queue lengths in a completely saturated switch with correlated arrival processes, as a corollary of \Cref{gs.thm:bounds}. This corollary by itself is a contribution of this paper because, to the best of our knowledge, the input-queued switch has been studied  only under independent arrivals assumption. However, it is known that in data centers this is not satisfied and, in fact, hot-spots are frequently observed. 

\begin{corollary}\label{cor.nxn.switch}
	Let $\vnu$ be an $N^2$-dimensional vector that satisfies $\langle\vcl,\vnu\rangle=\bl$ for all $\ell\in [2N]$, for $\vcl$ as defined in \eqref{gs.eq.def.cl.switch} and $\bl=1$ for all $\ell\in [2N]$. Consider a set of $N\times N$ input-queued switches as described above, parametrized by $\epsilon\in(0,1)$ as described in \Cref{gs.thm:bounds}. For each $i\in[N^2]$, let $\sigma_{a_i}^2=\left(\Sigma_a\right)_{i,i}$. Then,
	\begin{align*}
		& \lim_{\epsilon\dto 0}\epsilon \E{\sum_{i=1}^{N^2} \qbar^\peps_i }
		= \dfrac{1}{2N}\sum_{i=1}^{N^2} \left((2N-1)\sigma_{a_i}^2 + (N-1)\sum_{j\in row(i)\cup col(i)}(\Sigma_a)_{i,j}-\sum_{j\in other(i)}(\Sigma_a)_{i,j}\right).
	\end{align*}
\end{corollary}
\proof{Proof of Corollary \ref{cor.nxn.switch}.}
We use Remark \ref{gs.remark.limit.epar}. We first compute $\ve^{(i)}_{\parallel\cH}$ for each $i\in[N^2]$.  For any vector $\vy\in\bR^{N^2}_+$ we have $\vy_{\parallel\cH}$ has elements
\begin{align*}
	\vy_{\parallel \cH j}= \dfrac{1}{N}\sum_{j'\in row(j)\cup \{j\}}y_{j'} + \dfrac{1}{N}\sum_{j' \in col(j)\cup \{j\}}y_{j'}-\dfrac{1}{N^2}\sum_{j'=1}^{N^2} y_{j'}\qquad\forall j\in[N^2].
\end{align*}
Then, for each $i\in [N^2]$ the vector $\ve^{(i)}_{\parallel\cH}$ has elements
\begin{align*}
	\ve^{(i)}_{\parallel\cH j}= \begin{cases}
		\frac{2N-1}{N^2} & \text{, if }j=i \\
		\frac{N-1}{N^2} & \text{, if $j\in row(i)$ or $j\in col(i)$} \\
		-\frac{1}{N^2} & \text{, if $j\in other(i)$}
	\end{cases}
	\qquad\forall j\in[N^2].
\end{align*}
Using this expression in Remark \ref{gs.remark.limit.epar} we immediately obtain the result.
\Halmos
\endproof

\begin{corollary}\label{cor.nxn.switch.indep}
	Consider a set of $N\times N$ input-queued switches operating under MaxWeight, parametrized by $\epsilon\in(0,1)$ as described in \Cref{cor.nxn.switch}. Further, assume that the arrival processes to different queues are independent. Then,
	\begin{align*}
		\lim_{\epsilon\dto 0}\epsilon \E{\sum_{i=1}^{N^2} \qbar^\peps_i}= \left(1-\dfrac{1}{2N}\right)\sum_{i=1}^{N^2}\sigma_{a_i}^2.
	\end{align*}
\end{corollary}

The proof of \Cref{cor.nxn.switch.indep} is easy after considering \Cref{cor.nxn.switch}, since $(\Sigma_a)_{i,j}=0$ for all $i\neq j$ under the independent arrivals assumption. \Cref{cor.nxn.switch.indep} recovers the main result presented by \citeauthor{MagSri_SSY16_Switch} \cite{MagSri_SSY16_Switch}, where they explicitly set to zero the drift of $V_{\parallel\cH} (\vq)=\left\|\vq_{\parallel \cH} \right\|^2$ (similarly to our approach in the proof of \Cref{gs.thm:bounds}).


\subsection{Full-dimensional SSC.}\label{subsec:full-dim-SSC}

As mentioned in \Cref{sec:gen.switch.ssc}, if the point $\vnu$ is a vertex of the capacity region $\cC$, the cone $\cK$ is $n$-dimensional. In other words, $\cK$ is full-dimensional. In this section we explore this situation, and we present examples of SPNs where this phenomenon is observed. In particular, we present the case of a parallel-server system operating in discrete time in \Cref{subsec:parallel-server-system}, an $\cN$-system in \Cref{subsec:Nsystem}, and an ad hoc wireless network in \Cref{subsec:adhoc-wireless-network}. 
We first present the result in a general case.
\begin{corollary}\label{gen.switch.cor.full.ssc}
	Consider a set of generalized switches operating under MaxWeight, parametri-\newline zed by $\epsilon\in(0,1)$ as described in \Cref{gs.thm:bounds}. Let $P$, $\tilde{P}$ and $\vnu$ be as in \Cref{gs.thm:bounds} and suppose the cone $\cK$ is $n$-dimensional. Let $\sigma_{a_i}^2\defn(\Sigma_a)_{i,i}$ for each $i\in[n]$ and $\vsigma_a$ be a vector with elements $\sigma_{a_i}$. Then,
	\begin{align*}
		\lim_{\epsilon\dto 0} \epsilon\E{\langle\vqbar^\peps,\vw\rangle}=& \dfrac{1}{2}\left(\left\|\vsigma_a \right\|^2 + \vone^T\left((C^TC)^{-1}\circ \Sigma_B \right)\vone \right).
	\end{align*}
\end{corollary}

Observe that \Cref{gen.switch.cor.full.ssc} gives a rather surprising result. The right-hand side of the limit does not depend on the correlation among arrivals to different queues. In other words, in the heavy-traffic limit, these linear combinations of the queue lengths behave as if the arrival processes were independent if SSC is full-dimensional.

The proof of \Cref{gen.switch.cor.full.ssc} follows immediately from \Cref{gs.thm:bounds} because, if the cone $\cK$ is full-dimensional, then the subspace $\cH=\bR^n$ and, therefore, the projection matrix on $\cH$ satisfies $H=\bI$.

In the rest of this section, we present examples of SPNs that experience full-dimensional SSC.

\subsubsection{Parallel-server system.}\label{subsec:parallel-server-system}

\begin{figure}
	\centering
	\caption{Diagrams of examples of parallel-server systems. The dotted lines represent the compatibility between job-types and servers. \label{pic:parallel}}
	\vspace{0.1in}
	\begin{subfigure}[t]{0.3\linewidth}
		\centering 
		\includegraphics[width=\linewidth]{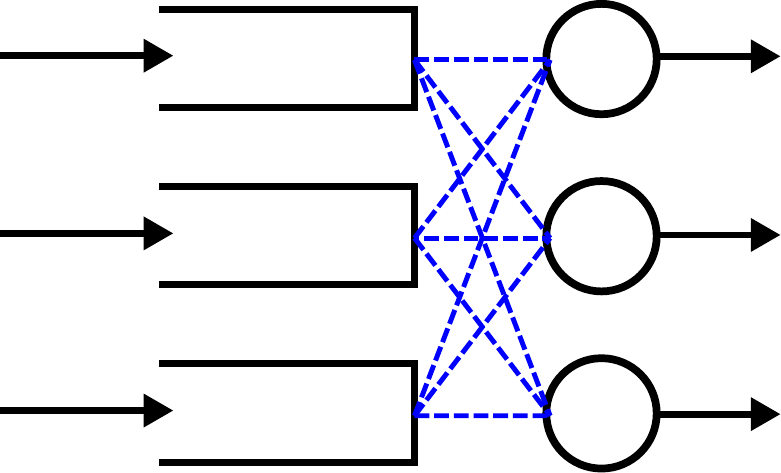}
		\caption{Fully flexible servers.\label{pic:parallel-general}}
	\end{subfigure}
	\hspace{.03\linewidth}
	\begin{subfigure}[t]{0.3\linewidth}
		\centering 
		\includegraphics[width=\linewidth]{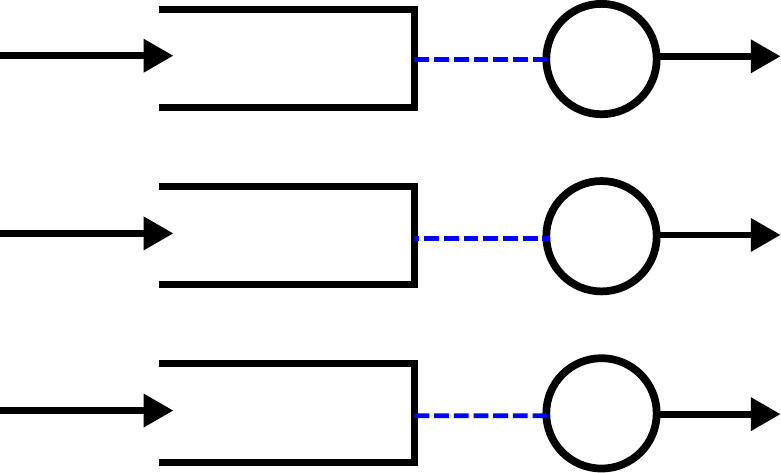}
		\caption{Dedicated servers.\label{pic:parallel-dedicated}}
	\end{subfigure}
	\hspace{.03\linewidth}
	\begin{subfigure}[t]{0.3\linewidth}
		\centering 
		\includegraphics[width=\linewidth]{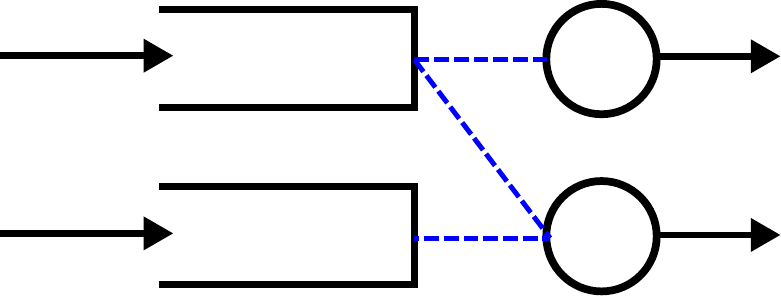}
		\caption{$\cN$-system.\label{pic:N-system}}
	\end{subfigure}
\end{figure}

Consider a parallel-server system as follows. There are $n$ types of jobs that arrive according to arrival processes as described in \Cref{sec:gen.switch.model}. Each job type can be processed by a subset of servers, and these subsets are modeled by a compatibility graph. In \Cref{pic:parallel} we present three examples of parallel-server systems, where the dotted lines represent the compatibility of the job types with the servers. In \Cref{pic:parallel-general}, all jobs can be served by all servers (fully flexible system), in \Cref{pic:parallel-dedicated} each job can be processed by only one server (dedicated system), and in \Cref{pic:N-system}, the jobs from the first queue can be processed by any server and the jobs from the second queue can only be processed by the second server ($\cN$-system, to be studied in \Cref{subsec:Nsystem}). 
The parallel-server systems (also called process flexibility) have received plenty of attention in the literature. For example, see the work by  \citeauthor{BellWill01-Nsystem} \cite{BellWill01-Nsystem}, \citeauthor{GarMan00-Nsystem-notes} \cite{GarMan00-Nsystem-notes}, \citeauthor{har_state_space} \cite{har_state_space}, \citeauthor{Williams_CRP} \cite{Williams_CRP}, and \citeauthor{zhong2019_process_flexibility} \cite{zhong2019_process_flexibility}, and the survey paper by \citeauthor{williams_survey_SPN} \cite{williams_survey_SPN}. However, most of the prior work is under the CRP condition. In this section we show that the parallel-server system can be studied as an immediate application of \Cref{gs.thm:bounds}, regardless of the CRP condition being satisfied.

To model a parallel-server system as a generalized switch, we assume that the service rate offered by each server in each time slot is a random variable that may depend on the service rate of other servers, but it is independent of the arrival and queueing processes. The joint distribution of the offered service rates is known, and we assume its state space is finite. Hence, the joint distribution of the offered service can be modeled as the channel state, and the compatibility graph determines the feasible service rate vectors in each time slot. Since we need the set of feasible service rate vectors in each channel state to be finite, we only consider the maximal vectors and their projection on the coordinate axes. Once the offered service rates are observed, the scheduler follows MaxWeight algorithm to decide which job types will be served and at which rate. We obtain the following result.

\begin{corollary}\label{gen.switch.cor.parallel-gral}
	Consider a set of parallel-server systems as described above, parametrized by $\epsilon$ as described in \Cref{gs.thm:bounds}. Suppose the capacity region $\cC$ has vertices that do not lie on the coordinate axes, and that $\vw$ is one of them. Let $\Sigma_B$ be as in \Cref{gs.thm:bounds} and $\vsigma_a$ be as in \Cref{gen.switch.cor.full.ssc}. Then,
	\begin{align*}
		\lim_{\epsilon\dto 0} \epsilon\E{\langle\vqbar^\peps,\vw\rangle}=& \dfrac{1}{2}\left(\left\|\vsigma_a \right\|^2 + \vone^T\left((C^TC)^{-1}\circ \Sigma_B \right)\vone \right).
	\end{align*}
\end{corollary}

The proof of \Cref{gen.switch.cor.parallel-gral} only requires modeling the parallel-server system as a generalized switch as we showed above, so we omit it. 

\begin{remark}
	In \Cref{gen.switch.cor.parallel-gral} we considered a vector $\vw$ in a vertex of the capacity region. However, \Cref{gs.thm:bounds} is immediately applicable for any $\vw$ in the boundary of the capacity region. Here we focused on a special case to illustrate the full-dimensional SSC result.
\end{remark}

Before finishing this section we present one of the simplest parallel-server systems to illustrate the result in \Cref{gen.switch.cor.parallel-gral}. Specifically, we work with a dedicated system, where every job type can be processed by exactly one server. A diagram with three job-types and three servers is presented in \Cref{pic:parallel-dedicated}.

Consider an SPN with $n$ servers, each with its own queue. Let $\left\{\hat{\vs}(k):k\in\bZ_+\right\}$ be a sequence of i.i.d. random vectors, such that $\hat{s}_i(k)$ is the potential service in queue $i$ in time slot $k$. Let $\vmu=\E{\hat{\vs}(1)}$ and $\Sigma_s$ be the covariance matrix of $\hat{\vs}(1)$. Suppose the vector $\hat{\vs}(1)$ has finite state space and that $\hat{s}_i(1)\leq \smax$ with probability 1 for all $i\in[n]$. Suppose $\min_{i\in[n]}\mu_i>0$.

The arrival process is defined as in \Cref{sec:gen.switch.model}, and we model heavy traffic as described there as well. Specifically, let $\epsilon\in(0,1)$ be the heavy-traffic parameter. Then, for each $\epsilon\in(0,1)$ and each $i\in[n]$, let the arrival process to the system be $\left\{\va^\peps(k):k\in\bZ_+\right\}$, which is a sequence of i.i.d. random vectors with mean $\vlambda^\peps=\E{\va^\peps(1)}=(1-\epsilon)\vmu$ and covariance matrix $\Sigma_a^\peps$. 

\begin{corollary}\label{gen.switch.cor.n-queues}
	Consider a set of dedicated parallel-server systems as described above, parametrized by $\epsilon\in(0,1)$ as described in \Cref{gs.thm:bounds}. Suppose $\lim_{\epsilon\dto 0}\Sigma_{a}^\peps=\Sigma_{a}$ component-wise. Let $\sigma_{a_i}^2=\left(\Sigma_a\right)_{i,i}$ and $\sigma_{s_i}^2=\left(\Sigma_s\right)_{i,i}$ for each $i\in[n]$. Then,
	\begin{align*}
		\lim_{\epsilon\dto 0}\epsilon\E{\sum_{i=1}^n \mu_i\qbar^\peps_i} = \dfrac{1}{2} \sum_{i=1}^n\left( \sigma_{a_i}^2 + \sigma_{s_i}^2 \right).
	\end{align*}
\end{corollary}

From the discussion after \Cref{gen.switch.cor.parallel-gral}, we expected that the correlation among the arrival processes would not be part of the right-hand side of the limit. However, observe that the correlation among the service processes does not appear in the answer either. Then, even though the arrival and potential service processes are correlated among queues, the mean linear combination of queue lengths behaves as if the queues were independent. Moreover, \Cref{gen.switch.cor.n-queues} recovers Kingman's bound. We present the proof of \Cref{gen.switch.cor.n-queues} below.
\proof{Proof of Corollary \ref{gen.switch.cor.n-queues}.}

The capacity region of this queueing system is
\begin{align*}
	\cC = \left\{\vx\in\bR_+^n: x_i\leq \mu_i \;,\, i\in[n] \right\}.
\end{align*}
To write it in the form of \eqref{gs.eq.cap.reg.pol}, we set $L=n$, and for each $i\in[n]$ we set $\vc^{(i)}=\ve^{(i)}$ and $b^{(i)}=\mu_i$. Therefore, the matrix $C$ is the identity matrix, which implies that the projection matrix $H$ is also the identity matrix.

Let $P=[n]$. Then, $\cap_{\ell \in P} \cF^\pl=\left\{\vmu\right\}$, and the left-hand side of \eqref{eq.gen.switch.thm.limit} yields
\begin{align*}
	\lim_{\epsilon\dto 0}\epsilon \E{\sum_{i=1}^n \mu_i\qbar^\peps_i}
\end{align*}

Since the projection matrix satisfies $H=\bI$, the first term on the right-hand side of \eqref{eq.gen.switch.thm.limit} yields
\begin{align*}
	\dfrac{1}{2}\vone^T\left(H\circ \Sigma_a\right)\vone =& \dfrac{1}{2}\vone^T\left(\bI \circ \Sigma_a\right)\vone = \dfrac{1}{2}\sum_{i=1}^n \sigma_{a_i}^2.
\end{align*}

To compute the second term of the right-hand side of \eqref{eq.gen.switch.thm.limit}, we consider the following interpretation of the channel state. Let $\cM$ be an enumeration of the elements of the state space of $\hat{\vs}(1)$, and $\vs^\pm$ be its $m\tth$ element for each $m\in\cM$. For each $m\in\cM$, let the set of feasible service rate vectors in channel state $m$ be
\begin{align*}
	\cS^\pm = \left\{\vs^\pm \right\}\cup \left\{\vs^\pm-s^\pm_i\ve^{(i)}:\; i\in[n] \right\},
\end{align*}
i.e., the set $\cS^\pm$ contains $\vs^\pm$ and its projection on the coordinate axes. We assume that MaxWeight breaks ties by choosing maximal schedules. Then, if the channel state is $m$ then the service rates vector is always $\vs^\pm$. With this assumption we lose some generality because arrivals occur after deciding the optimal schedule. However, we are interested in heavy-traffic analysis so this slight loss of generality does not affect our result. Then, the probability mass function of the channel state $\vpsi$ satisfies $\psi_m\defn \Prob{\hat{\vs}(1)=\vs^\pm}$ for each $m\in\cM$.

By definition of $\bml$ in \eqref{gs.eq.cl.weighted.def} and by definition of the sets $\cS^\pm$ and the vectors $\vcl$ above, we obtain that for each $\ell\in[n]$ we have
\begin{align*}
	\bml =& \langle\vcl,\vs^\pm\rangle = \langle\ve^\pl,\vs^\pl\rangle = s^\pm_\ell.
\end{align*}
Then, for each $\ell\in[n]$ the random variable $B_\ell(1)$ is such that $\Prob{B_\ell(1)=s_\ell^\pm}=\psi_m$ and $\E{B_\ell(1)}=\mu_\ell$. Therefore, the vectors $\left(B_1(1),\ldots,B_n(1) \right)$ and $\hat{\vs}(1)$ have the same distribution. Hence, $\left(\Sigma_B\right)_{i,j}=\Cov{\hat{\vs}_i,\hat{\vs}_j}$, and the second term in the right-hand side of \eqref{eq.gen.switch.thm.limit} becomes
\begin{align*}
	\dfrac{1}{2}\vone^T\left((C^TC)^{-1}\circ \Sigma_B \right)\vone &\stackrel{(a)}{=} \dfrac{1}{2}\vone^T\left(\bI \circ \Sigma_B \right)\vone	\stackrel{(b)}{=} \sum_{i=1}^n \sigma_{s_i}^2,
\end{align*}
where $(a)$ holds because $C=\bI$; and $(b)$ holds by definition of Hadamard's product and because the diagonal of $\Sigma_B$ contains the variance of $\hat{s}_i(1)$ for each $i\in[n]$. \Halmos
\endproof

\subsubsection{$\cN$-system.}\label{subsec:Nsystem}

The $\cN$-system model is a parallel-server system with two servers and two job types. One of the servers exclusively serves the jobs type 1, and the other server can process both. A diagram of the $\cN$-system is presented in \Cref{pic:N-system}. According to \citeauthor{GhamWard13-Nsystem} \cite{GhamWard13-Nsystem}, ``the $\cN$-system is one of the simplest parallel server system models that retains much of the complexity inherent in more general models''. Consequently, it has received plenty of attention over the years and there is vast literature that only focuses on its performance under the CRP condition. Examples can be found in the work by \citeauthor{BellWill01-Nsystem} \cite{BellWill01-Nsystem}, \citeauthor{GarMan00-Nsystem-notes} \cite{GarMan00-Nsystem-notes}, \citeauthor{har_state_space} \cite{har_state_space}, and \citeauthor{zhong2019_process_flexibility} \cite{zhong2019_process_flexibility}. \Cref{gs.thm:bounds} is immediately applicable to this system, and gives information about the mean queue lengths in both, the CRP and non-CRP cases. In this section we focus on the non-CRP case.

Let the arrival processes be as described in \Cref{sec:gen.switch.model} and suppose that each server processes jobs at rate 1. Then, the capacity region of this system is $\cC=\left\{ \vx \in\bR_+^2: x_1\leq 1,\, x_2\leq 1 \right\}$. We consider the heavy-traffic parametrization $\vlambda^\peps = (1-\epsilon)\vone$, for $\epsilon\in(0,1)$. Then, as $\epsilon\dto 0$, the arrival rate vector approaches a vertex of the capacity region and, hence, the $\cN$-system experiences full-dimensional SSC. We now present the result.

\begin{corollary}\label{cor:Nsystem}
	Consider a set of $\cN$-systems parametrized by $\epsilon\in(0,1)$,  as described above. Let $\vsigma_a$ be as in \Cref{gen.switch.cor.parallel-gral}. Then,
	\begin{align*}
		\lim_{\epsilon\dto 0} \epsilon \E{\qbar^\peps_1 + \qbar_2^\peps} = \dfrac{\sigma_{a_1}^2 + \sigma_{a_2}^2}{2}.
	\end{align*}
\end{corollary}

The proof is an immediate application of \Cref{gen.switch.cor.full.ssc}, so we omit it.

\subsubsection{Ad hoc wireless network.}\label{subsec:adhoc-wireless-network}

\begin{figure}
	\centering
	\caption{Diagram of ad hoc wireless networks.\label{pic:adhoc}}
	\begin{subfigure}[t]{0.3\linewidth}
		\centering 
		\includegraphics[width=1.3in]{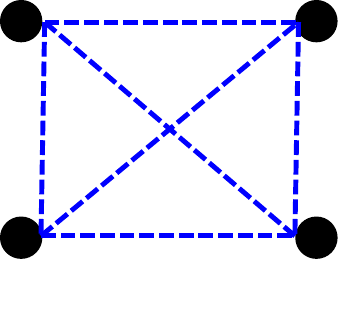}
		\caption{Example of an interference graph for an ad hoc wireless network with four links.\label{pic:adhoc-interference}}
	\end{subfigure}
	\hspace{.5in}
	\begin{subfigure}[t]{0.3\linewidth}
		\centering
		\includegraphics[width=\linewidth]{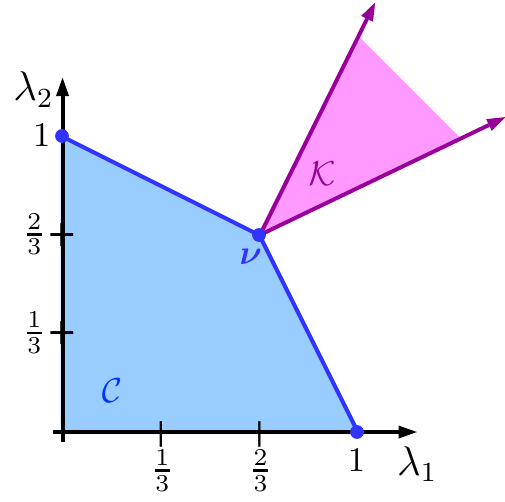}
		\caption{Capacity region and cone for SPN in \Cref{subsec:adhoc-wireless-network}.\label{pic:adhoc-C}}
	\end{subfigure}
\end{figure}

An ad hoc wireless network is composed by a set of nodes with no infrastructure for central coordination, and packets are transmitted between nodes (a transmitter and a receiver) if there is a link. The links interfere with each other and, therefore, not all of them can be active at the same time. These interference constraints are frequently represented with a graph, where the vertices represent links and an edge between two links represents interference. In \Cref{pic:adhoc-interference} we present an example of the interference graph of an ad hoc wireless network with four links, where all links interfere with each other. The packets to be transmitted arrive to each of the links and wait in line until they can be processed. This model has been studied in a long line of literature, including but not limited to the work by  \citeauthor{DimWal_06} \cite{DimWal_06}, \citeauthor{sivhajsri_11} \cite{sivhajsri_11}, \citeauthor{atilla} \cite{atilla}, and \citeauthor{KanWanJarYin_wireless-networks} \cite{KanWanJarYin_wireless-networks}, but in most of the cases the focus is on studying stability or optimality of the scheduling policy. Here we provide the heavy-traffic limit of linear combinations of the queue lengths under MaxWeight algorithm. A particular case of our result are the results obtained by \citeauthor{atilla} \cite{atilla}.

An ad hoc wireless network can be modeled as a generalized switch with fixed channel state. Then, \Cref{gs.thm:bounds} can be immediately applied. In this section we provide an example of an ad hoc wireless network that experiences full-dimensional SSC. We focus on a network with two links to illustrate the geometry of the capacity region and the cone where SSC occurs, but similar work can be done for larger networks.

Let $\sigma_1^2\defn \Var{\abar_1^\peps}$, $\sigma_2^2\defn \Var{\abar_2^\peps}$, and $\varphi\defn \Cov{\abar_1^\peps,\abar_2^\peps}$, where these three parameters do not depend on $\epsilon$. Suppose the set of feasible service rate vectors is $\cS = \left\{(1,0),\, (0,1),\,\left(\tfrac{2}{3},\tfrac{2}{3}\right)\right\}$. Then, the capacity region is $\cC=\left\{\vx\in\bR_+^2: x_1+2x_2\leq 2,\, 2x_1+x_2\leq 2 \right\}$. 
Applying \Cref{gs.thm:bounds} we obtain the following corollary.

\begin{corollary}\label{cor:adhoc}
	Consider an ad hoc wireless network as described above. Then,
	\begin{align*}
		\lim_{\epsilon\dto 0}\epsilon\E{\qbar_1^\peps+\qbar_2^\peps}= \dfrac{3}{4}\left(\sigma_1^2+\sigma_2^2\right)
	\end{align*}
\end{corollary}

In the proof of \Cref{cor:adhoc} we take the heavy-traffic limit as the vector of arrival rate approaches the point $\vnu=\frac{2}{3}(1,1)$ in the boundary of $\cC$. The proof is simple, so we omit it. In \Cref{pic:adhoc-C} we plot the capacity region, the point $\vnu$ and the cone $\cK$ where SSC occurs. Observe that the cone $\cK$ is two-dimensional and, therefore, this ad hoc wireless network experiences full-dimensional SSC.

\begin{remark}
	The input-queued switch can be modeled similarly to an ad hoc wireless network. However, the input-queued switch cannot experience full-dimensional SSC since all the vertices of its capacity region are on the coordinate axes. In other words, all the vertices in the capacity region of the input-queued switch require the arrival rate to (at least) one queue to be zero. This is equivalent to considering a queueing system where the zero-arrival rate queue does not exist, which already has a lower-dimensional state space.
\end{remark}

\section{Individual queue lengths and higher moments in the input-queued switch.}\label{sec:system.of.equations}
In this section we show that the drift method with polynomial test functions does not provide all the information that is necessary to compute the moments of all the linear combinations of the scaled queue lengths in systems that do not satisfy the CRP condition. We do this by presenting an alternate view of the drift method.

In the proof of \Cref{gs.thm:bounds} we use $V(\vq)=\left\|\vqparh\right\|^2$ as test function to obtain bounds on certain linear combinations of the queue lengths in a generalized switch. This choice of test function was first proposed by \citeauthor{MagSri_SSY16_Switch} \cite{MagSri_SSY16_Switch}, and the main reason to use it is that the term $\cT_4$ consisting of the `$qu$' terms (i.e., cross terms between the queue length and the unused service) converges to zero in the heavy-traffic limit. All of queueing theory in some sense is to get a handle on the unused service terms, and the drift method handles these terms by making sure that they `cancel out' in heavy traffic, using SSC and our choice of the test function. In this section, instead of trying to cancel out the `$qu$' terms, we consider them as unknowns and try to solve for them along with the mean queue lengths. We will see that this is impossible even if we use all possible quadratic test functions.

For simplicity of exposition, we present this result in the context of an input-queued switch, which is one of the simplest queueing systems that experience multidimensional SSC and it is a special case of the generalized switch, as shown in \Cref{subsec:switch}. The organization of this section is as follows. In \Cref{sec:soe.thm.syst.eqns} we present the main result, in \Cref{sec:soe.bounds} we use this result to compute bounds on the first moment of linear combinations of the scaled queue lengths and in \Cref{sec:soe.generalization} we discuss how to generalize this approach to other queueing systems that experience multidimensional SSC. 

\subsection{System of equations to compute linear combinations of the first moment of scaled queue lengths.}\label{sec:soe.thm.syst.eqns}

In this section we prove that the drift method with polynomial test functions is not sufficient to compute all the linear combinations of the first moment of the scaled queue lengths in queueing systems that do not satisfy the CRP condition. Specifically, we show that the use of polynomial test functions yields an under-determined system of equations.

In the drift method, one of the key challenges is to get a handle on the unused service. In general, when one sets to zero the drift of a polynomial test function in steady state, terms of the form $q_i(k+1)u_j(k)$ arise. The idea is to use a test function that captures the geometry of SSC so that we can show that all these cross terms are small. Therefore, the choice of the test function is important, and the region into which SSC happens must be used in this choice. The quadratic test function, $V(\vq)=\|\vqparh\|^2$ has been successfully used by \citeauthor{atilla} \cite{atilla}, \citeauthor{MagSri_SSY16_Switch} \cite{MagSri_SSY16_Switch}, \citeauthor{QUESTA_switch} \cite{QUESTA_switch} and \citeauthor{WeinaBandwidthJournal} \cite{WeinaBandwidthJournal} to obtain the mean sum of the queue lengths, similarly to Theorem \ref{gs.thm:bounds}. Typically one uses polynomial test functions of degree $(m+1)$ to get bounds on the expected value of the $m\tth$ power of the queue lengths. Therefore, in order to obtain bounds on the mean queue lengths, one must use quadratic test functions. In order to get all the linear combinations of the queue lengths, one can search through all the quadratic test functions, and this is equivalent to searching through all the quadratic monomials. The following theorem presents the result of using all the quadratic monomial test functions.

For ease of exposition, in this section we prove our result in the case of $N=2$ and independent arrivals, i.e., in the case of a  $2\times 2$ input-queued switch with independent arrivals. We present generalizations to this result in Appendix \ref{app:switch.generalizations}. Specifically, we present the case of a $2\times 2$ input-queued switch with correlated arrivals in Appendix \ref{app:switch.2x2.correlated}, and the case of an $N\times N$ input-queued switch with independent arrivals in Appendix \ref{app:switch.nxn}. The latter result can be easily generalized to the case of correlated arrivals, but we do not present the result here for ease of exposition.

\begin{theorem}\label{s.theorem:2switch}
	Consider a set of $\,2\times 2$ input-queued switches operating under MaxWeight, indexed by $\epsilon\in(0,1)$ as described in Corollary \ref{cor.nxn.switch.indep}. Let $(\Sigma_a^\peps)_{i,i}=\sigma_{a_i}^\peps$ and suppose $\lim_{\epsilon\dto 0}\sigma_{a_i}^\peps=\sigma_{a_i}$ for all $i\in[4]$. Then, the following system of equations is satisfied
	
	\begin{align}
	&\begin{aligned}\label{s.2.eq.1}
	&\lim_{\epsilon\dto 0} \epsilon\E{\qbar_1} \\
	&= \dfrac{9\sigma_{a_1}^2 + \sigma_{a_2}^2 + \sigma_{a_3}^2 + \sigma_{a_4}^2}{16} +\dfrac{1}{2}\lim_{\epsilon\dto 0}\E{\qbar_1^+\left(\ubar_2+\ubar_3 \right)}  -\dfrac{1}{2}\lim_{\epsilon\dto 0}\E{\left(\qbar_2^+ + \qbar_3^+ \right)\ubar_4}
	\end{aligned}\\[3pt]
	&\begin{aligned}\label{s.2.eq.2}
	&\lim_{\epsilon\dto 0} \epsilon \E{\qbar_2} \\
	&= \dfrac{\sigma_{a_1}^2 + 9\sigma_{a_2}^2 + \sigma_{a_3}^2 + \sigma_{a_4}^2}{16} + \dfrac{1}{2} \lim_{\epsilon\dto 0}\E{\qbar_2^+\left(\ubar_1-\ubar_3+\ubar_4\right)}
	\end{aligned}\\[3pt]
	&\begin{aligned}\label{s.2.eq.3}
	&\lim_{\epsilon\dto 0} \epsilon\E{\qbar_3}\\
	&= \dfrac{\sigma_{a_1}^2 + \sigma_{a_2}^2 + 9\sigma_{a_3}^2 + \sigma_{a_4}^2}{16} + \dfrac{1}{2}\lim_{\epsilon\dto 0} \E{\qbar_3^+\left(\ubar_1 -\ubar_2+\ubar_4\right)}
	\end{aligned}\\[3pt]
	& \begin{aligned}\label{s.2.eq.4}
	&\lim_{\epsilon\dto 0} \epsilon\E{\qbar_1+\qbar_2} \\
	&= \dfrac{3\sigma_{a_1}^2 + 3\sigma_{a_2}^2 - \sigma_{a_3}^2 - \sigma_{a_4}^2}{8} + \dfrac{1}{2}\lim_{\epsilon\dto 0}\E{\qbar_1^+ \left(3\ubar_2-\ubar_3\right)} + \dfrac{1}{2}\lim_{\epsilon\dto 0}\E{\qbar_2^+ \left(3\ubar_1+\ubar_3\right)} + \dfrac{1}{2}\lim_{\epsilon\dto 0}\E{\qbar_3^+ \ubar_4}
	\end{aligned} \\[3pt]
	& \begin{aligned}\label{s.2.eq.5}
	&\lim_{\epsilon\dto 0} \epsilon\E{\qbar_1+ \qbar_3} \\
	&= \dfrac{3 \sigma_{a_1}^2 - \sigma_{a_2}^2 + 3\sigma_{a_3}^2 - \sigma_{a_4}^2}{8}+ \dfrac{1}{2}\lim_{\epsilon\dto 0}\E{\qbar_1^+ \left(-\ubar_2+3\ubar_3\right)} + \dfrac{1}{2}\lim_{\epsilon\dto 0}\E{\qbar_2^+ \ubar_4} + \dfrac{1}{2}\lim_{\epsilon\dto 0}\E{\qbar_3^+ \left(3\ubar_1+ \ubar_2\right)}
	\end{aligned} \\[3pt]
	& \begin{aligned}\label{s.2.eq.6}
	&\lim_{\epsilon\dto 0} \epsilon\E{\qbar_2+ \qbar_3} \\
	&= \dfrac{\sigma_{a_1}^2 - 3\sigma_{a_2}^2 - 3\sigma_{a_3}^2 + \sigma_{a_4}^2}{8} + \dfrac{1}{2}\lim_{\epsilon\dto 0} \E{\qbar_2^+ \left(\ubar_1+3\ubar_3+\ubar_4\right)} + \dfrac{1}{2}\lim_{\epsilon\dto 0}\E{\qbar_3^+ \left(\ubar_1+3\ubar_2+\ubar_4\right)} ,
	\end{aligned}
	\end{align}
	where we omitted the dependence on $\epsilon$ of the variables for ease of exposition.
\end{theorem}

The proof of Theorem \ref{s.theorem:2switch} is presented in Section \ref{subsec:proof.switch}. Observe that in Theorem \ref{s.theorem:2switch} we have system of 6 equations and 11 variables, where the variables are
\begin{align*}
& \lim_{\epsilon\dto 0}\epsilon \E{\qbar_1}\;,\;\lim_{\epsilon\dto 0}\epsilon \E{\qbar_2}\;,\;\lim_{\epsilon\dto 0}\epsilon \E{\qbar_3}, \\
& \lim_{\epsilon\dto 0}\E{\qbar_1^+ \ubar_2}\;,\; \lim_{\epsilon\dto 0}\E{\qbar_1^+ \ubar_3}, \\
& \lim_{\epsilon\dto 0}\E{\qbar_2^+ \ubar_1}\;,\; \lim_{\epsilon\dto 0}\E{\qbar_2^+ \ubar_3}\;,\; \lim_{\epsilon\dto 0}\E{\qbar_2^+ \ubar_4}, \\
& \lim_{\epsilon\dto 0}\E{\qbar_3^+ \ubar_1}\;,\; \lim_{\epsilon\dto 0}\E{\qbar_3^+ \ubar_2} \;,\; \lim_{\epsilon\dto 0}\E{\qbar_3^+ \ubar_4}.
\end{align*}
Therefore, it cannot be solved uniquely. However, a specific linear combination of the scaled queue lengths can be obtained, as shown in the next Corollary.

\begin{corollary}\label{s.corollary}
	Consider a set of $\,2\times 2$ input-queued switches as described in Theorem \ref{s.theorem:2switch}. Then,
	\begin{align*}
	\lim_{\epsilon\dto 0}\epsilon \E{\qbar_2 + \qbar_3 }= \dfrac{3}{8}\left(\sigma_{a_1}^2+ \sigma_{a_2}^2+ \sigma_{a_3}^2 + \sigma_{a_4}^2\right)
	\end{align*}
\end{corollary}
\proof{Proof of Corollary \ref{s.corollary}.}
Consider the following linear combination of the equations in Theorem \ref{s.theorem:2switch}:
\begin{align*}
\eqref{s.2.eq.1} +\eqref{s.2.eq.2} +\eqref{s.2.eq.3} -\frac{1}{2}
\eqref{s.2.eq.4} -\frac{1}{2}\eqref{s.2.eq.5} +\frac{1}{2}\eqref{s.2.eq.6}.
\end{align*}
Then, reorganizing terms we obtain the result.
\Halmos
\endproof

Corollary \ref{s.corollary} can be also obtained as a consequence of Corollary \ref{cor.nxn.switch.indep} in the following way.

\proof{Alternative proof of Corollary \ref{s.corollary}.}
From Corollary \ref{cor.nxn.switch.indep} for $N=2$ we know
\begin{align}
& \lim_{\epsilon\dto 0} \epsilon\E{\qbar_1 + \qbar_2 + \qbar_3 + \qbar_4} = \dfrac{3}{4}\left(\sigma_{a_1}^2+ \sigma_{a_2}^2+ \sigma_{a_3}^2 + \sigma_{a_4}^2\right).
\end{align}

From SSC as proved in Proposition \ref{gs.prop:gen.switch.SSC} and by definition of the cone $\cK$ in \eqref{eq.switch.cone} we also know that for all $i\in[4]$ we have
\begin{align*}
\lim_{\epsilon\dto 0}\epsilon \E{\qbar_{\parallel\cH i}} = \lim_{\epsilon\dto 0}\epsilon \E{\qbar_i},
\end{align*}
where $\qbar_{\parallel\cH i}$ is the $i\tth$ element of $\vqbarparh$. Also, one interpretation of the cone $\cK$ presented by \citeauthor{MagSri_SSY16_Switch} \cite{MagSri_SSY16_Switch} is that for each vector in $\cK$, all schedules have the same weight in MaxWeight algorithm. This can be easily verified by definition of the cone $\cK$ in \eqref{eq.switch.cone}. Then,
\begin{align}\label{s.2.eq.all.schedules.same}
\qbar_{\parallel\cH 1}+\qbar_{\parallel\cH 4} = \qbar_{\parallel\cH 2}+\qbar_{\parallel\cH 3}.
\end{align}
Putting everything together we obtain the result in Corollary \ref{s.corollary}.
\Halmos
\endproof

In Theorem \ref{s.theorem:2switch} we prove that setting to zero the drift of all monomials of degree 2 leads to a system of 6 equations in 11 variables. Therefore, the solution is not unique. However, \citeauthor{MagSri_SSY16_Switch} \cite{MagSri_SSY16_Switch} and \citeauthor{QUESTA_switch} \cite{QUESTA_switch} obtain the limit of specific linear combinations of the scaled queue lengths. These linear combinations can be obtained because some of the variables cancel out, as shown in the first proof of \Cref{s.corollary}. However, to obtain other linear combinations of the expected heavy-traffic scaled queue lengths we need to actually work with all the variables of the system of equations. Therefore, we need additional equations. 

To better understand this argument, consider a tandem queue system with memoryless interarrival and service times in any (not necessarily heavy) traffic. We know that the steady-state joint distribution is product of two geometrics, and can be obtained using reversibility arguments. Using the drift approach described above, we get 3 equations in 4 unknowns. However, in addition to the drift arguments, if we use reversibility to separately prove that the queues are independent in steady state and impose it as an additional condition, we can solve for all the unknowns.

\subsection{Bounds on linear combinations of the scaled queue lengths in heavy-traffic.}\label{sec:soe.bounds}

In Section \ref{sec:soe.thm.syst.eqns} we presented a linear system of equations that the vector of queue lengths must satisfy in heavy-traffic. In this section we use this system of equations to obtain bounds on linear combinations of the expected scaled queue lengths in heavy traffic. A similar approach was studied by  \citeauthor{kumar_kumar_lineqns} \cite{kumar_kumar_lineqns} and \citeauthor{bertsimas_paschalidis_Tstitsiklis_optimization} \cite{bertsimas_paschalidis_Tstitsiklis_optimization}, where an under-determined set of linear systems of equations was obtained and linear programming was used to obtain bounds. However, the focus in those papers was on queueing networks under fixed arrival and service rates, as opposed to the heavy-traffic analysis in the current paper.

In the next theorem we provide an upper and a lower bound for the heavy-traffic limit of the expected value of any linear combination of the queue lengths in a $2\times 2$ input-queued switch.

\begin{theorem}\label{s.2.thm.lp}
	Consider the equations
	\begin{align}
	&v_1-\dfrac{w_1-w_2+w_5+w_8}{2}
	= \dfrac{9\sigma_{a_1}^2 + \sigma_{a_2}^2 + \sigma_{a_3}^2+\sigma_{a_4}^2}{16} \label{s.2.lp.eq.1}\\[3pt]
	&v_2-\dfrac{w_3+w_4-w_5}{2}
	= \dfrac{\sigma_{a_1}^2 + 9\sigma_{a_2}^2 + \sigma_{a_3}^2 + \sigma_{a_4}^2}{16} \label{s.2.lp.eq.2} \\[3pt]
	&v_3 - \dfrac{w_6+w_7-w_8}{2}
	= \dfrac{\sigma_{a_1}^2 +\sigma_{a_2}^2 + 9\sigma_{a_3}^2 + \sigma_{a_4}^2}{16} \label{s.2.lp.eq.3} \\[3pt]
	&\begin{aligned}\label{s.2.lp.eq.4}
	&v_1+v_2 - \dfrac{3w_+w_2-3w_3-w_4-w_8}{2} = \dfrac{3\sigma_{a_1}^2 + 3\sigma_{a_2}^2 -\sigma_{a_3}^2 - \sigma_{a_4}^2}{8}
	\end{aligned}  \\[3pt]
	&\begin{aligned}\label{s.2.lp.eq.5}
	&v_1+v_3 +\dfrac{w_1-3w_2-w_5-3w_6-w_7}{2} = \dfrac{3\sigma_{a_1}^2 - \sigma_{a_2}^2 + 3\sigma_{a_3}^2 - \sigma_{a_4}^2}{8}
	\end{aligned}\\[3pt]
	&\begin{aligned}\label{s.2.lp.eq.6}
	&v_2+v_3 - \dfrac{w_3-3w_4-w_5-w_6-3w_7-w_8}{2} = \dfrac{\sigma_{a_1}^2 - 3\sigma_{a_2}^2 - 3\sigma_{a_3}^2 + \sigma_{a_4}^2}{8}
	\end{aligned} \\
	& -v_1+v_2+v_3 \geq 0 \label{s.2.lp.ineq.1} \\
	& -w_1+w_7\geq 0 \label{s.2.lp.ineq.2} \\
	& -w_2 + w_4\geq 0 \label{s.2.lp.ineq.3}
	\end{align}
	and define $\mathcal{P}\defn\left\{(\vv,\vw)\in\bR^{3}_+\times \bR^8_+: \text{Equations \eqref{s.2.lp.eq.1}-\eqref{s.2.lp.ineq.3} are satisfied} \right\}$. For $\valpha\in\bR^3$, define
	\begin{align*}
	& \underline{f}(\valpha)\defn \min\left\{\langle\valpha,\vv\rangle:\; \exists \vw\text{ such that }(\vv,\vw)\in\mathcal{P} \right\}\\
	\text{and}\quad& \overline{f}(\valpha)\defn \max\left\{\langle\valpha,\vv\rangle:\; \exists \vw\text{ such that }(\vv,\vw)\in\mathcal{P} \right\}.
	\end{align*}
	Then,
	\begin{align}\label{s.2.lp.eq.bounds}
	\underline{f}(\valpha)\leq \lim_{\epsilon\dto 0}\epsilon \E{\langle \valpha,\vqbar^\peps\rangle}\leq \overline{f}(\valpha),
	\end{align}
	where $\epsilon$ and $\vqbar^\peps$ are defined as in Theorem \ref{s.theorem:2switch}. Furthermore, for any $B\in\bR_+$
	\begin{align}\label{s.2.lp.eq.Markov}
	\Prob{\lim_{\epsilon\dto 0}\epsilon\langle\valpha,\vqbar^\peps\rangle \geq B}\leq \dfrac{\overline{f}(\valpha)}{B}.
	\end{align}
\end{theorem}

\proof{Proof of Theorem \ref{s.2.thm.lp}.}
For ease of exposition we omit the dependence on $\epsilon$ of the variables. Let
\begin{align*}
& v_1=\lim_{\epsilon\dto 0}\epsilon \E{\qbar_1}\;,\; v_2=\lim_{\epsilon\dto 0}\epsilon \E{\qbar_2}\;,\; v_3=\lim_{\epsilon\dto 0}\epsilon \E{\qbar_3} ,\\
& w_1=\lim_{\epsilon\dto 0}\E{\qbar_1^+ \ubar_2}\;,\; w_2=\lim_{\epsilon\dto 0}\E{\qbar_1^+ \ubar_3}, \\
& w_3=\lim_{\epsilon\dto 0}\E{\qbar_2^+ \ubar_1}\;,\; w_4=\lim_{\epsilon\dto 0}\E{\qbar_2^+ \ubar_3}\;,\\
& w_5=\lim_{\epsilon\dto 0}\E{\qbar_2^+ \ubar_4} \;,\; w_6=\lim_{\epsilon\dto 0}\E{\qbar_3^+ \ubar_1}\;,\\
& w_7=\lim_{\epsilon\dto 0}\E{\qbar_3^+ \ubar_2} \;,\; w_8=\lim_{\epsilon\dto 0}\E{\qbar_3^+ \ubar_4}.
\end{align*}

Then, the proof of \eqref{s.2.lp.eq.bounds} follows from Theorem \ref{s.theorem:2switch} because the set $\mathcal{P}$ represents the system of equations presented there together with non-negativity constraints for all the variables. In particular, inequalities \eqref{s.2.lp.ineq.1}-\eqref{s.2.lp.ineq.3} represent non-negativity constraints associated to $\qbar_4$. These must be considered because, even though $\qbar_4$ does not appear in the system of equations explicitly, there are underlying constraints of the system related to $\qbar_4$ that affect its performance. Specifically, using \eqref{s.2.eq.all.schedules.same} and the definition of the variables above, we obtain that the inequalities
\begin{align*}
	\lim_{\epsilon\dto 0}\epsilon\E{\qbar_4}\geq 0\;,\quad \lim_{\epsilon\dto 0}\E{\qbar_4^+ \ubar_i}\geq 0\; \forall i\in\{1,2,3\}
\end{align*}
can be rewritten as \eqref{s.2.lp.ineq.1}, \eqref{s.2.lp.ineq.2}, \eqref{s.2.lp.ineq.3} and $w_3+w_6\geq 0$ but the last inequality is implied by $w_3\geq 0$ and $w_6\geq 0$, so we do not write it in the definition of $\mathcal{P}$.

Also, from Markov's inequality we know
\begin{align*}
\Prob{\lim_{\epsilon\dto 0}\epsilon\langle\valpha,\vqbar^\peps\rangle \geq B}\leq \dfrac{\lim_{\epsilon\dto 0}\epsilon \E{\langle\valpha,\vqbar^\peps\rangle}}{B}\leq \dfrac{\overline{f}(\valpha)}{B},
\end{align*}
where the last inequality holds by \eqref{s.2.lp.eq.bounds}.
\Halmos
\endproof

Theorem \ref{s.2.thm.lp} gives explicit bounds for all linear combinations of the expected scaled queue lengths. Similar linear programs can be written to obtain bounds on higher
moments, and consequently tighter tail probabilities.

In the rest of this section we present numerical results to compare the bounds that we obtain from the linear program presented in Theorem \ref{s.2.thm.lp} with the mean values that we obtain from simulation. We test four different objective functions, viz. $\lim_{\epsilon\dto 0}\epsilon \E{\qbar_i}$ for $i\in\{1,2,3\}$ and $\lim_{\epsilon\dto 0}\epsilon \E{\qbar_2+\qbar_3}$. We use the last function because in this case the system of equations has a unique solution, as shown in Corollary \ref{s.corollary}.

For simplicity, we assume tha t the arrivals to each queue are Bernoulli processes with mean $\lambda^\peps_i=\frac{1-\epsilon}{2}$ for all $i\in[4]$. We take $\epsilon\in \{0.01,0.05,0.1\}$ to evaluate the performance under different traffic intensities.

In the case of $\epsilon=0.01$ we ran the simulation for $5\times 10^8$ time slots and for the other values of $\epsilon$ we ran it for $2\times 10^8$ time slots. The reason is that, for smaller $\epsilon$, the system takes more time to reach steady state. In both cases we compute the mean value of the variables considering the last $2\times 10^4$ time slots. We present our results in Tables \ref{tab:s2.q2q3} and \ref{tab:s2.ind.q}. In Table \ref{tab:s2.q2q3} we present the right-hand side of the expression proved in Corollary \ref{s.corollary} and the results from the simulation. Specifically, we show the mean value of $\epsilon\left(\qbar_2^\peps+\qbar_3^\peps\right)$, a confidence interval constructed by adding and subtracting the standard deviation to the mean value, and the percentage error of the solution of the system of equations with respect to the simulation.

\begin{table}
	\TABLE
	{Numerical results for LP with objective function $\lim_{\epsilon\dto 0}\epsilon\E{\qbar_2^\peps+\qbar_3^\peps}$.\label{tab:s2.q2q3}}
	{\begin{tabular}{|c|c|c|c|c|}
			\hline
			$\epsilon$ & Solution to LP & Mean from simulation & Confidence interval & Error \\ \hline
			0.01  & 0.375 & 0.401 & (0.244, 0.557) & 6\% \\
			0.05  & 0.374 & 0.343 & (0.155, 0.531) & 9\% \\
			0.10   & 0.371 & 0.331 & (0.121, 0.540) & 12\% \\ \hline
		\end{tabular} }
	{}
\end{table}

Observe that in Table \ref{tab:s2.q2q3} the solution to the linear program is always in the confidence interval obtained from simulation, even in the case of $\epsilon=0.1$ which may be considered high for heavy traffic. The error is increasing with respect to $\epsilon$, but it is below 15\% in all cases.

In Table \ref{tab:s2.ind.q} we compute a lower and an upper bound to the mean individual queue lengths, and we compare these results with the mean value of $\epsilon\qbar_1$, $\epsilon\qbar_2$ and $\epsilon\qbar_3$ obtained from simulation. The reason to present only one optimal value for all the queue lengths is that solving the linear program presented in Theorem \ref{s.2.thm.lp} with objective function $\epsilon\E{\qbar_i^\peps}$ gives the same optimal value for all $i=1,2,3$, because of the symmetric arrival pattern.

\begin{table}
	\TABLE
	{Numerical results for individual queue lengths.\label{tab:s2.ind.q}}
	{ \begin{tabular}{|c|c|c|c|c|c|}
		\hline
		\multicolumn{1}{|c|}{\multirow{2}[2]{*}{Value of $\epsilon$}} & \multicolumn{1}{c|}{\multirow{2}[2]{*}{Minimum}} & \multicolumn{1}{c|}{\multirow{2}[2]{*}{Maximum}} & \multicolumn{3}{c|}{Simulation} \\
		&       &       & \multicolumn{1}{l|}{Mean $\epsilon \qbar_1$} & \multicolumn{1}{l|}{Mean $\epsilon \qbar_2$} & \multicolumn{1}{l|}{Mean $\epsilon \qbar_3$} \\
		\hline
		0.05  & 0.062 & 0.312 & 0.165 & 0.167 & 0.176 \\
		0.10   & 0.062 & 0.309 & 0.161 & 0.171 & 0.160 \\
		\hline
	\end{tabular} }
	{ }
\end{table}

Observe that for all the cases presented in Table \ref{tab:s2.ind.q}, the mean obtained by simulation is between the lower and upper bound obtained solving the LP. The bounds are not necessarily tight, but they can be computed very fast, as opposed to the mean values obtained from the simulation.

\subsection{Generalization to other queueing systems and higher moments.}\label{sec:soe.generalization}

In this section we focused on a $2\times 2$ input-queued switch in heavy traffic. We chose this system because it is one of the simplest queueing systems where the CRP condition is not satisfied. However, the same approach can be applied to any queueing system where the CRP condition is not met, which is what we discuss in this subsection. Specifically, we focus on a generalized switch with $n$ queues, where SSC occurs into a $d$-dimensional subspace.

\citeauthor{atilla} \cite{atilla} showed how to compute the moments of $\|\vq_{\parallel \cH}\|$ using the drift method in queueing systems that satisfy the CRP condition. In this case, setting to zero the drift of $V(\vq)=\|\vq_{\parallel \cH}\|^{m+1}$ in steady state and using SSC allows to compute the $m\tth$ moment because of the following reason. When one sets to zero the drift of $V(\vq)$, terms of the form $q_{\parallel\cH i}^+ u_{\parallel\cH i}$ arise and, since $\vq_{\parallel \cH}^+$ and $\vu_{\parallel\cH}$ belong to the same one-dimensional subspace, these terms can be approximated by $q_i^+ u_i$, which is zero by definition of unused service.

On the other hand, if the CRP condition is not satisfied, then $\vq$ lives in a $d$-dimensional subspace, where $d>1$. In this case, for each $i$, $q_{\parallel \cH i}^+ u_{\parallel\cH i}$ cannot be approximated by $q_i^+ u_i$ because of the following reason. In heavy traffic we only have the approximation (with some abuse of notation) $q_{\parallel \cH i}^+ u_{\parallel \cH i}\approx q_i^+(u_{k_1}+u_{k_2}+\cdots+u_{k_d})$, where $k_1,\ldots,k_d$ represent the $d$ dimensions that characterize SSC. In other words, cross terms arise exactly as the `$qu$' terms in Theorems \ref{s.theorem:2switch}, \ref{s.theorem:2switch.correlated} and \ref{s.theorem.nxn} for the input-queued switch. In the following analysis we present the number of equations and variables that appear in a general queueing system with $d$-dimensional SSC.

In order to obtain the $m\tth$ moment of the queue lengths, we should construct a system of equations that yields from setting to zero the drift of all the monomials of degree $m+1$. Since SSC occurs into a $d$-dimensional subspace, we need to consider all the possible monomials of degree $m+1$ in $d$ variables. Setting to zero the drift of each monomial will lead to an equation, so we will have $\binom{m+d}{d-1}$ equations. Now we count the number of `new' variables with respect to the system of equations that arises after setting to zero the drift of monomials of degree $k$, for all $k\leq m$. We say a variable is `new' for the system of equations that arises after setting to zero the monomials of degree $m+1$ if it does not appear in any system of equations of degree $k<m+1$. Observe that there are two types of `new' variables that do not vanish in the heavy-traffic limit. On one hand, we have the heavy-traffic limit of the expected value of products of the elements of $\vq_{\parallel \cH}$ and, on the other hand, we have the heavy-traffic limit of the expected value of the product between the elements of $\vq_{\parallel \cH}$ and of the vector of unused service. We will call them the `$q$' variables and the `$qu$' variables, respectively. Specifically, the `$q$' variables are all monomials of degree $m$ in $d$ variables, so there are $\binom{m+d-1}{d-1}$ `$q$' variables. The `$qu$' variables that do not vanish in heavy traffic are of degree $m$ in `$q$' and degree 1 in `$u$'. Also, the element corresponding to the unused service vector has to be different to the elements of the vector of queue lengths because the product between the queue length and the unused service of the same queue is zero by definition of unused service. Therefore, for each element of $\vu_{\parallel\cH}$ we need to consider all possible combinations of `$q$'s, i.e. all monomials of degree $m$ in $d-1$ variables. Therefore, there are $d\binom{m+d-2}{d-2}$ `$qu$' variables. Thus, in total we have $\binom{m+d-1}{d-1}+d\binom{m+d-2}{d-2}$ variables and this number is larger than the number of equations.

Summarizing, if we use the method introduced in this section to compute the $m\tth$ moment of the queue lengths of a queueing system that experiences $d$-dimensional SSC, we obtain a system of equations of $\binom{m+d}{d-1}$ equations and $\binom{m+d-1}{d-1}+d\binom{m+d-2}{d-2}$ variables. Therefore, it is under-determined. In other words, we need extra equations to find a unique solution to this system of equations. This analysis shows that the issues illustrated in Theorem \ref{s.theorem:2switch} arise in any queueing system with multidimensional SSC. 

\section{Proof of Theorems \ref{gs.thm:bounds} and \ref{s.theorem:2switch}.}\label{sec:proofs}
In this section we present the proof of the main theorems of this paper.

\subsection{Proof of Theorem \ref{gs.thm:bounds}.}\label{subsec:proof.gen.switch}

In this section we present the proof of the main theorem. We use the notation $\Em{\,\cdot\,}=\E{\,\cdot\,\left|\Mbar=m\right. }$, and we omit the dependence on $\epsilon$ of the variables for simplicity of exposition. Before presenting the proof of the theorem, we present two lemmas that formalize some intuition about the random variables $\Bbar_\ell$ and are essential in the proof of \Cref{gs.thm:bounds}.

Recall that $\Bbar_\ell\defn b^{(\Mbar,\ell)}$ and that, for each $m\in\cM$, $\bml$ is the maximum $\vcl$-weighted service rate in $\cS^\pm$. Similarly, $\bl$ can be interpreted as the maximum $\vcl$-weighted service rate in $\cC$, and $\vcl$ and $\bl$ define a facet of $\cC$. Hence, since the values of $\Bbar_\ell$ occur according to the probability mass function of the channel state, and the capacity region $\cC$ can be interpreted as the `expected capacity region' according to \eqref{gs.eq.cap.reg.def.weighted.CH}, we should expect $\E{\Bbar_\ell}=\bl$. Additionally, $\vcl$ and $\bml$ define a half-space that passes through the boundary of $ConvexHull(\cS^\pm)$ and, hence, there must exist a vector $\vnu^\pm$ such that $\bml=\langle \vcl,\vnu^\pm\rangle$. We formalize these results in \Cref{gs.lemma:bjl.bl}. 

\begin{lemma}\label{gs.lemma:bjl.bl}
	Let $\ell\in P$ and $m\in\cM$. Then, there exists $\vnu^\pm\in \cS^\pm$ such that $\bml=\langle \vcl,\vnu^\pm\rangle$. This implies that $\bl=\E{\Bbar_\ell}$ for all $\ell\in P$.
\end{lemma}

The proof of \Cref{gs.lemma:bjl.bl} follows immediately from the definition of the capacity region $\cC$ in \eqref{gs.eq.cap.reg.def.weighted.CH}, and of the parameters $\bml$ in \eqref{gs.eq.cl.weighted.def}. We present the details un Appendix \ref{app:lemma.gs.bjl.bl}.
 
As $\epsilon$ gets closer to zero, we know that $\vlambda^\peps$ gets closer to $\vnu$, and SSC implies that the vector of queue lengths can be approximated by its projection on $\cK$. In other words, as $\epsilon\dto 0$ the vector of queue lengths can be well approximated by a conic combination of the vectors $\vcl$ with $\ell\in P$. Therefore, since the scheduling problem is solved using MaxWeight algorithm, and given that the channel state is $m$, one should expect that $\langle \vcl,\vsbar\rangle = \bml$  with high probability. In the next lemma we formalize this intuition.

\begin{lemma}\label{gs.lemma:cl.s.bl}
	For each $m\in \cM$ and $\ell\in P$, define $\piml\defn \,\Prob{\left.\langle \vcl,\vsbar\rangle =\bml\,\right|\, \Mbar=m }$. Then, \newline $1-\piml$ is $O(\epsilon)$.
\end{lemma}

The proof of \Cref{gs.lemma:cl.s.bl} is a generalization of \cite[Claim 1]{atilla}, and we present it in Appendix \ref{app:gs.lemma.cl.s.bl} for completeness.

Now we prove \Cref{gs.thm:bounds}. 

\proof{Proof of \Cref{gs.thm:bounds}.}
	
	First observe that $\langle \vqbar,\vw\rangle = \langle\vqbarparh,\vnu\rangle$. To show this statement, define $\vw_\perp\defn \vw-\vnu$ for all $\vw\in\cap_{\ell\in P}\cF^\pl$, and observe that $\langle\vcl,\vw_\perp\rangle=0$ because both $\vnu,\vw\in \cF^\pl$ for all $\ell\in P$. Then, 
	\begin{align*}
		\langle\vqbarparh,\vnu\rangle = \langle \vqbarparh,\vw-\vw_\perp\rangle = \langle\vqbarparh,\vw \rangle \stackrel{(*)}{=} \langle\vqbar,\vw\rangle,
	\end{align*}
	where $(*)$ holds because $\vw\in\cap_{\ell\in P}\cF^\pl$ and because $\vqbarparh=\vqbar-\vqbarperph$. Hence, in the rest of the proof we focus on computing bounds for $\E{\langle \vqbarparh,\vnu\rangle}$.
	
	We set to zero the drift of $V_{\parallel \cH}(\vq)=\left\| \vq_{\parallel \cH}\right\|^2$, and bound separately each of the terms that arise. Before setting the drift to zero we need to make sure that $\E{V_{\parallel \cH}(\vqbarparh) }$ is finite. This result can be proved using the Foster-Lyapunov theorem with Lyapunov function $Z(\vq)=\left\|\vq \right\|^2$.  This proves that $\E{\left\|\vqbar \right\|^2}$ is finite. Then, since projection is nonexpansive we have $\E{\left\|\vqbarparh\right\|^2}$ is also finite. The proof is simple, so we omit the details for ease of exposition. Now, setting to zero the drift of $V_{\parallel \cH}(\vq)$ we obtain
	\begin{align}
		0 &= \E{\left\|\vqbarparh^+ \right\|^2 - \left\|\vqbarparh \right\|^2} \nonumber \\
		&\stackrel{(a)}{=} \E{\left\|\vabarparh - \vsbarparh\right\|^2 +2\langle \vqbarparh,\, \vabarparh-\vsbarparh \rangle -\left\|\vubarparh \right\|^2 + 2\langle\vqbarparh^+,\, \vubarparh\rangle} \label{gs.eq.drift.zero}
	\end{align}
	where $(a)$ holds by the dynamics of the queues presented in \eqref{gs.eq.dynamics.queues}, and reorganizing terms. Let 
	\begin{align*}
		& \cT_1\defn 2\E{\langle\vqbarparh,\vsbarparh-\vabarparh \rangle}, \quad \cT_2\defn \E{\left\|\vabarparh-\vsbarparh \right\|^2},\\
		& \cT_3\defn \E{\left\|\vubarparh \right\|^2}\quad \text{and} \quad \cT_4\defn 2\E{\langle \vqbarparh^+,\, \vubarparh\rangle}.
	\end{align*}
	Then, reorganizing the terms in \eqref{gs.eq.drift.zero} we obtain $\cT_1=\cT_2-\cT_3+\cT_4$. We compute each term separately. We start with $\cT_1$.
	\begin{align}
		\cT_1&
		\stackrel{(a)}{=} 2\E{\langle \vqbarparh,\vsbar-\vabar\rangle} \nonumber \\
		&\stackrel{(b)}{=} 2\epsilon\E{\langle\vqbarparh,\vnu\rangle} + \E{\langle\vqbarparh, \vsbar-\vnu\rangle} \nonumber \\
		& \stackrel{(c)}{=} 2\epsilon\E{\langle\vqbarparh,\vnu\rangle} + O(\sqrt{\epsilon}),
		\label{gs.eq.T1}
	\end{align}
	where $(a)$ holds by the orthogonality principle; $(b)$ holds because $\E{\vabar}=(1-\epsilon)\vnu$ and because $\vabar$ is independent of the vector of queue lengths; and $(c)$ holds by \Cref{gs.claim.T1.partial} stated below.
	\begin{claim}\label{gs.claim.T1.partial}
		Consider a set of generalized switches as described in \Cref{gs.thm:bounds}. Then, 
		\begin{align*}
			\left|\E{\langle \vqbarparh, \vsbar-\vnu\rangle} \right| \text{ is }O(\sqrt{\epsilon}).
		\end{align*}
	\end{claim}
	We present the proof of \Cref{gs.claim.T1.partial} in Appendix \ref{app:proof.claim.T1}. Now we compute $\cT_2$. Expanding the product we obtain
	\begin{align}
		\cT_2&= \E{\left\|\vabarparh-\vsbarparh \right\|^2} = \E{\left\|\vabarparh\right\|^2}+ \E{\left\|\vsbarparh \right\|^2}-2\E{\langle\vabarparh,\,\vsbarparh\rangle}. \label{gs.eq.T2.partial}
	\end{align}
	
	We compute each term in \eqref{gs.eq.T2.partial} separately. For the first two terms, we solve the least squares problem and we use the projection matrix on the subspace $\cH$, denoted as $H$. Let $h_{i,j}$ be its element $(i,j)$ for each $i,j\in[n]$. For the first term we have
	
	\begin{align}
		\E{\left\|\vabarparh\right\|^2} &= \E{\left\|H\,\vabar \right\|^2} \nonumber \\
		&\stackrel{(a)}{=} \sum_{i=1}^n\sum_{j=1}^n h_{i,j}\Cov{\abar_i,\abar_j} + \sum_{i=1}^n \sum_{j=1}^n h_{i,j}\E{\abar_i}\E{\abar_j} \nonumber \\
		&\stackrel{(b)}{=} \vone^T\left(H\circ \Sigma_a^\peps\right)\vone + (1-\epsilon)^2 \vnu^T H\vnu, \label{gs.eq.T2.a}
	\end{align}
	where $(a)$ holds solving the least squares problem, by definition of norm, because $H$ is a projection matrix (and therefore $H=H^T=H^2$), and by definition of covariance; and $(b)$ holds by definition of the Hadamard's product and because $\E{\abar_i}=\lambda_i^\peps=(1-\epsilon)\nu_i$ for each $i\in[n]$. For the second term in \eqref{gs.eq.T2.partial} we obtain
	\begin{align}
		\E{\left\| \vsbarparh\right\|^2} &= \E{\| H\vsbar\|^2} \nonumber \\
		&\stackrel{(a)}{=} \E{\vsbar^T C(C^TC)^{-1}C^T\vsbar} \nonumber \\
		&\stackrel{(b)}{=} \sum_{\ell_1\in\tilde{P}} \sum_{\ell_2\in \tilde{P}} (C^TC)^{-1}_{\ell_1,\ell_2}\E{\langle \vclone,\vsbar\rangle \langle\vcltwo,\vsbar\rangle} \nonumber \\
		&\stackrel{(c)}{=} \sum_{\ell_1\in\tilde{P}} \sum_{\ell_2\in \tilde{P}} (C^TC)^{-1}_{\ell_1,\ell_2}\sum_{m\in\cM}\psi_m\Em{\langle \vclone,\vsbar\rangle \langle\vcltwo,\vsbar\rangle}	\nonumber \\
		&\stackrel{(d)}{=} \vone^T \left((C^TC)^{-1}\circ \Sigma_B \right) \vone + \vnu^T H\vnu - O(\epsilon). \label{gs.eq.T2.s}
	\end{align}
	where $(C^TC)^{-1}_{\ell_1,\ell_2}$ is the element $(\ell_1,\ell_2)$ of the matrix $(C^TC)^{-1}$ for each $\ell_1,\ell_2\in \tilde{P}$. Here, $(a)$ holds solving the least squares problems, and because $H=C(C^TC)^{-1}C^T$ by definition of projection matrix; $(b)$ holds by definition of matrix multiplication, and because $C^T\vsbar$ is a vector with elements $\langle\vcl,\vsbar\rangle$ for $\ell\in \tilde{P}$; $(c)$ holds by law of total probability, conditioning on the channel state; and $(d)$ holds using \Cref{gs.lemma:cl.s.bl}, the definition of covariance and reorganizing the terms. 
	Now we compute the last term in \eqref{gs.eq.T2.partial}. We obtain
	\begin{align}
		-2\E{\langle\vabarparh,\vsbarparh\rangle}
		&\stackrel{(a)}{=}
		-2\E{\vabar^TH \vsbar} 
		\stackrel{(b)}{=} -2(1-\epsilon)\vnu^T\E{H\vsbar} 
		\stackrel{(c)}{=} -2(1-\epsilon)\vnu^T H\vnu+ O(\epsilon), \label{gs.eq.T2.as}
	\end{align}
	where $(a)$ holds because for any vector $\vx$, we have $\vx_{\parallel \cH}=H\vx$ by the solution of the least squares problem, and because $H$ is a projection matrix; $(b)$ holds because $\vabar$ is independent of $\vsbar$ and $\E{\vabar}=\vlambda^{\peps}=(1-\epsilon)\vnu$; and $(c)$ because $H=C\left(C^TC\right)^{-1}C^T$, because $C^T\vsbar$ has elements $\langle\vcl,\vsbar\rangle$ with $\ell\in \tilde{P}$, by \Cref{gs.lemma:bjl.bl} and \Cref{gs.lemma:cl.s.bl}, and  because $\vnu\in\cF^\pl$. 
	Therefore, using \eqref{gs.eq.T2.a}, \eqref{gs.eq.T2.s} and \eqref{gs.eq.T2.as} in \eqref{gs.eq.T2.partial} we obtain
	\begin{align}\label{gs.eq.T2}
		\left|\cT_2 - \left(\vone^T\left(H\circ \Sigma_a^\peps\right)\vone + \vone^T \left((C^TC)^{-1}\circ \Sigma_B \right)\vone + \epsilon^2\vnu^T H\vnu\right) \right| \text{ is } O(\epsilon).
	\end{align}

	Now we compute $\cT_3$. We obtain
	\begin{align*}
		0\leq \cT_3 &= \E{\left\|\vubarparh \right\|^2} 
		\stackrel{(a)}{\leq} \sum_{\ell\in P}\E{\langle\vcl,\vubar\rangle^2} 
		\stackrel{(b)}{\leq} n\smax C_{\max} \sum_{\ell\in P}\E{\langle\vcl,\vubar\rangle} \stackrel{(c)}{=}O(\epsilon),
	\end{align*}
	where $\ds C_{\max}=\max_{\ell\in P,i\in[n]}\left\{c^\pl_i\right\}$ and it is a finite constant. Here, $(a)$ holds because the vectors $\vcl$ are not necessarily orthogonal for all $\ell\in P$; $(b)$ holds because, by definition of the unused service, we know $\vubar\leq \vsbar\leq \smax\vone$ with probability 1; and $(c)$ holds by \Cref{gs.claim.Eu}. 
	\begin{claim}\label{gs.claim.Eu}
		Consider a set of generlized switches, as described in \Cref{gs.thm:bounds}. Then,
		\begin{align}\label{gs.eq.sum.cu}
			\sum_{\ell\in P} \E{\langle\vcl,\vubar\rangle} = \epsilon \sum_{\ell\in P}\bl - O(\epsilon).
		\end{align}
	\end{claim}
	We present the proof of \Cref{gs.claim.Eu} in Appendix \ref{app.proof.gs.claim.Eu}. Therefore, 
\begin{align}\label{gs.eq.T3}
\cT_3=O(\epsilon).
\end{align}

The last step is to prove that
\begin{align}\label{gs.eq.T4}
	|\cT_4|=O(\sqrt{\epsilon}).
\end{align}	
We provide the proof in Appendix \ref{app:gs.T4}. Putting equations \eqref{gs.eq.T1}, \eqref{gs.eq.T2}, \eqref{gs.eq.T3} and \eqref{gs.eq.T4} together we obtain
\begin{align*}
\begin{aligned}
\left|\E{\langle\vqbar^\peps, \vw\rangle}- \dfrac{1}{2\epsilon} \left(\vone^T\left(H\circ \Sigma_a^\peps\right)\vone + \vone^T\left((C^TC)^{-1}\circ \Sigma_B\right)\vone \right)\right| \leq K(\epsilon).
\end{aligned}
\end{align*}

This completes the proof.
\Halmos \endproof

Clearly, the above result and proof are much more general and more involved than the proof in the special case of an input-queued switch developed by \citeauthor{MagSri_SSY16_Switch} \cite{MagSri_SSY16_Switch} and \citeauthor{QUESTA_switch} \cite{QUESTA_switch}.
The bound in \Cref{gs.thm:bounds} is expressed in terms of a general projection of the second moments of arrival and service processes onto the space $\cH$. We would like to point out a couple of conceptual differences from the proof in the case of input-queued switch.
Firstly, in the proof of asymptotic upper bounds in an input-queued switch, the scheduling policy is not used. This means that for an input-queued switch, any scheduling policy that exhibits SSC also has the same asymptotic upper bounds.
In our proof here, we use the scheduling policy to upper bound the term  $\cT_1$ in \Cref{gs.claim.T1.partial}. Thus, we may not claim that any scheduling policy that exhibits SSC in \Cref{gs.prop:gen.switch.SSC} satisfies the bound in \Cref{gs.thm:bounds}.
Secondly, while SSC into the cone $\cK$ was established by \citeauthor{MagSri_SSY16_Switch} \cite{MagSri_SSY16_Switch} and \citeauthor{QUESTA_switch} \cite{QUESTA_switch} in the case of an input-queued switch, only the weaker result about collapse into the space $\cH$ was used to obtain heavy-traffic queue length bounds. In contrast, we use the collapse into the cone  $\cK$ in the proof of  \Cref{gs.thm:bounds} to lower bound the term  $\cT_1$. Both these differences are due to the fact that $\vsbarparh$ is constant for all maximal schedules $\vsbar \in \cS$ in the case of an input-queued switch, whereas in the case of the generalized switch this is not necessarily true.

\subsection{Proof of Theorem \ref{s.theorem:2switch}.}\label{subsec:proof.switch}

For ease of exposition, in this proof we use subscript $\parallel$ instead of $\parallel \cH$, since we only use projection on the subspace $\cH$ and not on the cone $\cK$.

\proof{Proof of Theorem \ref{s.theorem:2switch}.}
We know that SSC occurs into a subspace of dimension $2N-1=3$. Therefore, 3 variables are necessary to compute the most general quadratic polynomial. In fact, we know $\qbar_{\parallel 4}=\qbar_{\parallel 2}+\qbar_{\parallel 3}-\qbar_{\parallel 1}$. Then, we only need to consider the variables $\qbar_{\parallel 1}$, $\qbar_{\parallel 2}$ and $\qbar_{\parallel 3}$. The most general quadratic polynomial with these variables is
\begin{align*}
V(\vq)=& \alpha_1 \qbarpari{1}^2 + \alpha_2 \qbarpari{2}^2 + \alpha_3 \qbarpari{3}^2 + \alpha_4 \qbarpari{1}\qbarpari{2} + \alpha_5 \qbarpari{1}\qbarpari{3} + \alpha_6 \qbarpari{2}\qbarpari{3},
\end{align*}
where $\alpha_i\in\bR$ for all $i\in[6]$.

Setting to zero the drift of $V(\vq)$ is equivalent to setting to zero the drift of each monomial separately. Then, we set to zero the drift of the following 6 test functions:
\begin{align*}
& V_1(\vq)=\qbarpari{1}^2,\; V_2(\vq)=\qbarpari{2}^2,\; V_3(\vq)=\qbarpari{3}^2,\\
& V_4(\vq)=\qbarpari{1}\qbarpari{2},\; V_5(\vq)=\qbarpari{1}\qbarpari{3}\;\text{and}\; V_6(\vq)= \qbarpari{2}\qbarpari{3}.
\end{align*}

Before setting to zero the drift of $V_i(\vq)$ for $i\in[6]$ observe that, by definition of the cone $\cK$ in \eqref{eq.switch.cone} we have for any vector $\vy\in \bR^{4}$
\begin{align}
& \begin{aligned}\label{s.2.eq.par11}
y_{\parallel 1}=& \dfrac{y_{1}+y_{2}}{2}+ \dfrac{y_{1}+y_{3}}{2}-\dfrac{y_{1}+y_{2}+y_{3}+y_{4}}{4}
= \dfrac{3y_{1}+y_{2}+y_{3}-y_{4}}{4},
\end{aligned} \\[3pt]
& \begin{aligned}\label{s.2.eq.par12}
y_{\parallel 2}=& \dfrac{y_{1}+y_{2}}{2}+ \dfrac{y_{2}+y_{4}}{2}-\dfrac{y_{1}+y_{2}+y_{3}+y_{4}}{4} = \dfrac{y_{1}+3y_{2}-y_{3}+y_{4}}{4},
\end{aligned}  \\[3pt]
& \begin{aligned}\label{s.2.eq.par21}
y_{\parallel 3}=& \dfrac{y_{3}+y_{4}}{2}+ \dfrac{y_{1}+y_{3}}{2}-\dfrac{y_{1}+y_{2}+y_{3}+y_{4}}{4}
= \dfrac{y_{1}-y_{2}+3y_{3}+y_{4}}{4} .
\end{aligned}
\end{align}
Then, since the switch is completely saturated, we have
\begin{align}\label{s.2.eq.Ea}
\E{\abarpari{i}}=& \dfrac{1-\epsilon}{2}+\dfrac{1-\epsilon}{2}-\dfrac{2(1-\epsilon)}{4}=\dfrac{1-\epsilon}{2} \qquad\forall i\in[4]
\end{align}
and since $\vsbar$ is a maximal schedule we have
\begin{align}\label{s.2.eq.spar}
\sbarpari{i}=& \dfrac{1}{2}+\dfrac{1}{2}-\dfrac{2}{4}=\dfrac{1}{2}\qquad\forall  i\in[4].
\end{align}

We first set to zero the drift of $V_1(\vq)$. We obtain
\begin{align}
0&= \E{\left(\qbarpari{1}^+\right)^2-\qbarpari{1}^2} \nonumber \\
&= \E{\left(\qbarpari{1}^+-\ubarpari{1}+\ubarpari{1}\right)^2- \qbarpari{1}^2} \nonumber \\
&= \E{\left(\qbarpari{1}^+-\ubarpari{1}\right)^2 + \ubarpari{1}^2+2\left(\qbarpari{1}^+-\ubarpari{1}\right)\ubarpari{1} -\qbarpari{1}^2} \nonumber \\
&\stackrel{(*)}{=} \E{\left(\qbarpari{1}+ \abarpari{1}- \sbarpari{1}\right)^2 - \ubarpari{1}^2+ 2\qbarpari{1}^+ \ubarpari{1} -\qbarpari{1}^2} \nonumber \\
&\begin{aligned}\label{s.2.eq.q11.partial}
&= \E{\left(\abarpari{1}-\sbarpari{1} \right)^2 + 2\qbarpari{1}\left(\abarpari{1}-\sbarpari{1}\right)- \ubarpari{1}^2  + 2\qbarpari{1}^+ \ubarpari{1}},
\end{aligned}
\end{align}
where $(*)$ holds by \eqref{gs.eq.dynamics.queues} and reorganizing the terms. We compute each term separately. For the first term we have
\begin{align}
 \E{\left(\abarpari{1}-\sbarpari{1}\right)^2} 
&\stackrel{(a)}{=} \E{\left(\abarpari{1}-\dfrac{1}{2}\right)^2} \nonumber \\
&\stackrel{(b)}{=} \Var{\abarpari{1}}+ \left(\E{\abarpari{1}}\right)^2+ \dfrac{1}{4} -\E{\abarpari{1}} \nonumber \\
&= \Var{\abarpari{1}}+\left(\E{\abarpari{1}}-\dfrac{1}{2}\right)^2 \nonumber \\
&\stackrel{(c)}{=} \Var{\dfrac{3\abar_{1}+\abar_{2}+\abar_{3}-\abar_{4}}{4}}+ \dfrac{\epsilon^2}{4} \nonumber \\
&\stackrel{(d)}{=} \dfrac{9\left(\sigma_{a_1}^\peps\right)^2+\left(\sigma_{a_2}^\peps\right)^2 + \left(\sigma_{a_3}^\peps\right)^2 + \left(\sigma_{a_4}^\peps\right)^2}{16} + \dfrac{\epsilon^2}{4}, \label{s.2.eq.q11.as}
\end{align}
where $(a)$ holds by \eqref{s.2.eq.spar}; $(b)$ holds by definition of variance and reorganizing terms; $(c)$ holds by definition of $\abarpari{1}$ as in \eqref{s.2.eq.par11}, and by \eqref{s.2.eq.Ea}; and $(d)$ holds because the arrival processes to different queues are independent. For the second term we obtain
\begin{align}
2\E{\qbarpari{1}\left(\abarpari{1}-\sbarpari{1}\right)} \stackrel{(a)}{=}& 2\E{\qbarpari{1}\left(\abarpari{1}-\dfrac{1}{2}\right)}
\stackrel{(b)}{=} 2\E{\qbarpari{1}}\left(\E{\abarpari{1}}-\dfrac{1}{2}\right)  
\stackrel{(c)}{=} -\epsilon\E{\qbarpari{1}},
\label{s.2.eq.q11.lhs}
\end{align}
where $(a)$ holds by \eqref{s.2.eq.spar}; $(b)$ holds because the arrival processes are independent of the queue lengths; and $(c)$ holds by \eqref{s.2.eq.Ea}. For the third term, observe
\begin{align*}
0\leq \E{\ubarpari{1}^2}\leq \E{\left\|\vubar_{\parallel} \right\|^2}.
\end{align*}
From the proof of \Cref{gs.thm:bounds}, we know $\E{\left\|\vubar_{\parallel}\right\|^2}$ is $O(\epsilon)$ (see \eqref{gs.eq.T3}). 
Therefore,
\begin{align}\label{s.2.eq.q11.u}
\E{\ubarpari{1}^2}\text{ is }O(\epsilon).
\end{align}

Now we compute the last term. By definition of $\vqbarpar$ and $\vqbarperp$ we have
\begin{align*}
2\E{\qbarpari{1}^+ \ubarpari{1}}=& 2\E{\qbar_{1}^+ \ubarpari{1}} - 2\E{\qbar_{\perp 1}^+\ubarpari{1}}.
\end{align*}

\begin{claim}\label{switch.claim.qperp.upar}
	Consider the queueing system described in Theorem \ref{s.theorem:2switch}. Then,
	\begin{align*}
	\E{\qbar_{\perp 1}^+\ubarpari{1}}\text{ is }O(\sqrt{\epsilon}).
	\end{align*}
\end{claim}
The proof of Claim \ref{switch.claim.qperp.upar} is presented in Appendix \ref{app:switch.proofs.qperp.upar}. Then,
\begin{align*}
2\E{\qbarpari{1}^+ \ubarpari{1}} 
=& 2\E{\qbar_{1}^+ \ubarpari{1}} + O(\sqrt{\epsilon}) \\
\stackrel{(a)}{=}& \dfrac{1}{2}\E{\qbar_{1}^+\left(3\ubar_{1}+\ubar_{2}+\ubar_{3}-\ubar_{4} \right)} + O(\sqrt{\epsilon})  \\[3pt]
\stackrel{(b)}{=}& \dfrac{1}{2}\E{\qbar_{1}^+\left(\ubar_{2}+\ubar_{3}-\ubar_{4} \right)} + O(\sqrt{\epsilon})
\end{align*}
where $(a)$ holds by \eqref{s.2.eq.par11}; and $(b)$ holds by \eqref{gs.eq.qu}.

\begin{claim}\label{switch.claim.q11u22}
	Consider the queueing system described in Theorem \ref{s.theorem:2switch}. Then,
	\begin{align*}
	\E{\qbar_{1}^+\ubar_{4}}= \E{\qbar_{2}^+\ubar_{4}}+\E{\qbar_{3}^+\ubar_{4}} + O(\sqrt{\epsilon}).
	\end{align*}
\end{claim}
The proof of Claim \ref{switch.claim.q11u22} is presented in Appendix \ref{app:switch.proofs.q11u22}. Therefore, we obtain
\begin{align}\label{s.2.eq.q11.qu}
\begin{aligned}
& 2\E{\qbarpari{1}^+ \ubarpari{1}}
= \dfrac{1}{2}\E{\qbar_{1}^+\left(\ubar_{2}+\ubar_{3}\right)} - \dfrac{1}{2}\E{\qbar_{2}^+ \ubar_{4}} - \dfrac{1}{2}\E{\qbar_{3}^+\ubar_{4}}+O(\sqrt{\epsilon})
\end{aligned}
\end{align}

Using \eqref{s.2.eq.q11.as}, \eqref{s.2.eq.q11.lhs}, \eqref{s.2.eq.q11.u} and \eqref{s.2.eq.q11.qu} in \eqref{s.2.eq.q11.partial}, and reorganizing the terms we obtain
\begin{align*}
& \epsilon\E{\qbarpari{1}}\\
=& \dfrac{9\left(\sigma_{a_1}^\peps\right)^2+ \left(\sigma_{a_2}^\peps\right)^2 + \left(\sigma_{a_3}^\peps\right)^2 + \left(\sigma_{a_4}^\peps\right)^2}{16}+ \dfrac{1}{2}\E{\qbar_{1}^+\left(\ubar_{2}+\ubar_{3}\right)} - \dfrac{1}{2}\E{\qbar_{2}^+ \ubar_{4}}- \dfrac{1}{2}\E{\qbar_{3}^+\ubar_{4}} + \dfrac{\epsilon^2}{4}+O(\sqrt{\epsilon}) .
\end{align*}
Taking the limit as $\epsilon\dto 0$ on both sides we obtain \eqref{s.2.eq.1}. The proof of \eqref{s.2.eq.2} and of \eqref{s.2.eq.3} hold similarly, after setting to zero the drift of $V_2(\vq)$ and $V_3(\vq)$ respectively. We omit the details for brevity.

To obtain \eqref{s.2.eq.4} we set to zero the drift of $V_4(\vq)$. After similar manipulation as above, we obtain
\begin{align}
0 &= \E{\qbarpari{1}^+ \qbarpari{2}^+ - \qbarpari{1}\qbarpari{2}} \nonumber \\
& \begin{aligned}\label{s.2.eq.q11q12.partial}
&= \E{\qbarpari{1} \left(\abarpari{2}-\sbarpari{2}\right)} + \E{\qbarpari{2} \left(\abarpari{1}-\sbarpari{1}\right)} + \E{\left(\abarpari{1}-\sbarpari{1}\right)\left(\abarpari{2}- \sbarpari{2}\right)} \\
&\quad +\E{\qbarpari{1}^+\ubarpari{2}} +\E{\qbarpari{2}^+\ubarpari{1}} - \E{\ubarpari{1}\ubarpari{2}}.
\end{aligned}
\end{align}
We compute term by term. For the first term we have
\begin{align}\label{s.2.eq.q11q12.q11a12}
\E{\qbarpari{1} \left(\abarpari{2}-\sbarpari{2}\right)} 
&=-\dfrac{\epsilon}{2}\E{\qbarpari{1}},
\end{align}
where we used that $\sbarpari{2}=\frac{1}{2}$ and independence of the arrivals and queue lengths processes. Similarly, for the second term we obtain
\begin{align}\label{s.2.eq.q11q12.q12a11}
\E{\qbarpari{2}\left(\abarpari{1}-\sbarpari{1}\right)}=& -\dfrac{\epsilon}{2}\E{\qbarpari{2}}.
\end{align}

For the third term we have
\begin{align}
\E{\left(\abarpari{1}-\sbarpari{1}\right)\left(\abarpari{2}- \sbarpari{2}\right)} 
& \stackrel{(a)}{=} \E{\left(\abarpari{1}-\dfrac{1}{2}\right)\left(\abarpari{2}-\dfrac{1}{2}\right)} \nonumber \\
& \stackrel{(b)}{=} \Cov{\abarpari{1},\abarpari{2}} + \E{\abarpari{1}}\E{\abarpari{2}} - \dfrac{1}{2}\E{\abarpari{1}} - \dfrac{1}{2}\E{\abarpari{2}}+\dfrac{1}{4} \nonumber \\[3pt]
& \stackrel{(c)}{=} \Cov{\dfrac{3\abar_{1}+\abar_{2}+\abar_{3}-\abar_{4}}{4}, \dfrac{\abar_{1}+3\abar_{2}-\abar_{3}+\abar_{4}}{4}} + \dfrac{\epsilon^2}{4} \nonumber \\[3pt]
& \stackrel{(d)}{=} \dfrac{3\left(\sigma^\peps_{a_1}\right)^2 + 3\left(\sigma^\peps_{a_2}\right)^2 - \left(\sigma^\peps_{a_3}\right)^2 - \left(\sigma^\peps_{a_4}\right)^2}{16} + \dfrac{\epsilon^2}{4} \label{s.2.eq.q11q12.a}
\end{align}
where $(a)$ holds by \eqref{s.2.eq.spar}; $(b)$ holds by definition of covariance and reorganizing terms; $(c)$ holds by \eqref{s.2.eq.par11}, \eqref{s.2.eq.par12} and \eqref{s.2.eq.Ea}; and $(d)$ holds because the arrival processes to different queues are independent. For the fourth term we have
\begin{align}
\E{\qbarpari{1}^+ \ubarpari{2}} 
&\stackrel{(a)}{=} \E{\qbar_{1}^+ \ubarpari{2}} - \E{\qbar_{\perp 1}^+ \ubarpari{2}} \nonumber \\
&\stackrel{(b)}{=} \E{\qbar_{1}^+ \ubarpari{2}}+O(\sqrt{\epsilon}) \nonumber \\
&\stackrel{(c)}{=} \dfrac{1}{4}\E{\qbar_{1}^+ \left(\ubar_{1}+ 3\ubar_{2}-\ubar_{3} +\ubar_{4} \right)} + O(\sqrt{\epsilon}) \nonumber \\[3pt]
&\stackrel{(d)}{=} \dfrac{1}{4}\E{\qbar_{1}^+ \left( 3\ubar_{2}-\ubar_{3} +\ubar_{4} \right)} + O(\sqrt{\epsilon}) \nonumber \\[3pt]
& \stackrel{(e)}{=} \dfrac{1}{4}\E{\qbar_{1}^+\left(3\ubar_{2} - \ubar_{3} \right)} + \dfrac{1}{4}\E{\qbar_{2}^+\ubar_{4}} + \dfrac{1}{4}\E{\qbar_{3}^+\ubar_{4}} + O(\sqrt{\epsilon}), \label{s.2.eq.q11q12.q11u12}
\end{align}
where $(a)$ holds by definition of $\vqbarpar$ and $\vqbarperp$; $(b)$ holds similarly to \Cref{switch.claim.qperp.upar}; $(c)$ holds by \eqref{s.2.eq.par12}; $(d)$ holds by \eqref{gs.eq.qu}; and $(e)$ holds by Claim \ref{switch.claim.q11u22}. Similarly, for the fifth term we have
\begin{align}\label{s.2.eq.q11q12.q12u11}
\E{\qbarpari{2}^+\ubarpari{1}}=& \dfrac{1}{4} \E{\qbar_{2}^+\left(3\ubar_{1}+\ubar_{3}-\ubar_{4}\right)} +O(\sqrt{\epsilon}).
\end{align}

For the sixth term we have
\begin{align*}
0\leq \E{\ubarpari{1}\ubarpari{2}}\stackrel{(*)}{\leq}& \sqrt{\E{\ubarpari{1}^2}\E{\ubarpari{2}^2} } \\
\leq& \sqrt{\E{\left\|\vubar_{\parallel} \right\|^2} \E{\left\|\vubar_{\parallel} \right\|^2}} \\
=& \E{\left\|\vubar_{\parallel} \right\|^2}
\end{align*}
where $(*)$ holds by the Cauchy-Schwarz inequality. Also, since $\E{\left\|\vubar_{\parallel} \right\|^2}$ is $O(\epsilon)$, we obtain
\begin{align}\label{s.2.eq.q11q12.u}
\E{\ubarpari{1}\ubarpari{2}}\text{ is }O(\epsilon).
\end{align}

Using \eqref{s.2.eq.q11q12.q11a12}, \eqref{s.2.eq.q11q12.q12a11}, \eqref{s.2.eq.q11q12.a}, \eqref{s.2.eq.q11q12.q11u12}, \eqref{s.2.eq.q11q12.q12u11} and \eqref{s.2.eq.q11q12.u} in \eqref{s.2.eq.q11q12.partial}, and reorganizing terms we obtain
\begin{align*}
& \epsilon\E{\qbarpari{1}}+\epsilon\E{\qbarpari{2}} \\[1pt]
=& \dfrac{3\left(\sigma^\peps_{a_1}\right)^2 + 3\left(\sigma^\peps_{a_2}\right)^2 - \left(\sigma^\peps_{a_3}\right)^2 - \left(\sigma^\peps_{a_4}\right)^2}{8} + \dfrac{\epsilon^2}{2} +O(\sqrt{\epsilon}) \\
&+ \dfrac{1}{2}\E{\qbar_{1}^+\left(3\ubar_{2}- \ubar_{3} \right)} + \dfrac{1}{2}\E{\qbar_{3}^+\ubar_{4}} + \dfrac{1}{2} \E{\qbar_{2}^+\left(3\ubar_{1}+\ubar_{3}\right)}.
\end{align*}
Taking the limit as $\epsilon\dto 0$ on both sides we obtain \eqref{s.2.eq.4}. The proof of \eqref{s.2.eq.5} and \eqref{s.2.eq.6} hold similarly, after setting to zero the drift of $V_5(\vq)$ and $V_6(\vq)$, respectively. This completes the proof of \Cref{s.theorem:2switch}.	
\Halmos
\endproof

\section{Conclusion.}\label{sec:conclusion}
In this paper we studied one of the most general single-hop SPNs with control in service: the generalized switch. This model subsumes several queueing systems, such as the input-queued switch, parallel-server systems, ad hoc wireless networks, etc. Our result is widely applicable, since we do not assume the CRP condition, neither independence of the arrival processes. 

We showcase the generality of our result with three particular SPNs: the input-queued switch, parallel-server systems, and an ad hoc wireless network. Each of these results are interesting by themselves since they have been studied separately in the literature, and we can easily compute them as applications of \Cref{gs.thm:bounds}. 

Additionally, we prove that if the heavy-traffic limit is to a vertex of the capacity region, then SSC does not result in a reduction on the dimension of the state space. In other words, in this case we observe full-dimensional SSC. Under this condition, regardless of the correlation among arrival processes, the mean of the linear combinations of the queue lengths that we obtain behave as if the queues were independent in heavy traffic.

Our result is widely applicable to several SPNs, but it only allows to compute certain linear combinations of the queue lengths. In the case of an input-queued switch, this linear combination turns out to be the total queue length, and in parallel-server systems, the weights of the linear combination are the mean service rates. 

We also show that obtaining other linear combinations is a nontrivial problem, since using the drift method with polynomial test functions is equivalent to solving an under-determined system of linear equations. The results we obtain in this paper can be also obtained by taking specific linear combinations of these equations, such that some unknowns cancel out. An immediate line of future work is to extend the method so that all the linear combinations can be computed. This would allow us to also obtain higher moments and, eventually, the joint distribution of the queue lengths.

%
%
\section*{Acknowledgments.}

We acknowledge professor Mohit Singh for the recommendation of using least squares problem to prove Theorem \ref{gs.thm:bounds}.


\bibliographystyle{informs2014} 
\bibliography{biblio} 

\begin{thebibliography}{30}
\providecommand{\natexlab}[1]{#1}
\providecommand{\url}[1]{\texttt{#1}}
\providecommand{\urlprefix}{URL }

\bibitem[{Bell \protect\BIBand{} Williams(2001)}]{BellWill01-Nsystem}
Bell SL, Williams RJ (2001) Dynamic scheduling of a system with two parallel
  servers in heavy traffic with resource pooling: Asymptotic optimality of a
  threshold policy. \emph{Annals of Applied Probability} 11(3):608--649.

\bibitem[{Benson et~al.(2010)Benson, Akella, \protect\BIBand{}
  Maltz}]{datacenter_traffic}
Benson T, Akella A, Maltz DA (2010) Network traffic characteristics of data
  centers in the wild. \emph{Proceedings of the 10th ACM SIGCOMM conference on
  Internet measurement} 267--280.

\bibitem[{Bertsimas et~al.(1994)Bertsimas, Paschalidis, \protect\BIBand{}
  Tsitsiklis}]{bertsimas_paschalidis_Tstitsiklis_optimization}
Bertsimas D, Paschalidis IC, Tsitsiklis JN (1994) Optimization of multiclass
  queueing networks: Polyhedral and nonlinear characterizations of achievable
  performance. \emph{The Annals of Applied Probability} 43--75.

\bibitem[{Bertsimas \protect\BIBand{} Tsitsiklis(1997)}]{bertsimas_LPbook}
Bertsimas D, Tsitsiklis JN (1997) \emph{Introduction to linear optimization},
  volume~6 (Athena Scientific Belmont, MA).

\bibitem[{Dai \protect\BIBand{} Lin(2008)}]{dai2008max_pressure}
Dai J, Lin W (2008) Asymptotic optimality of maximum pressure policies in
  stochastic processing networks. \emph{The Annals of Applied Probability}
  18(6):2239--2299.

\bibitem[{Dimakis \protect\BIBand{} Walrand(2006)}]{DimWal_06}
Dimakis A, Walrand J (2006) Sufficient conditions for stability of longest
  queue first scheduling. \emph{Adv. Appl. Prob.} 505--521.

\bibitem[{Eryilmaz \protect\BIBand{} Srikant(2012)}]{atilla}
Eryilmaz A, Srikant R (2012) Asymptotically tight steady-state queue length
  bounds implied by drift conditions. \emph{Queueing Systems} 72(3-4):311--359,
  ISSN 0257-0130.

\bibitem[{Garnett \protect\BIBand{} Mandelbaum(2000)}]{GarMan00-Nsystem-notes}
Garnett O, Mandelbaum A (2000) An introduction to skills-based routing and its
  operational complexities. \emph{Teaching notes} .

\bibitem[{Ghamami \protect\BIBand{} Ward(2013)}]{GhamWard13-Nsystem}
Ghamami S, Ward AR (2013) Dynamic scheduling of a two-server parallel server
  system with complete resource pooling and reneging in heavy traffic:
  Asymptotic optimality of a two-threshold policy. \emph{Mathematics of
  Operations Research} 38(4):761--824.

\bibitem[{Gupta \protect\BIBand{} Shroff(2010)}]{gupta2010delay}
Gupta G, Shroff N (2010) Delay analysis for wireless networks with single hop
  traffic and general interference constraints. \emph{IEEE/ACM Transactions on
  Networking (TON)} 18(2):393--405.

\bibitem[{Hajek(1982)}]{hajek_drift}
Hajek B (1982) Hitting-time and occupation-time bounds implied by drift
  analysis with applications. \emph{Advances in Applied Probability} 502--525.

\bibitem[{Harrison(1998)}]{har_state_space}
Harrison J (1998) Heavy traffic analysis of a system with parallel servers:
  Asymptotic optimality of discrete review policies. \emph{Ann. App. Probab.}
  822--848.

\bibitem[{Harrison(2013)}]{harrison_2013_book}
Harrison J (2013) \emph{Brownian Models of Performance and Control} (Cambridge
  University Press),
  \urlprefix\url{http://dx.doi.org/10.1017/CBO9781139087698}.

\bibitem[{Harrison \protect\BIBand{} L{\'o}pez(1999)}]{harlop_state_space}
Harrison J, L{\'o}pez M (1999) Heavy traffic resource pooling in
  parallel-server systems. \emph{Queueing Systems} 339--368.

\bibitem[{Kandula et~al.(2009)Kandula, Sengupta, Greenberg, Patel,
  \protect\BIBand{} Chaiken}]{kandula_datacenter_traffic}
Kandula S, Sengupta S, Greenberg A, Patel P, Chaiken R (2009) The nature of
  data center traffic: {M}easurements \& analysis. \emph{Proceedings of the 9th
  ACM SIGCOMM conference on Internet measurement} 202--208.

\bibitem[{Kang \protect\BIBand{} Williams(2012)}]{kang2012diffusion}
Kang W, Williams R (2012) Diffusion approximation for an input-queued switch
  operating under a maximum weight matching policy. \emph{Stochastic Systems}
  2(2):277--321.

\bibitem[{Kang et~al.(2014)Kang, Wang, Jaramillo, \protect\BIBand{}
  Ying}]{KanWanJarYin_wireless-networks}
Kang X, Wang W, Jaramillo JJ, Ying L (2014) On the performance of
  largest-deficit-first for scheduling real-time traffic in wireless networks.
  \emph{IEEE/ACM Transactions on Networking} 24(1):72--84.

\bibitem[{Kingman(1962)}]{kingman}
Kingman J (1962) Some inequalities for the queue {GI/G/1}. \emph{Biometrika}
  315--324.

\bibitem[{Kumar \protect\BIBand{} Kumar(1994)}]{kumar_kumar_lineqns}
Kumar S, Kumar PR (1994) Performance bounds for queueing networks and
  scheduling policies. \emph{IEEE Transactions on Automatic Control}
  39(8):1600--1611, ISSN 0018-9286,
  \urlprefix\url{http://dx.doi.org/10.1109/9.310033}.

\bibitem[{Maguluri et~al.(2018)Maguluri, Burle, \protect\BIBand{}
  Srikant}]{QUESTA_switch}
Maguluri ST, Burle S, Srikant R (2018) Optimal heavy-traffic queue length
  scaling in an incompletely saturated switch. \emph{Queueing Systems}
  88(3-4):279--309.

\bibitem[{Maguluri \protect\BIBand{} Srikant(2016)}]{MagSri_SSY16_Switch}
Maguluri ST, Srikant R (2016) Heavy traffic queue length behavior in a switch
  under the {M}ax{W}eight algorithm. \emph{Stochastic Systems} 6(1):211--250,
  \urlprefix\url{http://dx.doi.org/10.1214/15-SSY193}.

\bibitem[{Meyn(2009)}]{Mey_08}
Meyn S (2009) Stability and asymptotic optimality of generalized maxweight
  policies. \emph{SIAM Journal on Control and Optimization} 47(6):3259--3294.

\bibitem[{Shah et~al.(2011)Shah, Tsitsiklis, \protect\BIBand{}
  Zhong}]{shah_switch_open}
Shah D, Tsitsiklis J, Zhong Y (2011) Optimal scaling of average queue sizes in
  an input-queued switch: an open problem. \emph{Queueing Systems}
  68(3-4):375--384, ISSN 0257-0130.

\bibitem[{Shi et~al.(2019)Shi, Wei, \protect\BIBand{}
  Zhong}]{zhong2019_process_flexibility}
Shi C, Wei Y, Zhong Y (2019) Process flexibility for multiperiod production
  systems. \emph{Operations Research} .

\bibitem[{Siva Theja~Maguluri \protect\BIBand{} Srikant(2011)}]{sivhajsri_11}
Siva Theja~Maguluri BH, Srikant R (2011) The stability of longest-queue-first
  scheduling with variable packet sizes. \emph{Proc. Conf. on Decision and
  Control}.

\bibitem[{Stolyar(2004)}]{stolyar2004maxweight}
Stolyar A (2004) Max{W}eight scheduling in a generalized switch: State space
  collapse and workload minimization in heavy traffic. \emph{Annals of Applied
  Probability} 1--53.

\bibitem[{Tassiulas \protect\BIBand{} Ephremides(1992)}]{TasEph_92}
Tassiulas L, Ephremides A (1992) Stability properties of constrained queueing
  systems and scheduling policies for maximum throughput in multihop radio
  networks. \emph{IEEE Transactions on Automatic Control} 37(12):1936--1948.

\bibitem[{Wang et~al.(2018)Wang, Maguluri, Srikant, \protect\BIBand{}
  Ying}]{WeinaBandwidthJournal}
Wang W, Maguluri S, Srikant R, Ying L (2018) Heavy-traffic insensitive bounds
  for weighted proportionally fair bandwidth sharing policies. \emph{arXiv
  preprint arXiv:1808.02120} .

\bibitem[{Williams(2000)}]{Williams_CRP}
Williams R (2000) On dynamic scheduling of a parallel server system with
  complete resource pooling. \emph{Fields Institute Communications}
  28(49-71):5--1.

\bibitem[{Williams(2016)}]{williams_survey_SPN}
Williams R (2016) Stochastic processing networks. \emph{Annual Review of
  Statistics and Its Application} 3:323--345.

\end{thebibliography}

\newpage

\renewcommand{\theHsection}{A\arabic{section}}
\begin{APPENDICES}
	

	\section{Proof of Proposition \ref{gs.prop.ulb}.}\label{app:prop.gs.ulb}
	\proof{Proof of Proposition \ref{gs.prop.ulb}.}
	We compute a lower bound for the expected queue length of the single server queue $\left\{\Phi^\peps(k):k\in\bZ_+ \right\}$ in steady state. 
	Let $\chi^\peps(k)$ the unused service in time slot $k$. We assume that in each time slot, arrivals occur before service. Then, for each $k\in\bZ_+$ we have
	\begin{align}\label{gs.eq.evolution.phi}
		\Phi^\peps(k+1)=\Phi^\peps(k)+\alpha^\peps(k)-\beta(k)+\chi^\peps(k).
	\end{align}
	
	Before computing the lower bound we need to make sure that the Markov chain $\left\{\Phi^\peps(k):k\geq 1 \right\}$ is positive recurrent for each $\epsilon\in(0,1)$. To do that we show that $\E{\beta(k)-\alpha^\peps(k)}>0$ for all $\epsilon\in(0,1)$. By definition of $\alpha^\peps(k)$ we have 
	\begin{align}
		\E{\alpha^\peps(k)}&= \E{\langle\vz,\va^\peps(k)\rangle} 
		= (1-\epsilon)\langle\vz,\vnu\rangle \label{gs.eq.Ealpha}
	\end{align}
	where the last equality holds because $\E{\vabar^\peps}=(1-\epsilon)\vnu$. 	
	By definition of $\beta(k)$ we have
	\begin{align}
		\E{\beta(k)} &= \sum_{\ell\in P}r^\pl \sum_{m\in\cM}\psi_m \bml \nonumber \\
		&\stackrel{(a)}{=} \sum_{\ell\in P}r^\pl\bl \nonumber \\
		& \stackrel{(b)}{=} \sum_{\ell\in P}r^\pl\langle\vcl,\vnu\rangle \nonumber \\
		& \stackrel{(c)}{=} \langle\vz,\vnu\rangle, \label{gs.eq.Ebeta}
	\end{align}
	where $(a)$ holds because $\bl= \sum_{m\in\cM}\psi_m \bml$ (as we prove later in \Cref{gs.lemma:bjl.bl}); $(b)$ holds because $\vnu\in\cF^\pl$ for all $\ell \in P$; and $(c)$ holds by definition of $\vz$.  Then, from \eqref{gs.eq.Ealpha} and \eqref{gs.eq.Ebeta} we obtain
	\begin{align*}
		\E{\beta(k)-\alpha^\peps(k)}=\epsilon \langle\vz,\vnu\rangle,
	\end{align*}
	which is a positive number. Let $\Phibar^\peps$ be a steady-state vector which is limit in distribution of $\left\{\Phi^\peps(k):k\geq 1 \right\}$, and $\left(\Phibar^\peps\right)^+ \defn \Phibar^\peps + \alphabar^\peps - \betabar^\peps + \chibar^\peps$, where $\alphabar^\peps$ and $\betabar^\peps$ are steady-state random variables with the distribution of $\alpha^\peps(1)$ and $\beta^\peps(1)$, respectively and $\chibar^\peps$ represents the unused service.
	
	Observe that, by definition of $\bml$, we know that $\beta(k)$ is stochastically greater than $\langle\vz,\vs(k)\rangle$ for all $k\in\bZ_+$. Therefore, $\Phi^\peps(k)$ is stochastically smaller than $\langle\vz,\vq^\peps(k)\rangle$ for all $k$, and hence,
	\begin{align*}
		\E{\langle\vz,\vqbar^\peps\rangle}\geq \E{\Phibar^\peps}.
	\end{align*}
	
	
	We omit the dependence on $\epsilon$ in the rest of this proof, for ease of exposition. It can be easily proved that $\E{\Phibar^2}<\infty$ (e.g., we can use Lemma \ref{lemma:hajek}), and we omit the proof for brevity. Then, we set to zero the drift of $V(\Phi)=\Phi^2$ in steady state, and we obtain
	\begin{align}
		0=& \E{\left(\Phibar^+\right)^2-\Phibar^2} \nonumber \\
		=& \E{\left(\Phibar^+-\chibar \right)^2 + \chibar^2 + 2\left(\Phibar^+-\chibar\right)\chibar-\Phibar^2} \nonumber \\
		\stackrel{(a)}{=}& \E{\left(\Phibar+\alphabar-\betabar\right)^2 -\chibar^2- \Phibar^2} \nonumber \\
		\stackrel{(b)}{=}& \E{\alphabar^2}+\E{\betabar^2}-2\E{\alphabar}\E{\betabar}-2\E{\Phibar}\E{\alphabar-\betabar}-\E{\chibar} \label{gs.eq.ulb.partial}
	\end{align}
	where $(a)$ holds after expanding the product, and because $\Phibar^+\chibar=0$ by definition of unused service; and $(b)$ holds after expanding the product and using independence of the arrival, service and queue length processes.
	
	We compute the terms in \eqref{gs.eq.ulb.partial} one by one. We already established that $\E{\alphabar}=(1-\epsilon)\langle\vz,\vnu\rangle$ and $\E{\betabar}=\langle\vz,\vnu\rangle$. Now we compute the quadratic terms. By definition of $\alphabar$ we have
	\begin{align*}
		\E{\alphabar^2}=& \E{\left(\sum_{i=1}^n z_i \abar_i\right)^2} \\
		=& \sum_{i=1}^n\sum_{j=1}^n z_iz_j\E{\abar_i\abar_j} \\
		\stackrel{(a)}{=}& \sum_{i=1}^n\sum_{j=1}^n z_iz_j\Cov{\abar_i,\abar_j} + \sum_{i=1}^n \sum_{j=1}^n z_iz_j\E{\abar_i}\E{\abar_j} \\
		\stackrel{(b)}{=}& \vz^T \Sigma_a^\peps \vz + (1-\epsilon)^2 \langle\vz,\vnu\rangle^2
	\end{align*}
	where $(a)$ holds by definition of covariance; and $(b)$ holds by definition of covariance matrix and because $\E{\vabar}=\vlambda^\peps=(1-\epsilon)\vnu$. For the service process, by definition of covariance matrix we obtain
	\begin{align*}
		\E{\betabar^2}	=& \vr^T\Sigma_B \vr +\langle\vz,\vnu\rangle^2,
	\end{align*}
	where the last equality holds because $\vnu\in\cF^\pl$ for all $\ell\in P$ and by definition of $\vz$. For the last term we compute an upper bound. By definition of unused service, we have $\chibar\leq \betabar$ with probability 1. Then, 
	\begin{align*}
		\E{\chibar^2}\leq& \E{\betabar\chibar} \\
		{\stackrel{(a)}\leq}& b_{\max}\left(\sum_{\ell\in P}r^\pl \right)\E{\chibar} \\
		\stackrel{(b)}{=}& b_{\max}\left(\sum_{\ell\in P}r^\pl  \right)\E{\betabar-\alphabar} \\
		\stackrel{(c)}{=}& \epsilon \langle\vz,\vnu\rangle b_{\max}\left(\sum_{\ell\in P}r^\pl \right),
	\end{align*}
	where $(a)$ holds by definition of $\betabar$ and $b_{\max}$; $(b)$ holds because $\E{\chibar}=\E{\betabar-\alphabar}$, which can be easily proved by setting to zero the drift of $V_l(\Phi)=\Phi$; and $(c)$ holds by \eqref{gs.eq.Ealpha} and \eqref{gs.eq.Ebeta}.
	
	Putting everything together in \eqref{gs.eq.ulb.partial} and rearranging terms we obtain the result.
	\Halmos
\endproof
	
	\section{Proof of Proposition \ref{gs.prop:gen.switch.SSC}.}\label{app:prop.gs.ssc}
	In this section we present the proof of Proposition \ref{gs.prop:gen.switch.SSC}. The proof is based on Lemma \ref{lemma:hajek}, which we state below for completeness. This lemma is a corollary of the results proved by \citeauthor{hajek_drift} \cite{hajek_drift}, and the version we present here was first stated by \citeauthor{atilla} \cite{atilla}.
\begin{lemma}\label{lemma:hajek}
	For an irreducible and aperiodic Markov Chain  $\left\{X(k):k\in\bZ_+\right\}$ over a countable state space $\cX$, suppose $Z:\cX\to \bR_+$ is a Lyapunov function and consider its drift at $x$, $\Delta Z(x)$. Suppose the following conditions are satisfied.
	
	\begin{enumerate}[label=(C\arabic*)]
		\item There exists $\eta>0$ and $\tau<\infty$ such that
		\begin{align*}
		\E{\Delta Z(x)\,|\, X(k)=x}\leq -\eta\qquad\text{for all $x\in \cX$ with $Z(x)\geq \tau$}.
		\end{align*}
		
		\item There exists $D<\infty$ such that
		\begin{align*}
		\Prob{|\Delta Z(x)|\leq D}=1\qquad\text{for all $x\in \cX$}.
		\end{align*}
	\end{enumerate}
	
	Then, there exist $\theta^\star>0$ and $C^\star<\infty$ such that
	\begin{align*}
	\limsup_{k\to\infty}\E{e^{\theta^* Z(X(k))}}\leq C^*.
	\end{align*}
	
	Further, if the Markov chain $\{X(k): k\in\bZ_+\}$ is positive recurrent, then it converges in distribution to a random variable $\overline{Z}$ for which
	\begin{align*}
	\E{e^{\theta^* \overline{Z}}}\leq C^*,
	\end{align*}
	which directly implies that the moments of $\overline{Z}$ exist and are finite.
\end{lemma}

\proof{Sketch of the proof of Proposition \ref{gs.prop:gen.switch.SSC}.}
	For ease of exposition, we omit the dependence on $\epsilon$ of the random variables in this proof. First observe that $\cK\subset \cH$ by definition. Therefore, for all $t=1,2,\ldots$ we have $\left\|\vqbarperph^\peps \right\|^t\leq\left\|\vqbarperpk^\peps \right\|^t$ with probability 1. This proves the first inequality.
	
	To prove the second inequality, we introduce the following notation. Let
	\begin{align*}
		V(\vq)\defn \|\vq\|^2,\quad V_\parallel(\vq)\defn \|\vqpark\|^2,\quad V_\perp(\vq)\defn \|\vqperpk\|^2\quad \text{and}\quad W_\perp(\vq)\defn \|\vqperpk\|.
	\end{align*}
	We use Lemma \ref{lemma:hajek} with Lyapunov function $W_\perp(\vq)$. We first prove that condition \textit{(C2)} of Lemma \ref{lemma:hajek} is satisfied. 	By definition of drift, we have
	\begin{align*}
		\big|\Delta W_\perp(\vq)\big| =& \big|W_\perp(\vq(k+1))-W_\perp(\vq(k))\big| \ind{\vq(k)=\vq} \\
		=& \left|\Big\|\vqperpk(k+1)\Big\|-\Big\|\vqperpk(k)\Big\| \right|\ind{\vq(k)=\vq} \\
		\stackrel{(a)}{\leq}& \Big\|\vqperpk(k+1)-\vqperpk(k) \Big\| \ind{\vq(k)=\vq} \\
		\stackrel{(b)}{=}& \left\|\vq(k+1)-\vq(k)-\big(\vqpark(k+1)-\vqpark(k) \big) \right\|\ind{\vq(k)=\vq} \\
		\stackrel{(c)}{\leq}& \left(\Big\|\vq(k+1)-\vq(k) \Big\|+\Big\|\vqpark(k+1)-\vqpark(k)\Big\| \right)\ind{\vq(k)=\vq} \\
		\stackrel{(d)}{\leq}& 2\Big\|\vq(k+1)-\vq(k) \Big\|\ind{\vq(k)=\vq} \\
		\stackrel{(e)}{=}& 2\Big\|\vq+\va(k)-\vs(k)+\vu(k)-\vq \Big\|\ind{\vq(k)=\vq} \\
		=& \left(2\sqrt{\sum_{i=1}^n \big|a_i(k)-s_i(k)+u_i(k)
			\big|^2 }\right)\; \ind{\vq(k)=\vq} \\
		\stackrel{(f)}{\leq}& 2\sqrt{n}\max\{\amax,\smax\}\quad\text{with probability 1}, \numberthis \label{gs.eq.ssc.proof.c2}
	\end{align*}
	where $(a)$ holds by triangle inequality; $(b)$ holds by definition of $\vqperpk$; $(c)$ holds by triangle inequality; $(d)$ holds because projection on the cone $\cK$ is non expansive; $(e)$ holds by the dynamics of the queues presented in  \eqref{gs.eq.dynamics.queues}; and $(f)$ holds because $a_i(k)\leq \amax$ with probability 1 and $s_i(k)\leq \smax$ for all $i\in[n]$ and all $k\geq 1$. Therefore, if we let $D=2\sqrt{n}\max\{\amax,\smax\}$ we have that condition \textit{(C2)} is satisfied. 
	
	Now we prove condition \textit{(C1)}. We start with an observation that was first used by \citeauthor{atilla} \cite{atilla}. Note that $W_\perp(\vq)=\sqrt{\|\vqperpk\|^2}$ and $f(x)=\sqrt{x}$ is a concave function. Then, using the definition of concavity and reorganizing terms we have
	\begin{align}\label{gs.eq.concavity}
		\Delta W_\perp(\vq)\leq \dfrac{1}{2\|\vqperpk\|}\left(\Delta V(\vq)-\Delta V_\parallel(\vq) \right).
	\end{align}
	
	We bound the conditional expectation of the terms in the brackets separately. We start with $\Evq{\Delta V(\vq)}$. We obtain
	\begin{align}
		& \Evq{\Delta V(\vq)} \nonumber \\
		=& \Evq{\big\|\vq(k+1)\big\|^2 - \big\|\vq(k)\big\|^2} \nonumber \\
		=& \Evq{\big\|\vq(k+1)-\vu(k)+\vu(k)\big\|^2 - \big\|\vq(k)\big\|^2} \nonumber \\
		=& \Evq{\big\|\vq(k+1)-\vu(k) \big\|^2 + \big\| \vu(k)\|^2 + 2\langle \vq(k+1)-\vu(k),\,\vu(k)\rangle - \big\|\vq(k)\big\|^2} \nonumber \\
		\stackrel{(a)}{=}& \Evq{\big\|\vq(k)+\va(k)-\vs(k) \big\|^2 - \big\|\vu(k)\big\|^2 - \big\|\vq(k)\big\|^2} \nonumber \\
		=& \Evq{\big\|\vq(k)\|^2 + \big\|\va(k)-\vs(k)\big\|^2 + 2\langle \vq(k),\va(k)-\vs(k)\rangle - \big\|\vu(k)\|^2 - \big\|\vq(k)\big\|^2} \nonumber \\
		\stackrel{(b)}{\leq}& \Evq{\big\|\va(k)-\vs(k) \big\|^2} + 2\Evq{\langle \vq(k),\,\va(k)-\vs(k)\rangle \vphantom{\big\| \big\|}}
		\label{gs.eq.cond1.i}
	\end{align}
	where $(a)$ holds by the dynamics of the queues presented in \eqref{gs.eq.dynamics.queues} and because, by definition of inner product and \eqref{gs.eq.qu}; and $(b)$ holds because $\Evq{\big\|\vu(k)\big\|^2}\geq 0$ by definition of norm.
	
	We bound the terms in \eqref{gs.eq.cond1.i} separately. First, observe
	\begin{align}
		\Evq{\big\|\va(k)-\vs(k)\big\|^2}=& \Evq{\sum_{i=1}^n \big(a_i(k)-s_i(k) \big)^2} \nonumber \\
		=& \Evq{\sum_{i=1}^n \left(a_i^2(k)+s_i^2(k)-2a_i(k)s_i(k) \right)} \nonumber \\
		\leq& n(\amax^2+\smax^2), \label{gs.eq.cond1.i.C1}
	\end{align}
	where the inequality holds because $2a_i(k)s_i(k)\geq 0$ with probability 1 for all $i\in[n]$ and all $k\in\bZ_+$; and because $0\leq a_i(k)\leq \amax$ with probability 1 and $0\leq s_i(k)\leq \smax$ for all $i\in[n]$ and all $k\in\bZ_+$. Let $\zeta_1\defn n(\amax^2+\smax^2)$.
	
	On the other hand,
	\begin{align}
		\Evq{\langle \vq(k),\, \va(k)-\vs(k)\rangle}=& \langle \vq,\vlambda^\peps\rangle - \Evq{\langle \vq(k),\vs(k)\rangle} \nonumber \\
		\stackrel{(a)}{=}& \langle \vq,(1-\epsilon)\vnu\rangle - \max_{\vx\in \cC}\langle \vq,\vx\rangle \nonumber \\
		=& -\epsilon \langle \vq,\vnu\rangle + \min_{\vx\in \cC}\langle \vq,\vnu-\vx\rangle \nonumber\\
		\leq& -\epsilon\langle \vq,\vnu\rangle + \langle\vq,\vnu-\vx^* \rangle \label{gs.eq.cond1.i.2.partial},
	\end{align}
	for any $\vx^*\in \cC$. Here, equality $(a)$ holds by definition of $\vlambda^\peps$ and by Lemma \ref{gs.lemma:schedule.in.C}.
	
	We pick $\vx^*=\vnu+\frac{\delta}{2\|\vqperpk\|}\vqperpk$. Before proceeding with the proof, we show that such $\vx^*\in \cC$. To do that, we show that $\langle\vcl,\vx^*\rangle\leq \bl$ for all $\ell\in[L]$. We have two cases. If $\ell\in P$, then
	\begin{align*}
		\langle \vcl,\vx^*\rangle &= \langle \vcl,\vnu\rangle + \dfrac{\delta}{\|\vqperpk\|} \langle \vcl,\vqperpk\rangle 
		\stackrel{(a)}{=} \langle\vcl,\vnu\rangle 
		\stackrel{(b)}{=} \bl
	\end{align*}
	where $(a)$ holds because $\langle \vcl,\vqperpk\rangle=0$ for all $\ell\in P$, by the orthogonality principle; and $(b)$ holds because $\vnu\in \bigcap_{\ell\in P} \cF^\pl$.
	
	If $[L]\setminus P\neq \emptyset$ and $\ell\in [L]\setminus P$ we have $\langle \vcl,\vnu\rangle < \bl$. Then, for each $\ell\notin P$ there exists $\delta^\pl>0$ such that $\vnu+\frac{\delta^\pl}{2\|\vqperpk\|}\vqperpk\in \cC$. Then, since there are finitely many hyperplanes defining $\cC$, we can pick $\ds\delta=\min_{\ell\in[L]\setminus P}\left\{\delta^\ell\right\}$.
	
	Then, from \eqref{gs.eq.cond1.i.2.partial} we obtain
	\begin{align}
		 \Evq{\langle \vq(k),\, \va(k)-\vs(k)\rangle} 
		&\leq -\epsilon\langle\vq,\vnu\rangle + \langle \vq,\vnu-\left(\vnu+\dfrac{\delta}{\big\|\vqperpk\big\|}\vqperpk\right)\rangle \nonumber \\
		&= -\epsilon\langle\vq,\vnu\rangle+\dfrac{\delta}{\big\|\vqperpk \big\|}\langle\vq,\vqperpk\rangle \nonumber \\
		&= -\epsilon\langle\vq,\vnu\rangle + \delta \|\vqperpk\|, \label{gs.eq.cond1.i.2}
	\end{align}
	where the last equality holds because $\vq=\vqpark+\vqperpk$ and $\langle\vqpark,\vqperpk\rangle=0$. Then, using \eqref{gs.eq.cond1.i.C1} and\newline
	\eqref{gs.eq.cond1.i.2} in \eqref{gs.eq.cond1.i} we obtain
	\begin{align}\label{gs.eq.cond1.i.complete}
		\Evq{\Delta V(\vq)}\leq \zeta_1-2\epsilon \langle\vq,\vnu\rangle - 2\delta\big\|\vqperpk\big\|.
	\end{align}
	
	To bound the second term in \eqref{gs.eq.concavity} we use properties of projection. We have
	\begin{align}
		&\Evq{\Delta V_\parallel(\vq)} \nonumber\\
		=& \Evq{\big\| \vqpark(k+1)\big\|^2 - \big\|\vqpark(k)\big\|^2} \nonumber\\
		=& \Evq{\langle \vqpark(k+1)+\vqpark(k),\, \vqpark(k+1)-\vqpark(k) \rangle } \nonumber \\
		=& \Evq{\big\|\vqpark(k+1)-\vqpark(k) \big\|^2} +2\Evq{\langle\vqpark(k),\,\vqpark(k+1)-\vqpark(k)\rangle} \nonumber \\
		\stackrel{(a)}{\geq}& 2\Evq{\langle\vqpark(k),\,\vqpark(k+1)-\vqpark(k)\rangle} \nonumber \\
		\stackrel{(b)}{=}& 2\Evq{\langle \vqpark(k),\, \vq(k+1)-\vq(k)\rangle} - 2\Evq{\langle \vqpark(k),\, \vqperpk(k+1)-\vqperpk(k)\rangle} \nonumber\\
		=& 2\Evq{\langle \vqpark(k),\, \va(k)-\vs(k)+\vu(k)\rangle} - 2\Evq{\langle \vqpark(k),\vqperpk(k+1)\rangle}+2\Evq{\langle \vqpark(k),\,\vqperpk(k)\rangle} \nonumber \\
		\stackrel{(c)}{\geq}& 2\Evq{\langle \vqpark(k),\,\va(k)-\vs(k)\rangle} \nonumber \\
		=& 2\langle \vqpark,(1-\epsilon)\vnu\rangle - 2\Evq{\langle \vqpark(k),\,\vs(k)\rangle} \nonumber \\ 
		=& -2\epsilon\langle \vqpark,\vnu\rangle + 2\Evq{\langle \vqpark,\vnu-\vs(k)\rangle} \label{gs.eq.cond1.ii.partial}
	\end{align}
	where $(a)$ holds because $\Evq{\big\|\vqpark(k+1)-\vqpark(k) \big\|^2}\geq 0$; $(b)$ holds because $\vqpark(k)=\vq(k)-\vqperpk(k)$ for all $k\in\bZ_+$; $(c)$ holds because, since $\vqpark(k)\geq 0$ and $\vu(k)\geq 0$ by definition, then $\langle\vqpark(k),\vu(k)\rangle\geq 0$, $\langle \vqpark(k),\vqperpk(k)\rangle=0$ because they are orthogonal by definition, and $\langle\vqpark(k),\vqperpk(k+1)\rangle\leq 0$ because $\vqpark(k)$ belongs to the cone $\cK$ and $\vqperpk(k+1)$ belongs to the polar cone of $\cK$, defined as $\cK^{\circ} \defn \left\{\vy\in \bR^n:\langle\vx,\vy\rangle\leq 0\;\;\forall \vx\in \cK \right\}$.
	
	Since $\vqpark(k)$ is the projection of $\vq(k)$ on $\cK$, there exist coefficients $\xi_\ell\geq 0$ with $\ell\in P$ such that 
	\begin{align*}
		\vqpark(k)=\sum_{\ell\in P}\xi_\ell \vcl.
	\end{align*}
	Then,
	\begin{align*}
		\Evq{\langle \vqpark,\,\vnu-\vs(k)\rangle}=& \sum_{\ell\in P}\xi_\ell \Evq{\langle \vcl,\vnu\rangle-\langle\vcl,\vs(k)\rangle} \\
		\stackrel{(a)}{=}& \sum_{\ell\in P}\xi_\ell \left(b^\ell-\langle \vcl,\argmax_{\vx\in \cC}\langle\vq,\vx\rangle \rangle \right) \\
		\stackrel{(b)}{=}& \sum_{\ell\in P}\xi_\ell \left(\bl - \langle \vcl,\overline{\vx}\rangle\right) \\
		\stackrel{(c)}{\geq}& \sum_{\ell\in P}\xi_\ell(\bl-\bl)=0
	\end{align*}
	for some $\overline{\vx}\in \cC$. Here, $(a)$ holds because $\vnu\in \bigcap_{\ell\in P} \cF^\pl$ and by Lemma \ref{gs.lemma:schedule.in.C}; $(b)$ holds for some $\overline{\vx}\in \cC$ because $\cC$ is a closed and bounded set, so the maximum is attained at some point in $\cC$; and $(c)$ holds because, since $\overline{\vx}\in \cC$, then $\langle \vcl,\overline{\vx}\rangle\leq \bl$ for all $\ell\in[L]$.
	
	Therefore, from \eqref{gs.eq.cond1.ii.partial} we obtain
	\begin{align}\label{gs.eq.cond1.ii.complete}
		\Evq{\Delta V_\parallel(\vq)}\geq -2\epsilon\langle \vqpark,\vnu\rangle.
	\end{align}
		
	Then, using \eqref{gs.eq.cond1.i.complete} and \eqref{gs.eq.cond1.ii.complete} in \eqref{gs.eq.concavity} we obtain
	\begin{align*}
		\Evq{\Delta W_\perp(\vq)}\leq& \dfrac{1}{2\big\|\vqperpk\big\|}\left(\zeta_1-2\epsilon\langle\vq,\vnu\rangle - 2\delta\big\|\vqperpk\big\|+2\epsilon \langle\vqpark,\,\vnu\rangle \right) \\
		=& \dfrac{\zeta_1}{2\big\|\vqperpk\big\|}-\delta + \dfrac{\epsilon}{\big\|\vqperpk\big\|}\left(\langle \vqpark,\vnu\rangle-\langle\vq,\vnu\rangle \right) \\
		\stackrel{(a)}{=}& \dfrac{\zeta_1}{2\big\|\vqperpk\big\|}-\delta + \dfrac{\epsilon}{\big\|\vqperpk\big\|}\left(-\langle\vqperpk,\vnu\rangle\right) \\
		\stackrel{(b)}{\leq}& \dfrac{\zeta_1}{2\big\|\vqperpk\big\|}-\delta +\epsilon\|\vnu\|
	\end{align*}
	where $(a)$ holds because $\vq=\vqpark+\vqperpk$; and $(b)$ holds by Cauchy-Schwarz inequality. Therefore, if $\epsilon\leq \frac{\delta}{2\|\vnu\|}$ we have
	\begin{align*}
		\Evq{\Delta W_\perp(\vq)}\leq\dfrac{\zeta_1}{2\big\|\vqperpk\big\|}-\dfrac{\delta}{2}.
	\end{align*}
	
	Further, if $\big\|\vqperpk\big\| \geq \frac{2 \zeta_1}{\delta}$ we have $\Evq{\Delta W_\perp(\vq)}\leq -\frac{\delta}{4}$. The last inequality verifies condition \textit{(C1)} with $\eta=\frac{\delta}{4}$ and $\tau=\frac{2 \zeta_1}{\delta}$. This completes the proof.
	\Halmos
\endproof

	
	
	\section{Generalizations of Theorem \ref{s.theorem:2switch}.}\label{app:switch.generalizations}
	
\subsection{System of equations for the $2\times 2$ input-queued switch with correlated arrivals.}\label{app:switch.2x2.correlated}

In this section we provide a generalization of Theorem \ref{s.theorem:2switch} to the case of a switch with correlated arrivals. We omit the proof, since it is similar to the proof of Theorem \ref{s.theorem:2switch}. 

\begin{theorem}\label{s.theorem:2switch.correlated}
	Consider a set of $2\times 2$ input-queued switches operating under MaxWeight, indexed by $\epsilon\in(0,1)$ as described in Corollary \ref{cor.nxn.switch}. Suppose $\Sigma^\peps$ is the covariance matrix of the arrival processes, and $\lim_{\epsilon\dto 0}\Sigma^\peps = \Sigma$ component-wise. Then, the following system of equations is satisfied
	\begin{align*}
		&\begin{aligned}
			&\lim_{\epsilon\dto 0} \epsilon\E{\qbar_1} \\
			&= \dfrac{ 9\Sigma_{1,1} + 6 \Sigma_{1,2} + 6\Sigma_{1,3} - 6\Sigma_{1,4} + \Sigma_{2,2} + 2\Sigma_{2,3} - 2\Sigma_{2,4} + \Sigma_{3,3} - 2\Sigma_{3,4} + \Sigma_{4,4} }{16} \\
			& \quad +\dfrac{1}{2}\lim_{\epsilon\dto 0}\E{\qbar_1^+\left(\ubar_2+\ubar_3 \right)}  -\dfrac{1}{2}\lim_{\epsilon\dto 0}\E{\left(\qbar_2^+ + \qbar_3^+ \right)\ubar_4}
		\end{aligned}\\[3pt]
		&\begin{aligned}
			&\lim_{\epsilon\dto 0} \epsilon \E{\qbar_2} \\
			&= \dfrac{\Sigma_{1,1} + 6\Sigma_{1,2} - 2\Sigma_{1,3} + 2\Sigma_{1,4} + 9\Sigma_{2,2} - 6\Sigma_{2,3} - 6\Sigma_{2,4} + \Sigma_{3,3} - 2 \Sigma_{3,4} + \Sigma_{4,4} }{16}  \\
			&\quad + \dfrac{1}{2} \lim_{\epsilon\dto 0}\E{\qbar_2^+\left(\ubar_1-\ubar_3+\ubar_4\right)}
		\end{aligned}\\[3pt]
		&\begin{aligned}
			&\lim_{\epsilon\dto 0} \epsilon\E{\qbar_3}\\
			&= \dfrac{\Sigma_{1,1} - 2\Sigma_{1,2} + 6\Sigma_{1,3} + 2\Sigma_{1,4} + \Sigma_{2,2} - 6\Sigma_{2,3} - 2\Sigma_{2,4} + 9\Sigma_{3,3} + 6\Sigma_{3,4} + \Sigma_{4,4}}{16} \\
			& \quad + \dfrac{1}{2}\lim_{\epsilon\dto 0} \E{\qbar_3^+\left(\ubar_1 -\ubar_2+\ubar_4\right)}
		\end{aligned}\\[3pt]
		& \begin{aligned}
			&\lim_{\epsilon\dto 0} \epsilon\E{\qbar_1+\qbar_2} \\
			&= \dfrac{3\Sigma_{1,1} + 18\Sigma_{1,2} - 6\Sigma_{1,3} + 6\Sigma_{1,4} + 3\Sigma_{2,2} - 2\Sigma_{2,3} + 2\Sigma_{2,4} - \Sigma_{3,3} + 2\Sigma_{3,4} - \Sigma_{4,4} }{8} \\
			& \quad + \dfrac{1}{2}\lim_{\epsilon\dto 0}\E{\qbar_1^+ \left(3\ubar_2-\ubar_3\right)} + \dfrac{1}{2}\lim_{\epsilon\dto 0}\E{\qbar_2^+ \left(3\ubar_1+\ubar_3\right)} + \dfrac{1}{2}\lim_{\epsilon\dto 0}\E{\qbar_3^+ \ubar_4}
		\end{aligned} \\[3pt]
		& \begin{aligned}
			&\lim_{\epsilon\dto 0} \epsilon\E{\qbar_1+ \qbar_3} \\
			&= \dfrac{3\Sigma_{1,1} - 6\Sigma_{1,2} + 18\Sigma_{1,3} + 6\Sigma_{1,4} - \Sigma_{2,2} + 6\Sigma_{2,3} + 2\Sigma_{2,4} + 3\Sigma_{3,3} + 6\Sigma_{3,4} - \Sigma_{4,4}}{8} \\
			&\quad + \dfrac{1}{2}\lim_{\epsilon\dto 0}\E{\qbar_1^+ \left(-\ubar_2+3\ubar_3\right)} + \dfrac{1}{2}\lim_{\epsilon\dto 0}\E{\qbar_2^+ \ubar_4} + \dfrac{1}{2}\lim_{\epsilon\dto 0}\E{\qbar_3^+ \left(3\ubar_1+ \ubar_2\right)}
		\end{aligned} \\[3pt]
		& \begin{aligned}
			&\lim_{\epsilon\dto 0} \epsilon\E{\qbar_2+ \qbar_3} \\
			&= \dfrac{\Sigma_{1,1} - 2\Sigma_{1,2} + 6\Sigma_{1,3} +2\Sigma_{1,4} -3\Sigma_{2,2} + 18\Sigma_{2,3} + 6\Sigma_{2,4} - 3\Sigma_{3,3} - 2\Sigma_{3,4} + \Sigma_{4,4}}{8} \\
			&\quad + \dfrac{1}{2}\lim_{\epsilon\dto 0} \E{\qbar_2^+ \left(\ubar_1+3\ubar_3+\ubar_4\right)} + \dfrac{1}{2}\lim_{\epsilon\dto 0}\E{\qbar_3^+ \left(\ubar_1+3\ubar_2+\ubar_4\right)} ,
		\end{aligned}
	\end{align*}
	where we omitted the dependence on $\epsilon$ of the variables for ease of exposition.
\end{theorem}

\subsection{System of equations for the $N\times N$ input-queued switch.}\label{app:switch.nxn}

For ease of exposition, in this section we use the matrix-shape interpretation of the switch and we assume the arrivals to different input ports are independent of each other. With a slight abuse of notation, we adhere to the notation used in Section \ref{sec:system.of.equations} for the variables, and we use two subscripts, one for the input port and one for the output port. For example, $q_{i,j}(k)$ is the number of packets in line at input port $i$ and output port $j$, for $i,j\in[N]$. 

Before presenting the theorem we introduce the following notation. For $i,j\in[N]$, define
\begin{align*}
	[N]_i\defn [N]\setminus\{i\} \quad\text{and}\quad [N]_{i,j}\defn [N]\setminus\{i,j\}
\end{align*}

\begin{theorem}\label{s.theorem.nxn}
	Consider a set of input-queued switches operating under MaxWeight, indexed by $\epsilon\in(0,1)$ as described in Corollary \ref{cor.nxn.switch.indep}. Further, for all $i,j\in[N]$ let $\sigma_{i,j}^\peps\defn\Var{\abar_{i,j}^\peps}$ and assume $\sigma_{i,j}^2 = \lim_{\epsilon\dto 0}\left(\sigma_{i,j}^\peps\right)^2$. Then, the following system of equations is satisfied, where we omit the dependence on $\epsilon$ of the variables by ease of exposition. 
	\begin{align}
		& \langle \vqbarpar , \boldsymbol{p}\rangle = \sum_{i=1}^n \qbar_{\parallel i,i}\qquad\forall \boldsymbol{p}\in\cS. \label{s.eq.n.0} \\
		&\begin{aligned}\label{s.eq.n.1}
			& \lim_{\epsilon\dto 0} \epsilon\E{\qbar_{1,j}} \\
			=& \dfrac{1}{2N^3}\left((2N-1)^2\sigma_{1,j}^2+ (N-1)^2\left(\sum_{i'\in[N]_1}\sigma_{i',j}^2+ \sum_{j'\in [N]_j}\sigma_{1,j'}^2\right)+ \sum_{i'\in[N]_1}\sum_{j'\in[N]_j}\sigma_{i',j'}^2 \right) \\
			& + \dfrac{1}{N}\lim_{\epsilon\dto 0} \E{(N-1)\left( \sum_{i\in[N]_1} \qbar_{1,j}^+ \ubar_{i',j} + \sum_{j'\in[N]_j}\qbar_{1,j}^+ \ubar_{1,j'}\right) - \sum_{i'\in[N]_1}\sum_{j'\in[N]_j}\qbar_{1,j}^+ \ubar_{i',j'}} \\
			& \pushright{\forall j\in[N]}
		\end{aligned}\\
		&\begin{aligned}\label{s.eq.n.2}
			& \lim_{\epsilon\dto 0}\epsilon\E{\qbar_{i,1}} \\
			=& \dfrac{1}{2N^3}\left((2N-1)^2\sigma_{i,1}^2 + (N-1)^2\left(\sum_{i'\in[N]_i}\sigma_{i',1}^2 +\sum_{j\in[N]_1} \sigma_{i,j'}^2 \right) +\sum_{i\in[N]_i}\sum_{j'\in[N]_1}\sigma_{i',j'}^2 \right) \\
			& + \dfrac{1}{N}\lim_{\epsilon\dto 0}\E{(N-1) \left(\sum_{i'\in[N]_i} \qbar_{i,1}^+\ubar_{i',1} + \sum_{j'\in[N]_1}\qbar_{i,1}^+\ubar_{i,j'} \right)- \sum_{i'\in[N]_i}\sum_{j'\in[N]_1}\qbar_{i,1}^+ \ubar_{i',j'}} \\
			& \pushright{\forall i\in[N]_1}
		\end{aligned} \\
		& \begin{aligned}\label{s.eq.n.5}
		& \lim_{\epsilon\dto 0} \epsilon\E{\qbar_{1,1}+ \qbar_{i,1}} \\
		=& \dfrac{(N-1)}{N^3}\left((2N-1)\left(\sigma_{1,1}^2+ \sigma_{i,1}^2 \right) + (N-1)\sum_{i\in[N]_{1,i}}\sigma_{i',1}^2 - \sum_{j'\in[N]_1}\sigma_{1,j'}^2 - \sum_{j'\in[N]_1}\sigma_{i,j'}^2 \right) \\
		& + \dfrac{1}{N^3}\sum_{i'\in[N]_{1,i}}\sum_{j\in[N]_1} \sigma_{i',j'}^2 \\
		& + \dfrac{1}{N}\lim_{\epsilon\dto 0}\E{\vphantom{\sum_{j\in[N]_1}} (2N-1)\qbar_{1,1}^+ \ubar_{i,1} + (N-1)\left(\sum_{i'\in[N]_{1,i}}\qbar_{1,1}^+ \ubar_{i',1} + \sum_{j'\in[N]_1}\qbar_{1,1}^+ \ubar_{i,j'}\right) } \\
		& + \dfrac{1}{N}\lim_{\epsilon\dto 0}\E{\vphantom{\sum_{j\in[N]_1}} (2N-1)\qbar_{i,1}^+ \ubar_{1,1} + (N-1)\left(\sum_{i'\in[N]_{1,i}}\qbar_{i,1}^+ \ubar_{i',1} + \sum_{j'\in[N]_1}\qbar_{i,1}^+ \ubar_{1,j'} \right) } \\
		& -\dfrac{1}{N}\lim_{\epsilon\dto 0} \E{\sum_{i'\in[N]_i}\sum_{j'\in[N]_1}\qbar_{1,1}^+ \ubar_{i',j'} + \sum_{i'\in[N]_1}\sum_{j'\in[N]_1}\qbar_{i,1}^+ \ubar_{i',j'} } \\
		& \pushright{\forall i\in[N]_1}
		\end{aligned} \\
		&\begin{aligned}\label{s.eq.n.3}
			& \lim_{\epsilon\dto 0}\epsilon\E{\qbar_{1,j}+\qbar_{1,m}} \\
			=& \dfrac{(N-1)}{N^3}\left((2N-1)\left(\sigma_{1,j}^2+ \sigma_{1,m}^2\right) +(N-1)\sum_{j'\in[N]_{j,m}}\sigma_{1,j'}^2 \right) \\
			&-\dfrac{(N-1)}{N^3}\left( \sum_{i'\in[N]_1}\sigma_{i',j}^2 - \sum_{i'\in[N]_m}\sigma_{i',m}^2 \right) +\dfrac{1}{N^3} \sum_{i'\in[N]_1}\sum_{j'\in[N]_{j,m}}\sigma_{i',j'}^2 \\
			&+\dfrac{1}{N}\lim_{\epsilon\dto 0}\E{ (2N-1)\qbar_{1,j}^+\ubar_{1,m} + (N-1)\left(\sum_{i'\in[N]_1} \qbar_{1,j}^+ \ubar_{i',m} + \sum_{j'\in[N]_{j,m}}\qbar_{1,j}^+ \ubar_{1,j'}\right) } \\
			& + \dfrac{1}{N}\lim_{\epsilon\dto 0} \E{(2N-1)\qbar_{1,m}^+ \ubar_{1,j} + (N-1)\left(\sum_{i'\in[N]_1}\qbar_{1,m}^+\ubar_{i',j} + \sum_{j'\in[N]_{j,m}}\qbar_{1,m}^+ \ubar_{1,j'} \right) } \\
			& -\dfrac{1}{N}\lim_{\epsilon\dto 0}\E{ \sum_{i'\in[N]_1}\sum_{j'\in[N]_m}\qbar_{1,j}^+ \ubar_{i',j'} + \sum_{i'\in[N]_1}\sum_{j'\in[N]_j}\qbar_{1,m}^+\ubar_{i',j'}} \\
			& \pushright{\forall (j,m)\in \mathcal{A}_1}
		\end{aligned}\\
		& \begin{aligned}\label{s.eq.n.4}
			&\lim_{\epsilon\dto 0}\epsilon\E{\qbar_{i,1} + \qbar_{l,1}} \\
			=& \dfrac{(N-1)}{N^3}\left((2N-1)\left(\sigma_{i,1}^2 + \sigma_{l,1}^2 \right) + (N-1)\sum_{i'\in[N]_{i,l}}\sigma_{i',1}^2\right)\\
			& - \dfrac{(N-1)}{N^3}\left(\sum_{j'\in[N]_1}\sigma_{i,j'}^2 - \sum_{j'\in[N]_1}\sigma_{l,j'}^2\right) +\dfrac{1}{N^3} \sum_{i'\in[N]_{i,l}}\sum_{j'\in[N]_1}\sigma_{i',j'}^2 \\
			& + \dfrac{1}{N}\lim_{\epsilon\dto 0} \E{\vphantom{\sum_{j\in[N]_1}} (2N-1)\qbar_{i,1}^+ \ubar_{l,1} + (N-1)\left(\sum_{i'\in[N]_{i,l}} \qbar_{i,1}^+ \ubar_{i',1} + \sum_{j'\in[N]_1}\qbar_{i,1}^+ \ubar_{l,j'} \right) } \\
			&+ \dfrac{1}{N}\lim_{\epsilon\dto 0} \E{\vphantom{\sum_{j\in[N]_1}} (2N-1)\qbar_{l,1}^+ \ubar_{i,1} + (N-1)\left(\sum_{i'\in[N]_{i,l}} \qbar_{l,1}^+ \ubar_{i',1} + \sum_{j'\in[N]_1}\qbar_{l,1}^+ \ubar_{i,j'} \right) } \\
			& -\dfrac{1}{N}\lim_{\epsilon\dto 0}\E{\sum_{i'\in[N]_l}\sum_{j'\in[N]_1}\qbar_{i,1}^+ \ubar_{i',j'} + \sum_{i'\in[N]_i}\sum_{j'\in[N]_1} \qbar_{l,1}^+ \ubar_{i',j'} } \\
			& \pushright{\forall (i,l)\in\mathcal{A}_2}
		\end{aligned} \\
		& \begin{aligned}\label{s.eq.n.6}
			& \lim_{\epsilon\dto 0}\epsilon \E{\qbar_{1,j}+ \qbar_{i,1}} \\
			=& \dfrac{1}{N^3}\left(-(2N-1)\left(\sigma_{1,j}^2 + \sigma_{i,1}^2 \right) + (N-1)^2\sigma_{1,1}^2 + \sum_{i'\in[N]_{1,i}}\sum_{j'\in[N]_{1,j}}\sigma_{i',j'}^2 \right)\\
			&- \dfrac{(N-1)}{N^3}\left(\sum_{i\in[N]_1}\sigma_{i',j}^2 + \sum_{j'\in[N]_{1,j}}\sigma_{1,j'}^2 + \sum_{i'\in[N]_{1,i}}\sigma_{i',1}^2 + \sum_{j'\in[N]_{1,j}}\sigma_{1,j'}^2 \right) \\
			& + \dfrac{1}{N}\lim_{\epsilon\dto 0} \E{\vphantom{\sum_{j\in[N]_1}} (2N-1)\qbar_{1,j}^+ \ubar_{i,1} + (N-1)\left(\sum_{i'\in[N]_i}\qbar_{1,j}^+ \ubar_{i',1} +\sum_{j'\in[N]_1}\qbar_{1,j}^+ \ubar_{i,j'}  \right) } \\
			& +\dfrac{1}{N}\lim_{\epsilon\dto 0}\E{\vphantom{\sum_{j\in[N]_1}} (2N-1) \qbar_{i,1}^+\ubar_{1,j} + (N-1)\left(\sum_{i'\in[N]_1}\qbar_{i,1}^+\ubar_{i',j} + \sum_{j'\in[N]_j}\qbar_{i,1}^+\ubar_{1,j'} \right) } \\
			& -\dfrac{1}{N}\lim_{\epsilon\dto 0}\E{\sum_{(i',j')\in [N]_i\times [N]_1\setminus \{(1,j)\}}\qbar_{1,j}^+ \ubar_{i',j'} + \sum_{(i',j')\in[N]_1\times [N]_j\setminus \{(i,1)\}} \qbar_{i,1}^+\ubar_{i',j'} } \\
			& \pushright{\forall i,j\in[N]_1}
		\end{aligned}
	\end{align}
	where $\mathcal{P}$ is the set of $N\times N$ permutation matrices and
	\begin{align*}
		\mathcal{A}_1 =& \left\{(x,y)\in[N]\times [N]:\; y\geq x+1 \right\} \\
		\mathcal{A}_2 =& \left\{(x,y)\in[N]\times [N]:\;  y\geq x+1 \;,\; 2\leq x\leq N-1 \right\}.
	\end{align*}
\end{theorem}

Equation \eqref{s.eq.n.0} is one interpretation of SSC, which says that all the schedules have the same weight in the cone $\cK$. Observe that in Theorem \ref{s.theorem:2switch} we did not have an equation of the form of \eqref{s.eq.n.0}. However, we used this condition in the proof to obtain a system of equations (see Claim \ref{switch.claim.q11u22}). In this case, we decided to write it as an equation to make explicit the use of SSC.
	
	We provide the details of the proof of the theorems.

\section{Details of the proof of Theorem \ref{gs.thm:bounds}.}\label{app:gs.proofs.all}

We prove the lemmas and the claims that we use in the proof of Theorem \ref{gs.thm:bounds}.

\subsection{Proof of Lemma \ref{gs.lemma:bjl.bl}.}\label{app:lemma.gs.bjl.bl}

\proof{Proof of Lemma \ref{gs.lemma:bjl.bl}.}
First, recall $\vnu\in \cC$ and, by definition of the capacity region we have 
\begin{align*}
\cC=\sum_{m\in \cM} \psi_m\,ConvexHull\left(\cS^\pm \right).
\end{align*}
Then, since each $\cS^\pm$ is finite, for each $m\in \cM$ there exists $\vnu^\pm\in \cS^\pm$ such that
\begin{align*}
\vnu= \sum_{m\in \cM} \psi_m \vnu^\pm
\end{align*}

Also, by definition of the $\ell\tth$ hyperplane, for each $\ell\in P$ we have
\begin{align*}
\bl =& \max_{\vx\in \cC}\langle \vcl,\vx\rangle \\
=& \sum_{m\in \cM}\psi_m \max\left\{\langle \vcl,\vx\rangle: \vx\in ConvexHull\left(\cS^\pm\right)\right\} \\
\stackrel{(a)}{=}& \sum_{m\in \cM}\psi_m \max_{\vx\in\cS^\pm} \langle\vcl,\vx\rangle \\
\stackrel{(b)}{=}& \sum_{m\in\cM}\psi_m \bml
\end{align*}
where $(a)$ holds because the objective function of the maximization is linear and, therefore, the optimal solution is an extreme point of $ConvexHull\left(\cS^\pm\right)$, which must be an element of $\cS^\pm$ by definition of convex hull; and $(b)$ holds by definition of $\bml$. This proves that $\bl=\E{\Bbar_\ell}$.

Observe that the last equality also implies that $\langle \vcl,\vnu^\pm\rangle=\bml$ for all $m\in\cM$, for the following reason. First, by definition of $\bml$ we know $\langle \vcl,\vnu^\pm\rangle\leq \bml$. Also, if there exists $m^*\in \cM$ with $\langle \vcl,\vnu^{(m^*)}\rangle<b^{(m^*,\ell)}$, then
\begin{align*}
\sum_{m\in \cM} \psi_m \langle\vcl,\vnu^\pm\rangle < \sum_{m\in\cM} \psi_m \bml.
\end{align*}
But
\begin{align*}
\langle \vcl,\vnu\rangle = \sum_{m\in \cM} \psi_m \langle\vcl,\vnu^\pm\rangle\quad\text{and}\quad \bl=\sum_{m\in\cM} \psi_m \bml.
\end{align*}
Therefore, we got a contradiction because $\vnu\in\bigcap_{\ell\in P} \cF^\pl$ and, hence, $\langle\vcl,\vnu\rangle=\bl$.
\Halmos
\endproof

\subsection{Proof of Lemma \ref{gs.lemma:cl.s.bl}.}\label{app:gs.lemma.cl.s.bl}

\proof{Proof of Lemma \ref{gs.lemma:cl.s.bl}.}
For ease of exposition, we omit the dependence on $\epsilon$ of the variables. Define
\begin{align*}
\gammam\defn&\, \min\left\{\bml-\langle\vcl,\vx\rangle:\; \langle\vcl,\vx\rangle<\bml ,\text{ for } \ell\in P,\, \vx\in \cS^\pm\right\}.
\end{align*}
Observe that, for each $m\in\cM$ we have $\gammam>0$ because each $\cS^\pm$ is a finite set and, therefore, $\bml-\langle\vcl,\vx\rangle$ cannot be arbitrarily close to zero.

From stability, for each $\ell\in P$ we have
\begin{align*}
\E{\langle\vcl,\vsbar\rangle}\geq& \E{\langle \vcl,\vabar\rangle}.
\end{align*}

Then, using Lemma \ref{gs.lemma:bjl.bl}, we obtain that for each $m\in \cM$
\begin{align*}
\Em{\langle\vcl,\vsbar\rangle} &\geq \E{\langle\vcl,\vabar\rangle} 
= (1-\epsilon)\langle\vcl,\vnu^\pm\rangle 
= (1-\epsilon)\bml 
\end{align*}

On the other hand, by definition of $\piml$ we have
\begin{align*}
& \Em{\langle\vcl,\vsbar\rangle}= \piml\bml + \left(1-\piml\right)\,\Em{\left. \langle\vcl,\vsbar\rangle \right| \bml\neq \langle\vcl,\vsbar\rangle}
\end{align*}

Putting these two results together we obtain
\begin{align}
\begin{aligned}\label{gs.eq.lemma.cl.s.bl.one}
&(1-\epsilon) \bml
\leq \bml\piml
+ \left(1-\piml\right)\, \Em{\langle\vcl,\vsbar\rangle\,\left|\, \bml\neq \langle\vcl,\vsbar\rangle\right.}
\end{aligned}
\end{align}

Also, by definition of $\gammam$ we have
\begin{align*}
\Em{\langle\vcl,\vsbar\rangle\,\left|\, \bml\neq \langle\vcl,\vsbar \rangle \right.}\leq \bml-\gammam.
\end{align*}

Replacing the last result in \eqref{gs.eq.lemma.cl.s.bl.one} and rearranging terms we obtain
\begin{align*}
1-\piml\leq \dfrac{\epsilon\bml}{\gammam}.
\end{align*}

Since the sets $P$ and $\cS^\pm$ for each $m\in \cM$ are finite, there exists $b_{\max}=\max_{m\in\cM,\,\ell\in P}\{\bml \}$ and $b_{\max}<\infty$. Therefore, we have $1-\piml$ is $O(\epsilon)$.
\Halmos
\endproof

\subsection{Proof of Claim \ref{gs.claim.T1.partial}.}\label{app:proof.claim.T1}
\proof{Proof of Claim \ref{gs.claim.T1.partial}.}
We condition on the channel state. Then, we have
\begin{align*}
\E{\langle \vqbarparh, \vsbar-\vnu\rangle} =& \sum_{m\in\cM}\psi_m \Em{ \langle \vqbarparh, \vsbar-\vnu^\pm\rangle} \\
\stackrel{(*)}{=}& \sum_{m\in\cM}\psi_m \Em{ \langle \vqbarparh, \vsbar-\vnu^\pm\rangle \ind{\langle\vcl,\vsbar\rangle\neq \bml}},
\end{align*}
where $\vnu^\pm$ is defined as in Lemma \ref{gs.lemma:bjl.bl}. Equality $(*)$ holds because $\vqbarparh = \sum_{\ell\in\tilde{P}}\tilde{\xi}_\ell \vcl$ for $\tilde{\xi}_\ell\in\bR$ for all $\ell\in\tilde{P}$ (by definition of projection on the subspace $\cH$) and if the channel state is $m$ we have
\begin{align*}
\langle\vcl,\vsbar-\vnu^\pm\rangle \ind{\langle\vcl,\vsbar\rangle= \bml} = \left(\bml-\langle\vcl,\vnu^\pm\rangle\right)\ind{\langle\vcl,\vsbar\rangle= \bml}=0,
\end{align*}
where the last equality holds by definition of $\vnu^\pm$.

It suffices to show that $\Em{ \langle \vqbarparh, \vsbar-\vnu^\pm\rangle \ind{\langle\vcl,\vsbar\rangle\neq \bml}}$ is $O(\sqrt{\epsilon})$ because $\vpsi = \left(\psi_m\right)_{m\in\cM}$ is a probability mass function and, therefore, each $\psi_m$ is bounded.

First observe that $\vqbar=\vqbarparh+\vqbarperph = \vqbarpark+\vqbarperpk$. Then,
\begin{align*}
\langle \vqbarparh,\vsbar-\vnu^\pm \rangle = \langle \vqbarpark,\vsbar-\vnu^\pm \rangle + \langle \vqbarperpk-\vqbarperph ,\vsbar-\vnu^\pm\rangle.
\end{align*}

Therefore,
\begin{align}
& \Em{ \langle \vqbarparh, \vsbar-\vnu^\pm\rangle \ind{\langle\vcl,\vsbar\rangle\neq \bml}} \nonumber \\
=& \Em{\langle \vqbarpark,\vsbar-\vnu^\pm \rangle\ind{\langle\vcl,\vsbar\rangle\neq \bml}} + \Em{\langle \vqbarperpk-\vqbarperph ,\vsbar-\vnu^\pm\rangle\ind{\langle\vcl,\vsbar\rangle\neq \bml}} \label{gs.eq.claimT1.partial}
\end{align}

We prove that each term in \eqref{gs.eq.claimT1.partial} is $O(\sqrt{\epsilon})$. For the first term first observe
\begin{align*}
& \Em{\langle \vqbarpark,\vsbar-\vnu^\pm \rangle\ind{\langle\vcl,\vsbar\rangle\neq \bml}} \leq 0
\end{align*}
because $\vqbarpark = \sum_{\ell\in P}\xi_\ell\vcl$ with $\xi_\ell\geq 0$ for all $\ell\in P$ (since the projection is on the cone $\cK$) and
\begin{align*}
\langle\vcl,\vsbar-\vnu^\pm\rangle = \langle\vcl,\vsbar\rangle-\bml \leq 0,
\end{align*}
by definition of $\vnu^\pm$ and $\bml$. Then, we have
\begin{align}
0\geq& \Em{\langle \vqbarpark,\vsbar-\vnu^\pm \rangle\ind{\langle\vcl,\vsbar\rangle\neq \bml}} \nonumber \\
\stackrel{(a)}{=}& \Em{\langle\vqbar,\vsbar-\vnu^\pm\rangle \ind{\langle\vcl,\vsbar\rangle\neq \bml} } - \Em{\langle\vqbarperpk, \vsbar-\vnu^\pm\rangle \ind{\langle\vcl,\vsbar\rangle\neq \bml}} \nonumber \\
\stackrel{(b)}{\geq}& - \Em{\langle\vqbarperpk, \vsbar-\vnu^\pm\rangle \ind{\langle\vcl,\vsbar\rangle\neq \bml}} \label{eq:MaxWtUB}\\
\stackrel{(c)}{\geq}& - \sqrt{\E{\left\|\vqbarperpk \right\|^2} \Em{\left\| \vsbar-\vnu^\pm \right\|^2 \ind{\langle\vcl,\vsbar\rangle\neq \bml}}} \nonumber \\
\stackrel{(d)}{\geq}& -\sqrt{T_2\, \Em{\left\| \vsbar-\vnu^\pm \right\|^2 \ind{\langle\vcl,\vsbar\rangle\neq \bml}} } \nonumber
\end{align}
where $(a)$ holds because $\vqbarpark=\vqbar-\vqbarperpk$; $(b)$ holds by MaxWeight as described in \eqref{gs.eq.MW} and because $\vnu^\pm\in\cS^\pm$; $(c)$ holds by Cauchy-Schwarz inequality; and $(d)$ holds by SSC in Proposition \ref{gs.prop:gen.switch.SSC}.

The last step is to prove that the first term in \eqref{gs.eq.claimT1.partial} is $O(\sqrt{\epsilon})$. 
We have
\begin{align}
0\leq& \Em{\left\| \vsbar-\vnu^\pm \right\|^2 \ind{\langle\vcl,\vsbar\rangle\neq \bml}} \nonumber \\
\stackrel{(a)}{=}& \Em{\left. \left\| \vsbar-\vnu^\pm \right\|^2 \right| \langle\vcl,\vsbar\rangle\neq \bml}\left(1-\piml \right) \nonumber\\
\stackrel{(b)}{\leq}& n\left(\smax^2 + V_{\max}^2\right)\left(1-\piml \right) \nonumber\\
\stackrel{(c)}{=}& O(\epsilon) \label{eq.gs.claimT1.last},
\end{align}
where $V_{\max} = \max_{m\in \cM, i\in[n]}\nu^\pm_i$ is a finite constant. Equality $(a)$ holds by definition of $\piml$ in Lemma \ref{gs.lemma:cl.s.bl}; $(b)$ holds because $\vsbar$ and $\vnu^\pm$ have bounded elements; and  $(c)$ holds by Lemma \ref{gs.lemma:cl.s.bl}.

The proof that the second term in \eqref{gs.eq.claimT1.partial} is $O(\sqrt{\epsilon})$ holds by linearity of dot product, Cauchy-Schwarz inequality and \eqref{eq.gs.claimT1.last}. We omit it for brevity.
\Halmos
\endproof

\subsection{Proof of Claim \ref{gs.claim.Eu}.}\label{app.proof.gs.claim.Eu}

\proof{Proof of Claim \ref{gs.claim.Eu}.}
	We set to zero the drift of $V_l(\vq)=\sum_{\ell\in P}\langle\vcl,\vq\rangle$. 
We obtain
\begin{align*}
	0=& \E{\sum_{\ell\in P} \langle\vcl,\vqbar^+\rangle- \sum_{\ell\in P}\langle\vcl,\vqbar\rangle}
	= \E{\sum_{\ell\in P}\langle\vcl,\vabar-\vsbar+\vubar\rangle},
\end{align*}
where the last equality holds by definition of $\vqbar^+$ and by the dynamics of the queues presented in \eqref{gs.eq.dynamics.queues}. Rearranging terms we obtain
\begin{align*}
	\sum_{\ell\in P}\E{\langle\vcl,\vubar\rangle}
	=& \sum_{\ell\in P}\E{\langle\vcl,\vsbar\rangle} - \sum_{\ell\in P}\E{\langle\vcl,\vabar\rangle} \\
	=& \sum_{\ell\in P}\E{\langle\vcl,\vsbar\rangle} - \sum_{\ell\in P}\langle \vcl,(1-\epsilon)\vnu\rangle.
\end{align*}

But
\begin{align*}
	& \sum_{\ell\in P}\E{\langle\vcl,\vsbar\rangle} \\
	&= \E{\sum_{\ell\in P}\left(\langle\vcl,\vsbar\rangle-\bl \right)}+\sum_{\ell\in P}\bl \\
	&\stackrel{(a)}{=} \sum_{\ell\in P}\sum_{m\in \cM} \psi_m \Em{\left. \left( \langle\vcl,\vsbar\rangle-\bml \right)\right|  \langle\vcl,\vsbar\rangle\neq\bml}\left(1-\piml \right)  +\sum_{\ell\in P}\langle\vcl,\vnu\rangle \\
	&\stackrel{(b)}{=} \sum_{\ell\in P}\langle\vcl,\vnu\rangle-O(\epsilon),
\end{align*}
where $(a)$ holds because $\langle\vcl,\vnu\rangle=\bl$ for all $\ell\in P$.; and $(b)$ holds by \Cref{gs.lemma:cl.s.bl}. Then, since $\langle \vcl,\vnu\rangle=\bl$, we have
\begin{align*}
	\sum_{\ell\in P} \E{\langle\vcl,\vubar\rangle} = \epsilon \sum_{\ell\in P}\bl - O(\epsilon).
\end{align*}
\Halmos\endproof

\subsection{Proof of Equation \eqref{gs.eq.T4}.}\label{app:gs.T4}
\proof{Proof of Equation \eqref{gs.eq.T4}.}
We compute bounds on $\cT_4$. To do that we use the same approach that is used by \citeauthor{atilla} in \cite{atilla}. For each $\ell\in P$, let $\cL_+^\pl\defn \left\{i\in[n]:\; c_i^\pl>0 \right\}$ and define
\begin{align*}
\vcltilde=\left(c_i^\pl\right)_{i\in\cL_+^\pl} ,\, \vqbartildel=\left(\qbar_i\right)_{i\in\cL_+^\pl} \;\text{and}\; \vubartildel=\left(\ubar_i\right)_{i\in\cL_+^\pl}.
\end{align*}
Then,
\begin{align*}
0\leq \left|\dfrac{\cT_4}{2} \right|=& \left|\E{\langle\vqbarparh^+,\,\vubarparh\rangle} \right| \\
\stackrel{(a)}{=}& \left|\E{\langle\sum_{\ell\in P}\langle\vcl,\vqbar^+\rangle\vcl,\, \sum_{\ell\in P}\langle\vcl,\vubar\rangle\vcl \rangle} \right| \\
=& \left|\E{\langle\sum_{\ell\in P}\langle\vcltilde,\left(\vqbartildel\right)^+\rangle\vcltilde,\, \sum_{\ell\in P}\langle\vcltilde,\vubartildel\rangle\vcltilde \rangle} \right| \\
=& \left|\E{\langle\left(\vqbarparhtilde\right)^+,\,\vubarparhtilde\rangle} \right| \\
=& \left|\E{\langle \left(\vqbartildel\right)^+,\vubarparhtilde\rangle - \langle \left(\vqbarperphtilde\right)^+,\vubarparhtilde\rangle} \right| \\
=& \left|\E{\langle \left(\vqbartildel\right)^+,\vubartildel-\vubarperphtilde\rangle - \langle \left(\vqbarperphtilde\right)^+,\vubarparhtilde\rangle } \right| \\
\stackrel{(b)}{=}& \left|\E{\langle\left(\vqbarparhtilde\right)^+ +\left(\vqbarperphtilde\right)^+,\,-\vubarperphtilde \rangle- \langle \left(\vqbarperphtilde\right)^+,\vubarparhtilde\rangle  } \right| \\
\stackrel{(c)}{=}& \left|\E{-\langle\left(\vqbarparhtilde\right)^+,\,\vubarperphtilde\rangle - \langle \left(\vqbarperphtilde\right)^+,\,\vubartildel\rangle } \right| \\
\stackrel{(d)}{=}& \left|\E{- \langle \left(\vqbarperphtilde\right)^+,\,\vubartildel\rangle } \right| \\
\leq& \sqrt{\E{\left\|\left(\vqbarperphtilde\right)^+ \right\|^2}\E{\left\|\vubartildel\right\|^2}}
\end{align*}
where $(a)$ holds by definition of the projection on $\cH$; $(b)$ holds by equation \eqref{gs.eq.qu} and because $\left(\vqbartildel\right)^+=\left(\vqbarparhtilde\right)^+ +\left(\vqbarperphtilde\right)^+$; $(c)$ is obtained reorganizing terms; $(d)$ holds because $\left(\vqbarparhtilde\right)^+$ and $\vubarperphtilde$ belong to orthogonal subspaces; and $(f)$ holds by Cauchy-Schwarz inequality. Observe
\begin{align*}
\E{\left\|\left(\vqbarperphtilde\right)^+ \right\|^2}\leq \E{\left\|\vqbarperph^+\right\|^2}\leq T_2
\end{align*}
and
\begin{align*}
0\leq& \E{\left\|\vubartildel\right\|^2} \\
=& \sum_{i\in\cL_+^\pl}\E{\widetilde{u}_i^2} \\
\stackrel{(a)}{\leq}& \sum_{\ell\in P}\sum_{i\in \cL_+^\pl} \dfrac{\ctildel_i}{\ctildel_i} \E{\widetilde{u}_i^2} \\
\leq& \dfrac{\smax}{\widetilde{c}_{\min}} \sum_{\ell\in P}\sum_{i\in\cL_+^\pl} \ctildel_i\E{\widetilde{u}_i} \\
=& \dfrac{\smax}{\widetilde{c}_{\min}} \sum_{\ell\in P} \E{\langle\vcltilde,\,\vubartildel\rangle} \\
\stackrel{(b)}{=}& \dfrac{\smax}{\widetilde{c}_{\min}} \sum_{\ell\in P}\E{\langle\vcl,\vubar\rangle} \\
 \stackrel{(c)}{=}&  O(\epsilon)
\end{align*}
where $\ds\widetilde{c}_{\min}=\min_{\ell\in P, i\in[n]}\{\ctildel_i\}$ and $|P|$ is the cardinality of set $P$. Inequality $(a)$ holds because the terms in the summation are all non-negative; $(b)$ holds by definition of $\widetilde{\vc}^\pl$; and $(c)$ holds by \Cref{gs.claim.Eu}. This proves that $\cT_4=O(\sqrt{\epsilon})$. 
\Halmos
\endproof

\section{Details of the proof of Theorem \ref{s.theorem:2switch}.}\label{app:switch.proofs.all}

We prove the claims we made in the proof of Theorem \ref{s.theorem:2switch}. 

\subsection{Proof of Claim \ref{switch.claim.qperp.upar}.}\label{app:switch.proofs.qperp.upar}

\proof{Proof of Claim \ref{switch.claim.qperp.upar}.}
Observe
\begin{align*}
\E{\qbar_{\perp 1}^+\ubarpari{1}}& \leq  \E{\left|\qbar_{\perp 1}^+\right| \left|\ubarpari{1}\right|} \\ & \leq \E{\sum_{i=1}^4 \left|\qbar_{\perp i}^+\right|\left|\ubarpari{i}\right|} \\
& \stackrel{(a)}{\leq} \sqrt{\E{\left\|\vqbarperp \right\|^2} \E{\left\|\vubar_{\parallel} \right\|^2}} \\
& \stackrel{(b)}{\leq} \sqrt{T_2}\sqrt{\E{\left\|\vubar_{\parallel} \right\|^2}},
\end{align*}
where $(a)$ holds by Cauchy-Schwarz inequality; and $(b)$ holds by Proposition \ref{gs.prop:gen.switch.SSC}. Similarly,
\begin{align*}
\E{\qbar_{\perp 1}^+\ubarpari{1}}\geq& - \E{\left|\qbar_{\perp 1}^+\right| \left|\ubarpari{1}\right|}
\geq -\sqrt{T_2}\sqrt{\E{\left\|\vubar_{\parallel} \right\|^2}}.
\end{align*}
Then,
\begin{align*}
\left|\E{\qbar_{\perp 1}^+\ubarpari{1}} \right|\leq \sqrt{T_2}\sqrt{\E{\left\|\vubar_{\parallel} \right\|^2}} 
\end{align*}
and $\E{\left\|\vubar_{\parallel} \right\|^2}$ is $O(\epsilon)$. This proves the claim.
\Halmos
\endproof

\subsection{Proof of Claim \ref{switch.claim.q11u22}.}\label{app:switch.proofs.q11u22}

\proof{Proof of Claim \ref{switch.claim.q11u22}.}
We use Claim \ref{switch.claim.qperp.upar}. We obtain
\begin{align*}
\E{\qbar_{1}^+\ubar_4}=& \E{\qbarpari{1}^+\ubar_4} + O(\sqrt{\epsilon}) \\
\stackrel{(a)}{=}& \E{\left(\qbarpari{2}^+ +\qbarpari{3}^+ - \qbarpari{4}^+\right)\ubar_4} + O(\sqrt{\epsilon}) \\
\stackrel{(b)}{=}& \E{\left(\qbar_{2}^+ +\qbar_{3}^+ - \qbar_{4}^+\right)\ubar_{4}} + O(\sqrt{\epsilon}) \\
\stackrel{(c)}{=}& \E{\left(\qbar_{2}^+ +\qbar_{3}^+  \right)\ubar_{4}} + O(\sqrt{\epsilon})
\end{align*}
where $(a)$ holds by SSC; $(b)$ holds by Claim \ref{switch.claim.qperp.upar}; and $(c)$ holds by \eqref{gs.eq.qu}.
\Halmos
\endproof

	%

\end{APPENDICES}


\end{document}